\input amstex 
\scrollmode\NoBlackBoxes
\magnification=1100

\comment
\font\rm=Times at 10pt

\font\bf=TimesB

\font\it=TimesI at 10pt
\font\sc=Times at 7pt

\def\Sc #1{{\sc \uppercase{#1}}}
\endcomment

\long\def\Pf{\par\noindent {\it Proof.} }
\def\({\left(}
\def\){\right)}
\def\st{such that }
\def\qed{\hfill$\bullet$\vskip 4pt}

\def\brcs#1{\left\{ #1\right\}}
\def\det{\text{det}}

\def\Set#1#2{\brcs{#1 \left|\vphantom{#1 #2} \right.#2}}

\def\C{\text{\bf C}}

\def\Tr{\text{Tr}}

\def\wrt{with respect to }
\long\def\Lemma #1. #2\par{\noindent {\Sc  {#1.}} {
#2}\vskip 2pt}
\def\Arrow #1;#2.{#1\:#2 \to }
\def\spec{\text{spec}\,}


\def\P{{\Cal P}}
\def\R{\text{\bf R}}
\def\N{\text{\bf N}}
\def\Z{\text{\bf Z}}
\def\Q{\text{\bf Q}}
\def\Aff{\text{Aff}\,}
\def\Ag{{\Cal AG}}
\def\op{{}^{\text{op}}}
 
\def\slfrac#1#2{{\raise -.07 ex\hbox{$^{#1}$}}\!/\raise .35 ex \hbox{${}_{#2}$}}
\def\ssf #1/#2{\slfrac {#1}{#2}}

\long\def\Rmk{\par \noindent{\it Remark. }}
\def\diag{\text{diag}\,}

\def\T{{\Cal T}}
\def\I{\text{I}}
\def\TT #1,#2.{\T^c_{#1} \cap \T^c_{#2}}

\let\Gl= \gl 
\let\GL = \gl


   \long\def\Lem
#1.#2\par{\write1{#1,
p\,\folio\par}\vskip4pt{\baselineskip=13pt
   \noindent {\rm \uppercase{#1}} #2\vskip3pt

   }} 

   \long\def\Title #1\par {\noindent{ #1}\vskip 9pt}
   \long\def\SubT #1.{\noindent {\it #1\/} }
   \long\def\SecT #1\par{\vskip 4pt \noindent {\bf #1}\vglue1pt
   \noindent}
\def\Arg{\text{Arg}\,}
\def \Mn #1{\text{M}_{#1}}
\def\C{\text{\bf C}}
\def\I{\text{\bf I}}
\def\G{{\Cal G}}

\def\GL #1{\text{GL}(#1)}

\def\J{{\Cal J}}
\def\M #1,#2{{\Cal M}_{#1,#2}}
\def\Mm{{\Cal M}}
\def\H{{\Cal H}}
\def\Nn{{\Cal N}}
\def\T{{\Cal T}}
\def\Tt{\text{\bf T}}
\def\rk{\text{rk\,}}
\def\Op{^{\text{op}}}
\def\S{{\Cal S}}

\def\tQsr{{\text{tQsr}}}
\def\Ss{xx}
\let\iso=\cong


\def\oneone{1.4}
\def\onetwo{1.5}
\def\onethr{1.1}
\def\onefou{1.2}
\def\onefiv{1.3}

\def\twoone{2.1}

\def\throne{3.1}
\def\thrtwo{3.3}
\def\thrthr{3.5}
\def\thrfou{3.7}
\def\thrfiv{3.2}
\def\thrsix{3.4}
\def\thrsev{3.6}
\def\threig{3.8}
\def\thrnin{3.9}

\def\fouone{4.1}
\def\foutwo{4.2}
\def\fouthr{4.3}

\def\fivone{5.1}
\def\fivtwo{5.2}
\def\fivthr{5.3}
\def\fivfou{5.4}
\def\fivfiv{5.1}
\def\fivsix{5.2}
\def\fivsev{5.3}

\def\sixone{6.1}
\def\sixtwo{5.4}

\def\sevone{6.1}
\def\sevtwo{6.2}
\def\sevthr{6.3}
\def\sevfou{6.5}
\def\sevfiv{6.6}
\def\sevsix{6.4}
\def\sevsev{6.2}

\def\sevele{6.7}

\def\eigone{7.1}
\def\eigtwo{7.2}

\def\ninone{9.1}
\def\nintwo{9.2}

\def\Aone{A.1}
\def\Atwo{A.2}
\def\Athr{A.3}
\def\Afou{A.4}
\def\Afiv{A.5}

\def\Asev{A.7}
\def\Aeig{A.6}
\def\Anin{A.7}

\def\Atwe{A.9}
\def\Athi{A.10}
\def\Aftn{A.11}
\def\Affn{A.12}
\def\Asxn{A.13}
\def\Asvn{A.14}
\def\Aegn{A.15}
\def\Antn{A.16}
\def\Atty{A.17}
\def\Atwn{A.18}
\def\Atwt{A.23}
\def\Atwh{A.17}

\def\Atff{A.8}

\def\Asix{B.1}
\def\Aele{B.3}
\def\Aten{B.2}

\def\Bsev{B.4}
\def\Beig{B.6}
\def\Bnin{B.7}
\def\Bten{B.5}

\def\Bone{C.1}
\def\Btwo{C.2}
\def\Bthr{C.3}
\def\Bfou{C.4}

\def\Ctwo{D.2}
\def\Cthr{D.3}
\def\Cfou{D.4}
\def\Cfiv{D.5}

\def\Atwf{D.1}

\def\Ord{\text{\rm Ord}\,}
\def\ord{\text{\rm ord}\,}
\def\det{\text{det}\,}
\def\flt{{FLT}}
\def\Mn{\text{M}}
\def\Nn{\Cal N}
\def\tr{\text{tr\,}}
\def\diag{\text{diag\,}}
\def\PE #1.{$\text{PE}(2,#1)$}
\def\pe{$\text{PE}(2,R)$}
\def\petwo{$\text{PE}_2(2,R)$}%
\def\gl{\text{GL}}
\let\Gl= \gl 
\let\GL = \gl 

\def\T{{\Cal T}}
\def\I{\text{I}}
\def\TT #1,#2.{\T^c_{#1} \cap \T^c_{#2}}

\def\II{{\Cal I}}

 \def\lgui
 {{\rm(\!(}}
  \def\rgui
  {{\rm)\!)}}
  \def \gui #1{\lgui $#1$\rgui}
  \def\P{\Cal P}
\def\m #1,#2{m_{#1,#2}}
\def\aac{\'a}
\long\def\bu{\noindent $\bullet$\hskip 1 em}
   \def\pe{\text{PE}(2,R)}
   \def\peone{\text{PE}_1(2,R)}
   \def\petwo{\text{PE}_2(2,R)}

\Title {\bf Abstract fractional linear  transformations}%

{\rightskip 1 true in \leftskip=.5 true in
\parindent= 1 em
\baselineskip=10pt 

\noindent {\it Abstract.} 
We begin with (densely-defined) fractional linear transformations (\flt) on (some) Banach algebras and their relatives. This leads to Wedderburn's continued fractions (recursively-defined noncommutative polynomials) for any ring. Along the way, we discover a one-parameter family of (noncommutative) polynomials \st if one of them is invertible, then read in the opposite order, the corresponding polynomial is also invertible (extending the well known $1+ab$ is invertible if $1+ ba$ is, and the not-so-well-known, $a + abc + c$ and $a + cba + c$).

This in turn leads to a  definition of \flt\ for general rings $R$, which turns out to be $\pe$ (the projective elementary group). Using Wedderburn's polynomials, this permits us to define a length function on $\pe$, which suggests a stable range type condition (for $n  =1$, it {\it is\/} stable range one, but higher values do not correspond. 

Again using the length results, we prove the expected results for $\pe$: under very modest conditions on $R$, the commutator subgroup of $\pe$ is perfect and of index one or two. 

Along the same lines, we also prove results on simplicity of the commutator subgroup: we require the usual generative properties on the simple ring $R$, as well either the very restrictive $1$ in the range, or a mild condition about invertibles, involving intersections of three translates of $\gl(1,R)$. This last property is  explored in the appendices, which give examples (and non-examples). 

Numerous questions suggest themselves throughout. 
\comment 
For suitable subrings of (noncommutative) Banach algebras, we construct 
a concrete fractional linear transformation group. Based on the relations among the generators,
 we extend the definition to general rings (and correspondingly, the projective elementary group). Via a notion of length (roughly counting the minimal number of inverses that appear in a factorization of a group element), we also establishj connections between stable range-like conditions and Wedderburn's noncommutative  polynomials. Based on the latter we uncover a one-parameter sequence of noncommutative polynomials with properties extending the well known $1+ab$ is invertible iff $1+ba$ is invertible.  We also examine perfectness and simplicity in the groups closely connected to these constructions. This in turn leads to a study (in the appendix) of rings for which given elements $b,c$, there exist an invertible $u$ \st $u+b$ and $u+c$ are both invertible, among other things. 
\endcomment

}

\vskip 16pt

\noindent {\it David Handelman}\plainfootnote{${}$}{Keywords and phrases: fractional linear transformations, noncommutative continued fractions, stable range conditions, Hermite ring, elementary divisor ring.}\vskip 4pt

\noindent {}\plainfootnote{}{AMS (MOS) 2020 classification: 19C99, 20H99}

\vskip -8 pt

\SecT 0 Introduction

\noindent Motivated in part by two-sided fractional linear transformations (\flt) on rings of matrices [H] for suitable Banach algebras and subrings thereof, $R$ we construct a group of densely-defined transformations acting as \flt. This group can be identified with \pe{}, the projective group generated by size two elementary matrices. 

For an $n \times n$ matrix (with entries from $\C$), $X$, and six matrices $A,B,C,D,E,F$ with $A,B,D,E$ invertible, we can define the  transformation, $X \mapsto (AXB + C)(DXE+F)^{-1}$ on a dense open subset of $\Mn_n \C$. These can be composed to make even more complicated strings of partially defined transformations; however, it turns out that every such expression can be simplified to involve only two {${}^{-1}$} terms. For Banach algebras with stable range one, the same result holds, and the {\flt} forms a group (which it turns out is just \pe{}).

For general rings $R$, there is no apparent realization of \pe\  as a group of fractional linear transformations (although we can come close if $R$ satisfies a condition studied in the Appendix), so we consider {\pe} itself. And Weddurburn's noncommutative fractions [W] plays a crucial role. These are defined recursively, $P_n (a_1, a_2, \dots, a_n)$, and $Q_n (a_1, a_2, \dots, a_n)$ as polynomials in noncommuting variables. En route to establishing the connection between these and \pe{}, we come across a sequence, $1 + ab, a + abc + c, \dots$ with the property if one of them is invertible, then the one obtained by reversing the order of the variables is also invertible. Thus $1+ab$ is invertible if $1+ ba$ is (well known!), $a + abc + c$ is invertible if $a + cba + c$ is invertible (appears in one reference in the literature, as far as I can tell, [MV]), and all the subsequent ones are apparently new). 

When we view {\pe{}} via the natural set of generators suggested by the {\flt} point of view, we see that is a natural length function on elements, called $\ord (g)$. Then $R$ has stable range one if and only $\ord (g) \leq 2\slfrac12$. Larger bounds yield a stable range-like condition, but it is not stable range. 

Also using the $\ord$ function, we give a sufficient condition for {\pe{}} to be perfect or simple. A sample result: if $R$ is simple, has one in the stable range, has a centre with at least four elements, and is generated as a ring by the centre and the  commutator group of the invertibles, then {\pe{}} is simple. The relatively restrictive  stable range one condition can be replaced by the condition mentioned above on invertibles, motivating its study in the appendices. 
Much simpler is a condition guaranteeing that {\pe{}} is perfect.

The Appendices contain results (mostly elementary)  of the type, if $a(1), a(2), \dots, a(n)$ belong to $R$, then the intersection of their translates of $\gl(1,R)$ is nonempty; that is, if $\gl (1,R) - a $ denotes the set $\Set{u-a}{u \in \gl(1,R)}$, then 
$$
\bigcap_{j=1}^n (\gl (1,R) -a(j)) \neq \emptyset \tag {\gui n}
$$

This leads to a number of interesting questions, expecially when $n =3$ (which is what is needed, absent stable range one, in the group results above), some of which we can answer. 

Section 1 deals with the construction of {\flt} on suitable Banach algebras, and motivates the introduction of Wedderburn's $p$s and $q$s (which we have relabelled with capital letters). Sections 2 and 3 analyze them over general rings, and we obtain the invertibility result cited above (Corollary \thrthr), and a similar but more restrictive one for zero divisors (Corollary \thrtwo). We then establish the isomorphism between the group of {\flt} (when defined) and {\pe}, (Proposition \thrnin)

Section 4 introduces a notion of length on elements of {\pe{}} Sections 5 and 6 deal with normal subgroups via the length function. For example, under modest conditions, the commuator subgroup of {\pe{}} is identified, and shown to be perfect (Proposition \fivfou). One is in the stable range of $R$ if and only if for every $g$ in {\pe{}}, $\ord (g) \leq 2\slfrac12$ (Proposition \sevsix). The simplicity result cited above also appears in this section (Proposition \sevele).

\comment
Section 7 deals with the abelianization of {\pe{}} [
]. Finally, Section 8 provides direct proofs in the case of Banach algebras with stable range one, of some of the earlier results. [
]. 
\endcomment

The Appendices deals with classes of rings satisfying the intersection property \gui n, described above. This was motivated by results in Section 1 about constructing {\flt}, and more importantly in Section 6, dealing with simplicity. Among other results, if $F_q$ is a finite field with $q$ elements, then for every $n >1$, $\Mn_n F_q$ satisfies \gui {q} (Proposition \Aele). This naturally leads to number theoretic questions; for example, how many primes do we have to invert so that $\Z[1/2p_1p_2\cdots p_k]$ satisfies the very simple condition \gui 2?). Does $\Mn_3 \Z$ satisfy \gui 3?

Appendix C is devoted to proving that for all $n \geq 3$, $\Mn _n F_2$ satisfies \gui 3. Appendix D considers a weakening of \gui {2} (called \gui{1}) and shows that for $F$ a field that is not algebraic over a finite field, $F[x,P^{-1}]$ fails \gui{1} for all polynomials $P$ in the polynomial ring.

In this paper, {\it invertible\/} (element of a ring) means {\it two-sided invertible}. Rings are not necessarily commutative. 
 The group of invertible elements of the $n \times n$ matrix ring over $R$ is denoted $\gl(n,R)$, and this applies when $n =1$, that is, $\gl(1,R)$; the last is often denoted $\gl(R)$. 
 
 We 
denote the commutator subgroup of a group $H$ by $H'$ or $[H,H]$; the abelianization of $H$, $H/[H,H]$,
 is denoted $H^{\text{ab}}$.
 
 Many elementary techniques are similar to those of lower K-theory. 
 
 \comment
Let by noncommutative fractional linear
transformations [H], but turned out to be an elaboration of a
remark of Kolster [K]. The latter concerns  recursively-defined
noncommuting polynomials (with coefficients in an arbitrary ring). They
have properties resembling  those of continued fractions; one sequence of
them begins  $Q_0 = 1$, $Q_1(a_1) = a_1$, $Q_2 (a_1,a_2) = 1 + a_2 a_1$,
$Q_3 (a_1,a_2, a_3) = a_3 + a_1 + a_3 a_2 a_1$, \dots. They have the
interesting property that $Q_k(a_1,a_2,\dots,a_k)$ is invertible in
the ring (for a sequence   of elements of the ring)
if and only if $Q_k(a_k,a_{k-1},a_{k-2},\dots,a_1)$ is invertible
(generalizing a well-known observation about $Q_2(a_1,a_2) = 1+ a_2 a_1$
versus $1+ a_2 a_1$). 

The idea is first to obtain a concrete model of a group of transformations
that includes those of the form (densely defined) $X \mapsto (AXB + C)(DXE +F)^{-1}$ (for $X,A,B,C,D,E,F$ 
in the ring); this only applies to dense subalgebras of Banach algebras for
which the invertibles are relatively dense. However, observing the relations 
among the generators, we can construct an abstract form via generators and relations
that applies to any (unital) ring. 

The noncommuting polynomials enable us to attach three groups, $\Ag
(R)$, $\Ag_1 (R)$, and $\Ag_2 (R)$, to an arbitrary unital ring $R$. The
assignment of $R$ to any one of these is functorial only for ring
homomorphisms which map the centre to the centre. When $R$ has $1$ in its
stable range,  $\Ag (R)/\Ag_2 (R)$ is easy to describe, but other
stable range conditions seem to be irrelevant---suggesting that
$1$-stable range can be generalized in a different manner from the usual
one. Generically, $\Ag_2 (R)$ is perfect, and simplicity can be obtained 
in some cases. 

Many of the techniques are similar to those of lower $K$-theory. 
\endcomment 

\SecT 1 Concrete construction

Under suitable hypotheses on a ring $R$, we can construct a group of partially defined maps from subsets of $R$ to $R$.

Define the following three classes of (partial) maps.
$$\eqalign{
\Arrow T_s;R.R \qquad \qquad & T_s(r)= r+s \quad \text{for each
$s$ in $R$}\cr
\Arrow \J;\GL R.\GL R  \qquad \qquad & \J(r)= r^{-1} \cr
\Arrow \M x,y;R.R  \qquad \qquad& \M x,y (r) = xry \quad \text{for each
pair $(x,y)$ in
$\GL R \times \GL R\op$} \cr }$$

For invertible $r,s$ in $R$,  and for arbitrary $a$ in $R$, $\M x,y$ and $T_a$ are invertible functions from $R$ to $R$. They  generate a subgroup $\H$, which is a semidirect product. We may form formal strings, $h_1 \J h_2 \J \dots \J h_k$, where $h_i$ belong to $\H$, and we can assume that if $h_i$ is the identity transformation, then either $i = 1$ or $k$. We will see that these strings are defined on intersections of translates of $\gl (1,R)$. Unfortunately, in general, these could be empty, or even if not empty, there is another problem which we will discuss.

For this section, we make the following assumption, for all $n$, $R$ satisfies, for all sequences  $\brcs{d_i}_{i=1}^{\infty}$ inside $R$
$$
\cap_{i =1}^n (\gl (1,R) - d_i) \neq \emptyset  \tag {\gui n} 
$$
For some rings, even the countable intersections are always nonempty, and this includes Banach algebras (see Appendices A--C where these properties are discussed in detail.)

Our immediate aim is to show that if $R$ satisfies \gui {n} for all $n \geq 2$, then we can form a group out of the fractional linear transformations defined above. 

\Lem Lemma \onethr. For $h \in \H$ and $d$ in $R$, there exists $b$ in $R$ \st $h^{-1}((\gl (1,R) - d) $ contains $\gl(1,R) - b$.

\Rmk The ambiguous use of $h^{-1}$ either as a set function or a transformation is deliberate, as $h $ and $h^{-1`}$ are defined on all of $R$. 

\Pf Write $h = V_1 V_2 \cdots V_m$ where each $V_i$ is either of the form $T_a$ or $\M r,s$ for appropriately varying choices of elements. We examine $V_m^{-1} V_{m-1}^{-1}\cdots V_1^{-1} (\gl (1,R) - d)$, and proceed by induction on $m$. 

If $V_1 = T_a$, then $V_1^{-1} (\gl (1,R) - d) = T^{-1}(\gl (1,R) - d)$, and this is $\gl(1,R) - (a+d)$. If $V_1 = \M r,s$ (with invertible $r,s$), then $V_1^{-1} (\gl (1,R) - d) = (\gl (1,R) - r^{-1}ds^{-1})$. In either event, the inverse image is a single translate of $\gl (1,R)$, and induction applies. \qed

Let $a(1), a(2), \dots$ be elements of $R$. Define the (possibly empty) subsets inductively via 
$$\eqalign{
U_0 & = \gl(1,R) \cr
U_i (a(1), a(2), \dots, a(i) )& = \(\gl (1,R) \cap (U_{i-1}) (a(1), a(2), \dots, a(i-1))  -  a(i))\)^{-1}
}$$
The inverse at the right means the set of inverses of the elements in the set. If we have a formal string $w = \J h_n \J h_{n-1} \J \dots h_2 \J h_1$ (with $h_i$ in $\H$), then $w$ is defined on $U_n (a(1), \dots, a(n))$ for suitable ring elements $a(i)$. The definition of the $U_i$ is motivated by the fact that 
$$
h_1^{-1} \J \dots \J h_n^{-1} (\gl (1,R))
$$ 
is of the form $U_i$  for a suitable set of elements. 

In particular, we can compose two strings by taking the intersection of their domains, and, assuming \gui {n}  holds for all  $n$, the intersections will be nonempty. This will give a group structure to the compositions, but there is still a problem. Suppose two such strings are defined on the intersection of their domains, and they agree there.  Now suppose each is defined on a possibly different set---do they agree on it?
If so, then we can define the group of fractional linear transformations in the obvious way, as the set of equivalence classes of ordered pairs $(w,U)$ where $U$ contains one of the $U_i$, and equivalence is equality on the intersection. 
This permits the product operation, stringing together, forming a group, which we denote $\G(R)$. It is the free product $\H \circ \Z_2$ modulo some relations. 

So we require a condition such as
\item{$\bullet$} If two strings are defined on some $U_n{a(1), \dots, a(n)}$ and they agree on it, then they agree on any set on which they are both defined.

This is in general horrible to verify, but there is one situation in which it is straightfoward. Suppose that $R$ is a Banach algebra; then \gui {n} holds for all $n$. In particular, any finite intersection of $U$s is not empty.  

A Banach algebra $R$ has one in the stable range if $\gl (1,R)$ is dense in $R$ (this is not the definition, it is a theorem). Since $\gl(1,R)$ is always open in Banach algebras, any finite intersection of translates of $\gl(1,R)$ is thus open, and if $1$ is in the stable range, dense. Since all of the transformations are continuous on their domains, they are uniquely determined by their effect on any of them. 
Hence $(\bullet)$ holds (and a bit more: any countable intersection is still dense). 

Marginally more generally, here is a set of conditions on a topological ring $R$, which together result in $(\bullet)$ holding. 

\item{(A)} $\gl(1,R)$ is open in $R$;
\item{(B)} the map $x \mapsto x^{
-1}$ is continuous on $\gl(1,R)$
\item{(C)} $\gl (1,R)$ is dense in $R$
\item{(C$^{\prime}$)} for every $d$ in $R$, the set $\gl(1,R) \cap (\gl(1,R) - d)$ is dense in $\gl(1,R)$. 

Obviously, (C) implies (C${}^{\prime}$) in the presence of (A). At first glance, (C$^{\prime}$) represents a new condition for Banach algebras. But it doesn't. 

\Lem Lemma \onefou.Let $R$ be a metrizable topological ring. If $R$ satisfies (C$^{\prime}$), then it satisfies (C).

\Pf Select $c$ in $R$ and write $c = 1 + (c-1)$, and set $b = c-1$. There exist $f_n $ in $(\gl (1,R) - b) \cap \gl (1,R)$ \st $f_n \to 1$. Set $h_n = f_n + b$; these belong to $\gl(1,R)$ and converge to $1 + b = c$. \qed

(All that is needed is that every limit point of a subset is a limit point of sequence in that subset; this is weaker than metrizability.)

Banach algebras satisfy (A) and (B); they satisfy (C) if and only if one is in their stable range, and this is a large class.

If the topological ring $R$ satisfies (A--C), then each of the open sets $U_n (\cdot)$ is dense, so that the intersection of any finite collection of them is also dense (and open). Continuity of the transformations on their respective domains yields ($\bullet$). 

\Lem Corollary \onefiv.  Let $S$ be a Banach algebra with $1$ in its stable range, and let $R$ be a dense subring of $S$ \st $\gl(1,R) = R \cap \gl(1,S)$. Then $\G (R)$ exists.

Dense subalgebras of C*-algebras, with the indicated property appear regularly.  A different type of example (with {\it subring\/} rather than {\it subalgebra\/}) is the following. Let $T$ be a field of characteristic zero \st   $|T| \leq |\C|$. Then $T$ is embeddable in $\C$, and thus its closure (\wrt the inherited topology) is either $\R$ or $\C$. If we take $R = \Mn_n T \subset \Mn_n \C = S$, then its closure (inside $\Mn_n \C$) is $\Mn_n \R$ or $\Mn_n \C$. It is trivial to verify that $\gl (1,R) = R \cap \gl(1,S)$ holds. So the corollary applies in these cases. 

We obtain another class with direct limits of finite-dimensional (not necessarily semisimple) subalgebras over a subfield of $\C$; there is a norm on it for which the completion, $S$, is a Banach algebra, and since the direct limit is algebraic over some field, it is immediate that the condition on invertibles is satisfied (and $S$ has one in the stable range). 

Presumably, there are other examples of topological rings satisfying (A,B,C) unrelated to Banach algebras, but I haven't found any. For example, adic topologies usually have $\gl (1,R)$ closed in $R$ (so that (C) fails); regular rings with a rank function metric will similarly fail (C). Rings of ad\`eles fail (B), and the weak or strong topology on von Neumann algebras fail (A) (of course, if we treat the last class as C*-algebras, the we would use the norm topology). 

\comment
Let $\Mm$ denote the subgroup generated by all the $\M x,y$; it is easy
to see that $\Mm $ may be identified with $\GL R \times \GL
R^{\text{op}}/\brcs{z,z^{-1}}$, where $z$ runs over the centre of $\GL
R$. Let $\H$ denote the group generated by $\Mm $ and $\J$; because $\M
x,y \J = \J \M y^{-1},{x^{-1}}$, $\H$ is the semidirect product of $\Mm$
by the order two element $\J$. Let $\T$ denote the group generated by the
translations, $T_s$; so $\T$ can be identified with the additive group of
$R$. From the relation $\M x,y T_s = T_{xsy} \M x,y$, it follows that the
group generated by $\T$ and $\Mm$ is a semidirect product of $\T$ by
$\Mm$.

\endcomment

We notice that if $F,D,E$ are elements of $R$, then $\J T_F \M D,E $ is the (densely-defined) 
transformation $X \mapsto (DXE + F)^{-1}$ (essentially the 
denominator of the  general fractional linear transformation described above); 
when $R = \Mn \C$, this was studied in detail in [H], particularly \wrt its fixed points. Everything
of the form $X \mapsto (AXB + C)(DXE +F)^{-1}$ (for $X,A,B,C,D,E,F$ 
in the ring) can be obtained from a word the generators requiring at most two $\J$s. 

\comment
The equation
$$\J T_1 \J T_2 \dots \J T_m = \M x,y
$$ can be solved explicitly in terms of $x,y, T_i$.

Every element of $\G$ can be written in at least one of the following
forms
$$\eqalign{
& T_s \quad \M x,y \quad \J \quad T_s \M x,y \quad \J \M x,y
\quad T_s
\J \M x,y
\quad \J T_s \J \M x,y \cr 
& T_s \J T_r \M x,y \quad
\quad \J T_s \J T_r  \M x,y \quad \J T_s \J T_r \M x,y  \quad
\quad T_s \J T_r \J T_u\M x,y \cr
}$$ for some choice of the parameters $r,s,u, x,y$ (with $x$ and $y$
invertible).

Define the obvious equivalence relation on elements of $\G$. For $g$ and
$g'$ in $\G$, we declare $g \sim g'$ (relative to $\G$) if there exists
$f$ in $\G$ \st $fgf^{-1} = g'$.

We show that if $R= \Mn n F$ for some subfield $F$ of $\C$, then  $T_a
\sim T_b$ if and only if $\rk a = \rk b$; alternatively, all the
equivalences are implemented by elements of $\H$. On the other hand, the
conjugacy problem for
$\M r,s$ is somewhat more complicated---in general, not all the
equivalences can be implemented by elements of $\H$, but the exceptional
set (consisting of elements which have the three T's in the display) is
weird.

\endcomment

We observe some simple relations among the elements of $\G$. When we
write $\M x,y$, it is implicit that $x$ and $y$ are elements of $\GL R$.

$$\eqalign{  T_a T_b & = T_{a+b} \qquad \M u,v \M w,x = \M {uw},{xv} \cr
\M x,y \J & = \J \M y^{-1},{x^{-1}} \qquad \M x,y T_a = T_{xay} \M x,y.
\cr  }$$

The relations between $\J$ and elements of $\T$ are far more complicated (or interesting, depending on your point of view).
To begin  analyzing them, we require an elementary lemma.

\Lem Lemma \oneone. Let $R$ be a topological ring satisfying conditions (A--C), and suppose there exists a central element $k$ \st $k(k-1)$ is not a zero divisor. Suppose that $a,b,c,d,r,s$ are elements of $R$ with the latter two invertible, and there exists a dense open subset $U$ of $R$ \st for all $x$ in $U$, $cx+ d$ is invertible, and 
$$
(ax + b) (cx + d)^{-1} = rxs. \tag 1
$$ 
Then $b = c = 0$ and there exists central  invertible $\lambda$ \st $a = \lambda r$ and $b = \lambda s^{-1}$. 

\Pf We have 
$$
(ax + b) = rxs (cx + d) \tag 2
$$
for all $x$ in $U$. Since $U$ is dense and ring elements act continuously, we have that (2) holds for all $x$ in $R$. 

Setting $x = 0$, we deduce $b = 0$. Setting $x = 1$ and then $x =k$, we deduce $k(k-1) rsc = 0$. By hypothesis, $rsc = 0$, and since $r,s$ are invertible, we deduce $c = 0$. 

Thus $ax = rxsd$ for all $x$. Since $cx + d$ is invertible for at least one value of $x$, and $c = 0$, it follows that $d $ is invertible. Setting $x = 1$, we infer $a = rsd$, so that $a$ is invertible, and $r^{-1} a= sd$. Then $r^{-1}ax = xsd$ for all $x$, and so $r^{-1}a = sd:= \lambda $ central and invertible. 
\qed

The condition that such a $k$ exist is extremely weak.

\comment
\Lem Lemma \oneone. Suppose the centre of $R$ contains at least two
invertible elements, and $R$ is generated additively by its invertible
elements. Suppose that $a,b,c,d$ are elements of $R$ and $r,s$ are in
$\GL R$ \st for all $x$ in a dense subset of $R$, $cx +d$ is invertible
and $(ax+b)(cx+d)^{-1} = rxs$. Then $b = c = 0$ and there exists a
central invertible element $\lambda$ \st $a = \lambda r$ and $d = \lambda
s^{-1}$.

\Pf Post-multiplying by $cx +d$, we deduce that $ax+b = rxs(cx+d)$ for a
dense set of $x$s; by continuity of multiplication, this equality holds
for all $x$ in $R$. Substituting $x= 0$ yields $b = 0$. Let $\gamma$ be
any element of the centre. Substituting $x = \gamma$, we deduce $\gamma a
= \gamma^2 rsc + \gamma rsd$. If $\gamma $ is not a zero divisor, we may
divide by it, and thus deduce $a - rsd = \gamma rsc$. As more than one
such $\gamma$ exists, we deduce that $rsc = 0$; as $rs$ is invertible,
$c= 0$ follows. 

Therefore, for all $x$ in $R$, $ax = rxsd$. The hypothesis on $cx+d$
implies that $d$ is invertible, so that $a$ is as well. Thus $x = a^{-1}r
x sd$ for all $x$ in $R$. Setting $x = 1$, we obtain $(sd)^{-1} = a^{-1}r$, from 
which $a^{-1}r$ is central and invertible.
\qed
\endcomment

Now we obtain  a recursive formula for $\S_k \equiv \S_k (a_1,a_2, \dots,
a_k): = \J T_{a_k} \J T_{a_{k-1}}\cdots \J T_{a_1}$. We first observe
that $\S_1$ sends $x$ to $(x+ a_1)^{-1}$  (now viewing $\S_1$ as a 
densely defined transformation), and $\S_2 (x) = ((x +a_1)^{-1} +
a_2)^{-1}$; this simplifies to $x \mapsto (x+a_1)(a_2 x + a_2 a_1 +
1)^{-1}$. We attempt to define the  polynomials (in the noncommuting
variables $a_1, \dots, a_k$), $p_k$, $q_k$, $P_k$, $Q_k$ via the relation
$\S_k (x) = (p_k x + q_k)(P_k x + Q_k)^{-1}$. There is no a priori reason
why these should exist; however, we show they do exist by means of the
recurrence relations they satisfy. These are in fact  the noncommutative polynomials introduced by Wedderburn [W]. We use upper case $P$s and $Q$s rather than his lower case $p$s and $q$s. 

For $k =1$ and $k=2$, they exist, and we have $p_1 = 0$, $q_1 = 1$, $P_1
= 1$, $Q_1= a_1$,  and $p_2 = 1$, $q_2 = a_1$, $P_2 = a_2$, $Q_2 = 1 +
a_2  a_1$. Assume the $p$s and $q$s (upper and lower case) have all been
defined for $0\leq i \leq k$ (for some $k\geq 2$). Then $\S_{k+1} (x) = \J
T_{a_{k+1}}\S_k (x)$, and this is $((p_k x + q_k)(P_k x + Q_k)^{-1} +
a_{k+1})^{-1}$. This equals $(P_k x + Q_k)((p_k + a_{k+1}P_k)x + (q_k +
a_{k+1}Q_k))$. We infer that $p_{k+1}$, $q_{k+1}$, $P_{k+1}$, $Q_{k+1}$
exist if they are defined as 
$$\eqalign{  p_{k+1} := P_k \qquad\qquad & q_{k+1} := Q_k \cr  P_{k+1} :=
P_{k-1} + a_{k+1}P_{k} \qquad\qquad & Q_{k+1}:= Q_{k-1} + a_{k+1}Q_k. \cr
}$$

The first few terms in the $P$ series are $1, a_2, 1 + a_3 a_2, a_2 + a_4
+ a_4 a_3 a_2$; the corresponding $Q$s are $a_1, 1 + a_2 a_1, a_1 + a_3 +
a_3 a_2 a_1, 1 + a_2 a_1 + a_4a_1 + a_4 a_3 + a_4 a_3 a_2 a_1$.

From Lemma\,\oneone, we deduce that  for a sequence$ a = (a(1),a(2),\dots, a(k))$ and elements $r$ and
$s$ of GL$(R)$, the equation (as operators) $\S_k (a) = \M r,s$ implies all of the following:
$$\eqalign{
Q_{k-1}(a) & =0 = P_k (a), \quad \text{and for some central invertible $\lambda$,} \cr
P_{k-1}(a) & = \lambda r \cr
Q_k(a) & = \lambda s^{-1}\cr
}$$

\Lem Proposition \onetwo. Suppose that $\S_k (a_1,\dots, a_k) = \M x,y$
for some $k \geq 3$. Then there exists a central invertible element
$\lambda$ \st $y^{-1} = \lambda^{-1} Q_{k-2} $ and $x =
\lambda^{-1}P_{k-1}^{-1}$; in particular, $Q_{k-2}$ and $P_{k-1}$  are
invertible. Moreover, $a_{k-1} = - Q_{k-3}Q_{k-2}^{-1}$ and $a_k = -
P_{k-2}P_{k-1}^{-1}$.  Finally, $x = P_{k-3} - Q_{k-3} Q_{k-2}^{-1}
P_{k-2}$.

\Rmk A consequence is that if $\S_k = \M x,y$, then $\M x,y$, $a_{k}$,
and $a_{k-1}$ are uniquely determined by $\brcs{a_1, \dots, a_{k-2}}$.
The converse holds as well, as will be seen Proposition \fouthr. 

\Pf As a consequence of  Lemma\, \oneone\ and the definitions,  $Q_{k-1}
= 0$, $P_{k} = 0$, and then it follows that $P_{k-1} = \lambda x$ and
$Q_k = \lambda y^{-1}$. From the recurrence relations, $Q_{k-2} = Q_k$,
so $Q_{k-2}$ is invertible. 

Since $0 = Q_{k-1} = Q_{k-3} + a_{k-1} Q_{k-2}$, we deduce that $a_{k-1}
= - Q_{k-3}Q_{k-2}^{-1}$, whence $a_{k-1} = -\lambda Q_{k-3} y$. Since $0
= P_{k} = P_{k-2} + a_k P_{k-1}$, we similarly deduce $a_k = -
\lambda^{-1}P_{k-2} x^{-1}$. Finally, $x = P_{k-1} = P_{k-3} +
a_{k-1}P_{k-2}$ yields the final equality.
\qed

In some cases $\M -1,1$ is already in the group generated by $\J$ and
$\T$. This will happen if $\lambda^2 =-1$ for some central element $\lambda$, we can find $a_i$ \st t $\S_3 (a_1,a_2,a_3) $ is the identity, and so
every element of the group can be written as a string with an even number
of $\J$s. 

\comment
For future reference (Section 8), we define $\G_1 (R)$ to be the subgroup of $\G (R)$ generated by $\brcs{\J, \T}$, and $\G_2 (R)$ to be the subgroup generated by $\Set{\J T_a \J T_b}{a,b \in R}$. (The latter corresponds to $\text{PSL}(2,R)$.)
\endcomment

\SecT 2 Noncommutative continued fractions

Let $R$ be a unital ring (not necessarily commutative), and let
$\brcs{a_i}_{i=1}^{\infty}$ be a sequence of elements therein. Define two
sequences of noncommutative polynomials recursively, as given in section\,1. These were introduced by Wedderburn [W] (who used lower case $p$s and $q$s). These are called {\it continuants\/} in [C, p\,25]. 

$$\eqalign{ P_1(a_1) &= 1, \quad P_2 (a_1,a_2) = a_2, \qquad
P_k(a_1,\dots,a_k) = P_{k-2} + a_k P_{k-1} \cr Q_1(a_1) & = a_1 \quad Q_2
(a_1,a_2) = 1+ a_2 a_1, \qquad Q_k(a_1,\dots,a_k) = Q_{k-2} + a_k
Q_{k-1}. \cr }$$ We have suppressed the arguments of the $P$s and $Q$s 
on the right for convenience; unless otherwise noted, $P_k $ represents
$P_k (a_1,
\dots, a_k)$. It is routine to see that $P_k (a_1, \dots, a_k) =
Q_{k-1}(a_2, \dots, a_k)$. We similarly define their opposites,
$P_k\Op$ and
$Q_k\Op$,
$$\eqalign{ P_1\Op(a_1) &= 1, \quad P_2\Op (a_1,a_2) = a_2, \qquad
P_k\Op(a_1,\dots,a_k) = P_{k-2}\Op + P_{k-1}\Op a_k\cr Q_1\Op(a_1) & =
a_1 \quad Q_2\Op (a_1,a_2) = 1+ a_1 a_2, \qquad Q_k\Op(a_1,\dots,a_k) =
Q_{k-2}\Op +  Q_{k-1}\Op a_k. \cr }$$

We observe that $Q_k\Op (a_1,\dots, a_k) = Q_k (a_k,\dots, a_1)$ (as
follows easily by induction), although the corresponding relation for the
$P$s is more complicated. More interesting are the following relations,
which are reminiscent of continued fractions.

\Lem Proposition \twoone. (Mostly in [W, Theorem 3]) The following hold: 
\item{(i)} $P_k Q_k\Op = Q_k P_k\Op$;
\item{(ii)} $P_{k-1} Q_k\Op - Q_{k-1}P_k\Op = (-1)^k$;
\item{(iii)} $Q_k P_{k-1}\Op - P_k Q_{k-1}\Op = (-1)^k$.

\Pf The proof is by induction on all three parts together. It is
straightforward to verify all the cases for $k\leq 3$.
We assume that all these are true for all $k < n$ (with $n\geq 4$), and
verify each of the assertions for $k= n$.

\noindent For (i), we expand $P_n Q_n\Op$:
$$\eqalign{ P_n Q_n\Op & = (P_{n-2} + a_n P_{n-1})(Q_{n-2}\Op +
Q_{n-1}\Op a_n)\cr & = P_{n-2}Q_{n-2}\Op + a_n P_{n-1}Q_{n-1}\Op a_n + a_n
P_{n-1}Q_{n-2}\Op + P_{n-2}Q_{n-1}\Op a_n \cr & = Q_{n-2} P_{n-2}\Op +
a_n Q_{n-1}P_{n-1} a_n + a_n (Q_{n-1}P_{n-2}\Op - (-1)^{n-1})  +
((-1)^{n-1} + Q_{n-2}P_{n-1}\Op)a_n \cr & = Q_{n-2} P_{n-2}\Op + a_n
Q_{n-1}P_{n-1} a_n + a_n Q_{n-1}P_{n-2}\Op + Q_{n-2}P_{n-1}\Op a_n  \cr &
= (Q_{n-2} + a_n Q_{n-1})(P_{n-2}\Op + P_{n-1}\Op) = Q_n P_n\Op. \cr }$$

\noindent The remaining two are even simpler.
$$\eqalign{ P_{n-1}Q_n\Op - Q_{n-1}P_n \Op & = P_{n-1} (Q_{n-2}\Op + 
Q_{n-1}\Op a_n) - Q_{n-1}(P_{n-2}\Op+ P_{n-1}\Op a_n)\cr & = (P_{n-1}
Q_{n-2}\Op - Q_{n-1}P_{n-2}\Op) +  (P_{n-1}Q_{n-1}\Op -
Q_{n-1}P_{n-1}\Op) a_n \cr & = - (-1)^{n-1} + 0 = (-1)^n .\cr }$$

$$\eqalign{ Q_n P_{n-1}\Op - P_n Q_{n-1}\Op & =  (Q_{n-2} + a_n
Q_{n-1})P_{n-1}\Op - (P_{n-2}+a_n P_{n-1})Q_{n-1}\Op\cr & = (Q_{n-2}
P_{n-1}\Op - P_{n-2}Q_{n-1}\Op) + a_n (Q_{n-1}P_{n-1}\Op -
P_{n-1}Q_{n-1}\Op) \cr & = - (-1)^{n-1} + 0 = (-1)^n .\cr }$$
\qed

\SecT  3 Matrix trickery

In this section, we obtain more idenities on the $P$s and $Q$s (most of which can be found in [W]) via $2 \times 2$ matrices. This leads to a one-parameter family of noncommutating polynomials $Q_i$ \st $Q_i$ is invertible if $Q_i\op$ is, Corollary \thrthr\ (the basic example is $1 + ab$ and $1 + ba$). Under a direct finiteness assumption, we have a similar result for equality to zero, Corollary \thrtwo.

 For $a(1), a(2), \dots, a(k)$ in $R$, define the matrices $\P_k \equiv
\P_k (a(1), \dots, a(k))$ and $\P_k\Op$ via
$$\eqalign{
\P_k &= \(\matrix 
P_{k-1} & Q_{k-1} \\
P_{k} & Q_{k} \\
\endmatrix\)  = \(\matrix 0 & 1 \\ 1 & a(k) \\ \endmatrix\) \P_{k-1} \cr
\P_k\Op &= \(\matrix 
P_{k-1}\Op & P_{k}\Op \\
Q_{k-1}\Op & Q_{k}\Op \\
\endmatrix\)  = \P_{k-1}\Op \cdot \(\matrix 0 & 1 \\ 1 & a(k) \\
\endmatrix\), \cr
}$$
where $P_k $ is short for $ P_k (a(1),\dots, a(k))$ etc. Note that $\P_k
\Op$ is obtained by applying ${}\Op$ to the entries and then transposing.
These lead to the factorizations
$$\eqalign{
\P_{k} &= \(\matrix 0 & 1 \\ 1 & a(k) \\ \endmatrix\) \(\matrix 0 & 1 \\ 1
& a(k-1) \\ \endmatrix\)  \cdots \(\matrix 0 & 1 \\ 1 & a(1) \\
\endmatrix\)\cr
\P_{k}\Op &= \(\matrix 0 & 1 \\ 1 & a(1) \\ \endmatrix\) \(\matrix 0 & 1
\\ 1 & a(2) \\ \endmatrix\)  \cdots \(\matrix 0 & 1 \\ 1 & a(k) \\
\endmatrix\)   \cr 
}$$

Briefly referring to the concrete case ($\G (R)$ is defined and there exists central $k$ \st $k(k-1)$ is a non-zero-divisor), Lemmas \,\oneone,   together with this definition 
yields that for a sequence $a = (a(1),\dots a(n))$, and invertibles $r,s$, 
the equation $\S_k (a) = \M r,s$ implies there exists central invertible $\lambda$ \st
 $\P_k (a) = \lambda \(\smallmatrix r & 0 \\ 0 & s^{-1} \\ \endsmallmatrix\)$.

Since $\(\smallmatrix 0 & 1 \\ 1 & a \\ \endsmallmatrix\)   =
\(\smallmatrix 0 & 1 \\ 1 & 0 \\ \endsmallmatrix\) \(\smallmatrix 1 & a
\\ 0 & 1 \\ \endsmallmatrix\) $, we have that each of $\P_k$ and
$\P_k\Op$ is invertible, and moreover
$$
\P_k \equiv \P_k\Op \equiv \(\matrix 0 & 1 \\ 1 & 0 \\ \endmatrix\)^k
\(\matrix 1 &
\sum_{i=1}^k a(i) \\ 0 &1 \\ \endmatrix\) \mod \GL {2,R}'.
$$
(Examine the homomorphic images of the factorizations in any abelian
group.)

The identity 
$$\eqalign{
\M r,s T_b (\M r,s)^{-1} & = T_{rbs} \quad \text{yields} \cr
\M r,s T_b (\M r,s)^{-1}T_{-b} & = T_{rbs -b}, \cr
}$$
showing that $T_a$ is a commutator in $\H$ (and thus in $\G$) if there exist $b$ in $R$ and $r,s$ in $\gl (1,R)$ \st $rbs-b = a$. 

This corresponds to the equation
$$
\(\matrix r & 0\\ 0 & s^{-1}\\ \endmatrix \)
\(\matrix 1 & b\\ 0 & 1\\ \endmatrix \)
\(\matrix r^{-1} & 0\\ 0 & s\\ \endmatrix \)
\(\matrix 1& -b\\ 0 & 1\\ \endmatrix \) = \(\matrix 1 & rbs-b\\ 0 & 1\\ \endmatrix \),
\tag 3$$
which is easy enough to  verify directly (the expression $\M r,s$ corresponds to $\diag (r,s^{-1})$). 

The group $\text{E}(2,R) $ is the subgroup of $\gl (2,R)$ generated by 
$$
\(\matrix 1& a\\ 0 & 1\\ \endmatrix \),\(\matrix r & 0\\ 0 & s\\ \endmatrix \), \(\matrix 0 & 1\\ 1 & 0\\ \endmatrix \),
$$
where $a$ varies over all elements of $R$, and $r,s$ vary over $\gl (1,R)$. 

The following consists of  very slight refinement of  standard results in lower K-theory.

\Lem Lemma \throne. Let $R$ be a ring \st for all $a $ in $R$, there exist $b$ in $R$ and invertibles $r,s$ \st $rbs-b = a$. \noindent 
Then \item{(i)}$\(\smallmatrix 1 & a \\ 0 & 1 \\
\endsmallmatrix\)$ is a commutator in upper triangular matrices;
\item{(ii)} $\(\smallmatrix 0 & -1 \\ 1 & 0 \\
\endsmallmatrix\)$ 
is in the commutator subgroup of $\text{E}(2,R)$;
\item{(iii)} if additionally, there exists a central $\lambda$ in $R$ \st $\lambda^2 = -1$, then $\(\smallmatrix 0 & \lambda \\ \lambda & 0 \\
\endsmallmatrix\)$ is in the commutator subgroup of $\text{E}(2,R)$. 

\Rmk If $R = \Z$, the conclusion is not true---consider the
image modulo $2$ of the matrix with $a=1$---it is not in the commutator
subgroup of  $\gl(2,\Z_2) \iso S_3$.

\Pf (i) A consequence of $(3)$. 

\noindent (ii) We have the equation 
$$
\(\matrix 1 & a \\ 0 & 1 \\ \endmatrix\)
 \(\matrix 0 & -1 \\ 1 & 0 \\ \endmatrix\) 
 = \(\matrix a & -1 \\ 1 & 0\\ \endmatrix\) ,
$$
so that $\(\smallmatrix 0 & -1 \\ 1 & 0 \\
\endsmallmatrix\)$ is in the commutator subgroup of $\text{E}(2,R)$ if $\(\smallmatrix 1 & -1 \\ 1 & 0 \\
\endsmallmatrix\)$ is ($a=1$). We also have 
$$
\(\matrix 1& 0\\ 1 & 1\\ \endmatrix \)
\(\matrix 0 & 1\\ 1 & 0\\ \endmatrix \)
\(\matrix 1 & 0\\ -1 & 1\\ \endmatrix \)
\(\matrix 0& 1\\ 1 & 0\\ \endmatrix \) = 
\(\matrix 1 & -1\\ 1 & 0\\ \endmatrix \),
$$
so we are done.

\noindent (iii) Parallelling (ii), we have 
$$\eqalign{
 \(\matrix 0 & \lambda \\ \lambda & 0 \\ \endmatrix\) 
 &= \(\matrix 1 & \lambda \\ 0 & 1 \\ \endmatrix\)\(\matrix 1 & \lambda \\ \lambda & 0\\ \endmatrix\) \cr
 \(\matrix 1& 0\\ \lambda & 1\\ \endmatrix \)
 \(\matrix 0 & -1\\ 1 & 0\\ \endmatrix \)
\(\matrix 1 & 0\\ -\lambda & 1\\ \endmatrix \)
\(\matrix 0& 1\\ -1 & 0\\ \endmatrix \) &= 
\(\matrix 1 & \lambda\\ \lambda & 0\\ \endmatrix \).
}$$
\qed

In particular, if $k$ is even and $R$ satisfies any of the weak conditions 
below, then $\P_k$ belongs to $\gl(2,R)'$ for any choice of
$a(1), a(2), \dots, a(k)$ in $R$. 

\Lem Lemma \thrfiv. Any of the following conditions on $R$ imply that for all $a$, there exist $b$ in $R$ and invertible $r,s$ \st $a = rbs-b$. 
\item{(a)} There exist invertible $r$ and central invertible $s$ in $R$ \st $r- s^{-1}$ is invertible.
\item{(b)} There exists $r$ in $R$ \st both $r$ and $r-1$ are invertible.
\item{(c)} The centre of $R$ contains a field with at least three elements.
\item{(d)} $2$ is invertible in $R$.
\item{(e)} There exists a unital subring $R_0 \subseteq R$ \st $R_0$ is isomorphic to a finite direct product $\prod \Mn_{n(i) }R_i$ for some rings $R_i$
 with all $n(i) \geq 2$.

 \Pf (a) For $s$ central, $rbs-b = (rs-1) b$; since $rs-1 = (r-s^{-1})s$, it is invertible. Set $b = (rs-1)^{-1}a$. 
 
 \noindent (b) Apply (a) with $s = 1$.
 
 \noindent (c) If $K$ is a field in the centre, the polynomial $x(x-1)$ only has roots $0,1$ (in $K$); if we pick any $r$ in $K \setminus \brcs{0,1}$, then  $r(r-1)$ will be invertible in $K$ (and thus in $R$), so (b) applies.
 
 \noindent (d)  Apply (b) with $r = 2$.
 
 \noindent (e) We verify (b); it is enough to prove it for $R_0$. For each $i$, let $A_i$ be the companion matrix (with entries $0,\pm 1$) of the polynomial $x^{n(i)} - x + 1$, viewed as an integer matrix of size $n(i)$. Then $A_i (A_i^{n(i)}-\I_{n(i)}) = \I _{n(i)} = (A_i^{n(i)}-\I_{n(i})) A_i$, so that $A_i$ is invertible in $\Mn_{n(i)} \Z$. Since $A_i -\I_{n(i)}= A_i^{n(i)}$, it follows that $A_i -\I_{n(i)}$ is invertible.
 
 Let $\Arrow \phi_{i}; \Z . R_i$ be the unital homomorphism; this extends to a ring homomorphism 
 $$\Arrow \Phi_{i}; \Mn _{n(i)}\Z . \Mn_{n(i)}R_i,
 $$
  and set $r_i = \Phi_i (A_i)$. It is routine to check that  $r:= (r_1, \dots, r_n) \in R_0$ is invertible and so is $r-1$. 
 \qed

 The recurrence relations
among the
$P$s and
$Q$s also yield 
$$
\P_k^{-1} =  (-1)^k \(\matrix Q_k\Op
& -Q_{k-1}\Op\\ -P_k\Op & P_{k-1}\Op \\ \endmatrix\)
$$
---actually the recurrence relations only yield that the matrix on the
right is a right inverse for $\P_k$, but $\P_k$ is already known to be
invertible, owing to its factorization above. In particular,
$$
\(\matrix 0 & -1 \\ 1 & 0 \\ \endmatrix \) \P_k^{-1} \(\matrix 0 & 1\\ -1
& 0
\\
\endmatrix \) = (-1)^k \P_k\Op,
$$
a conjugacy.

From $\P_k^{-1} \P_k = \I$, we deduce four more relations (which can also
be proved inductively, just as the first batch were), this time with the
${}\Op$ terms on the left:

\item{(**)} $Q_{k}\Op P_{k-1} - Q_{k-1}\Op P_{k} = (-1)^k$  
\item{} $P_{k-1}\Op Q_{k} - P_{k}\Op Q_{k-1} = (-1)^k$  
\item{} $Q_{k}\Op Q_{k-1} = Q_{k-1}\Op Q_{k} $ 
\item{} $P_{k}\Op P_{k-1} = P_{k-1}\Op P_{k} $.

If $Q_{k-1} = 0$, then $\P_k = \(\smallmatrix a & 0 \\ b & c  \\
\endsmallmatrix\)$ is invertible. From (two-sided) invertibility of
$\P_k$, we quickly  deduce that $a$  is right invertible (that is, $ad
=1$ for some $d$ in $R$) and $c$ is left invertible in $R$. If $R$ is
{\it directly finite\/} ($xy = 1$ implies $yx =1$), then it is easy to see
now that the upper right entry of $\P_k^{-1}$ must be zero and both
diagonal entries, $a$ and $c$, are invertible. Similar results hold if we
assume that $Q_{k-1}\Op$ is zero.

\Lem Corollary \thrtwo. Suppose that  $R$ is directly finite. Let $m$ be  a
positive integer and $a(1), \dots, a(m)$ be elements of $R$.
Then  $Q_m (a(1), \dots , a(m))= 0$  
if and only if $Q_m\Op (a(1), \dots , a(m))= 0$, and this entails that
$Q_{m-1}(a(1), \dots , a(m-1)) $ and $P_m (a(1), \dots , a(m))$ are
invertible. 

For $ m = 2$, this is simply a restatement of direct finiteness, $a_1 a_2 = -1$ entails
$a_2 a_1 = -1$. If $m=3$, it is  somewhat more interesting:  if
$a$, $b$, and $c$ are any three elements in a directly finite ring \st $a +
c + cba = 0$, then $a + c + abc = 0$. This by itself implies direct  finiteness---set $b =1$, 
and notice that $a + c + ac = 0$ is equivalent to $(a+1) (c+1) = -1$. More generally (but not
very interestingly), if $Q_m = 0$ implies $Q_m\Op = 0$ for {\it some\/} $m \geq 2$, then it implies
weak finiteness, hence the implication for {\it all\/} $m$. To see this, simply observe
$$\eqalign{
Q_m(a_1,a_2,1,0,0\dots,0) &= \cases Q_2(a_1,a_2) & \text{if $m$ is even}\\
Q_3(a_1,a_2,1) & \text{if $m$ is odd}\\
\endcases \cr
{}
Q_m\Op(a_1,a_2,1,0,0,\dots,0) &= \cases Q_2\Op(a_1,a_2) & \text{if $m$ is even}\\
Q_3\Op(a_1,a_2,1) & \text{if $m$ is odd}\\
\endcases
}$$

If $R$ is not directly finite, weird things can happen. 

\Lem Lemma \thrsix. Suppose $R$ is not directly finite. Then there is an invertible upper triangular matrix in $\Mn_2 R$ whose inverse is not upper triangular.  

\Rmk The converse is trivial.

\Pf Suppose that $xy = 1 \neq yx$. Define 
$$
A = \(\matrix x & 0\\ 1-yx & y \endmatrix\) \qquad 
B = \(\matrix y & 1-yx\\ 0 & x \endmatrix\).
$$
It is straightforward to verify that $AB = \I = BA$. \qed

For this example, if $R$ is  the C$^*$-algebra of bounded linear operators on $l^2(\N)$, and $x,y$ are the left shift and the right shift respectively, the matrices $A$ and $B$ are even unitary (and this of course extends any non-directly finite C$^*$-algebra). These matrices are in the elementary subgroup. 

\Lem Corollary \thrthr. Suppose that $n \geq 2$  and $Q_{n}(a_1, \dots,
a_{n})$ is invertible in $R$.  Then $Q_n\Op$ is invertible, with inverse
$(-1)^n(P_{n-1}- Q_{n-1}Q_n^{-1}P_n)$. In particular, $Q_n (a_n,
a_{n-1},a_{n-2},\dots, a_1 )$ is invertible if $Q_n (a_1, \dots, a_n)$ is.

\Pf We expand
$$\eqalign{ (P_{n-1}- Q_{n-1}Q_n^{-1}P_n) Q_n\Op & = P_{n-1} Q_{n}\Op -
Q_{n-1}Q_n^{-1}P_n Q_n\Op  \cr & = P_{n-1} Q_{n}\Op - Q_{n-1}Q_n^{-1}Q_n
P_n\Op \cr & = P_{n-1} Q_{n}\Op - Q_{n-1} P_n\Op \cr & = (-1)^n.\cr 
}$$
This yields that $Q_n\Op$ has a left inverse. However, we can now use the
relations appearing in (**) to similarly show that the product the other
way around is also $(-1)^n$, hence $Q_n\Op$ is invertible.

Finally, $Q_n\Op = Q_n (a_n, a_{n-1},a_{n-2},\dots, a_1 )$, as
follows from an easy induction.
\qed

When $k=2$, we have the well-known result that $1 + a_2 a_1$ is
invertible if and only if $1+ a_1 a_2$ is. This is usually used to prove
that the Jacobson radical, defined as the intersection of maximal right
ideals, is also the intersection of the maximal left ideals, and is thus
a two-sided ideal. It also yields that if $R$ is a Banach algebra, then
$\spec (a_1 a_2) \cup \brcs{0} = \spec (a_2 a_1) \cup
\brcs{0}$. 

When $k =3$, the result is that $a_1 + a_3 + a_3 a_2 a_1$ is
invertible if and only if $a_1 + a_3 + a_1 a_2 a_3$ is; this was observed in [MM, where Wedderburn's polynomials are called {\it continuants}] and subsequently in
[MV, just prior to Lemma 1.5]. 

When $k =4$, the expressions are
$1 + a_4 a_1 +  a_2 a_1 + a_4 a_3 + a_4 a_3 a_2 a_1$ and $1 +  a_1 a_4 +
a_3 a_4 + a_1 a_2 + a_1 a_2 a_3 a_4$. Unlike the situation with $Q_2$, there exist
examples over $R = \Mn _2 \C$ with $Q_3(a)$ and
$Q_3(a)\Op$ having different nonzero spectra.

\Lem Lemma \thrsev. Suppose that $R$ satisfies the hypotheses of Lemma \throne. If $Q_k (a_1, a_2, \dots, a_k)$ is invertible, then 
$$
\(\matrix Q_k \cdot (Q_k\op)^{-1}  & 0 \\
0 & 1 \\
\endmatrix\) \in \gl (2,R)'.
$$
If $1$ is in the stable range of $R$, then $Q_k \cdot (Q_k\op)^{-1} \in \gl (1,R)'$.

\Pf 
Factor $\P_k = \(\smallmatrix Q_k
& 0 \\ 0 & Q_k \\ \endsmallmatrix\)  \(\smallmatrix a & b \\ c & 1 \\
\endsmallmatrix\)  $ where $a = Q_k^{-1} P_{k-1}$, $b = Q_k^{-1}
Q_{k-1}$, and $c = Q_k^{-1} P_{k-1}$. The factorization
$$
\(\matrix a & b \\ c & 1\\ \endmatrix\) = \(\matrix 1
& b \\ 0 & 1 \\ \endmatrix\) \(\matrix a-bc & 0 \\ 0 & 1 \\ \endmatrix\)
\(\matrix 1 & 0 \\ c & 1 \\ \endmatrix\) 
$$
yields that $a-bc$ is invertible in $R$; of course, $a- bc = Q_k^{-1}\cdot
(P_{k-1} - Q_{k-1} Q_k^{-1} P_k) = (-1)^k Q_k^{-1}\cdot (Q_k\Op)^{-1}$. 

Now $\(\smallmatrix Q_k & 0 \\ 0 & Q_k \\
\endsmallmatrix\) = \(\smallmatrix Q_k & 0 \\ 0 & 1 \\
\endsmallmatrix\)\(\smallmatrix 1 & 0 \\ 0 & Q_k \\
\endsmallmatrix\)$, hence is congruent to $\(\smallmatrix Q_k^2 & 0 \\ 0 &
1 \\
\endsmallmatrix\)$ modulo $\gl(2,R)'$. We conclude that $\(\smallmatrix
(-1)^kQ_k \cdot (Q_k\Op)^{-1} & 0 \\ 0 & 1 \\
\endsmallmatrix\)$ is congruent to $\(\matrix 0 & 1 \\ 1 & 0 \\
\endmatrix\)^k$ modulo $\gl(2,R)'$. If $k$ is even, this means that $\(\smallmatrix
Q_k \cdot (Q_k\Op)^{-1} & 0 \\ 0 & 1 \\
\endsmallmatrix\)$ is in $\gl(2,R)'$. If $k$ is odd, then $\(\smallmatrix
Q_k \cdot (Q_k\Op)^{-1} & 0 \\ 0 & 1 \\
\endsmallmatrix\)$ is congruent to $\(\smallmatrix
0 & -1 \\ 1 & 0 \\ \endsmallmatrix\)$ modulo $\GL {2,R}'$. The last
matrix is a commutator, and so regardless of the parity of $k$,
$\(\smallmatrix Q_k \cdot (Q_k\Op)^{-1} & 0 \\ 0 & 1 \\
\endsmallmatrix\)$ is in $\gl(2,R)'$. 
\qed

This computation was performed for
$k = 3$ in [MV, Proposition 1.5].

Although we do not know that $Q_k \cdot (Q_k\Op)^{-1}$ is in $\GL (R)'$, we
have that its image in $\gl(2,R)$ is in the commutator subgroup. For
many rings $R$, the map $\GL (1,R) \to \gl(2,R)$ given by $r \mapsto \(\smallmatrix
r & 0 \\ 0 & 1 \\ \endsmallmatrix\)$ induces an {\it isomorphism\/} 
$\GL (R)^{\text{ab}} \to \gl (2,R)^{\text{ab}}$ on their abelianizations.
For such rings (and assuming $R$ satisfies the modest conditions of the
lemma, and of course that $Q_k$ is invertible), we deduce that  $Q_k \cdot
(Q_k\Op)^{-1}$ is in
$\gl (2,R)'$. For example, this holds if $R = \Mn _m D$ for any  $ m> 3$ and
any Dedekind domain $D$. If $R$ is
commutative, $Q_k = Q_k
\Op$, so the result is trivial. This stable result also holds for rings
with $1$ in their stable range (but satisfying the modest conditions). 

\Lem Corollary \thrfou. Let $C$ be any commutative ring, and form $S = \Mn
_N C$. For any $k$-sequence $(a_1,
\dots, a_k)$ of elements of $S$, 
$$
\det\, Q_k (a_1,\dots,a_k) = \det\, Q_k\Op (a_1,\dots,a_k).
$$

\Pf First, we prove it for $C = \C$, the complexes. In this case, $R = S
= \Mn _N \C$, the invertibles are dense, the commutator subgroup of the
invertibles is closed and equal to the kernel of the determinant (when
the latter is restricted to $\gl (n,\C)$). We have just seen that if $Q_k$ is
invertible, then so is $Q_k\Op$, and $Q_k^{-1} Q_k\Op \oplus \I$ is in
$\GL (2N,\C)'$. Hence $\det\, Q_k^{-1} Q_k\Op \oplus \I = 1$, from which 
$\det\, Q_k  = \det\, Q_k\Op $. Next, if $Q_k$ is not invertible, then
neither is $Q_k\Op$, hence they simultaneously have determinant zero.
Thus
$\det\,
 Q_k$  and
$\det\,
 Q_k\Op$ agree for all choices of
$N \times N$ matrices. 

Now we use a standard method to go from $\C$  to an arbitrary unital
commutative ring
$C$. Begin with a sequence of $k$ matrices in $S$, $(a_1,\dots, a_k)$.
These involve $l:= kN^2$ entries, and define variables $y_j$, one for each
of the entries. Form the ring $B = \Z [y_1,\dots, y_l]$, the polynomial
ring with integer coefficients. Define the unital ring homomorphism
$\Arrow
\pi;B.C$ that sends $y_i$ to the corresponding element of $C$ (the  
entry of the appropriate matrix). Let $A_i$ be the matrices with entries
from $B$ \st $\pi (A_i) = a_i$ (the entries of $A_i$ constitute a subset
of the set of variables).

We embed $B $ in $\C$ (by taking any set of $l$ numbers in $\C$ that is
algebraically independent over the rational numbers). As elements of $R =
\Mn _N \C$, we have $\det\, Q_k (A_1,\dots,A_k) = \det\,  Q_k\Op
(A_1,\dots,A_k)$; since the definitions of $\det$ and the $Q$s are the
same whether defined over $B$ or $\C$, this equality also holds viewing
the $A_j$  as elements of $\Mn _N B$. Moreover, the definitions of $\det$
and $Q_k$ are compatible with the ring homomorphism $\pi$, so $\det\, Q_k
(a_1,\dots,a_k) = \det\,  Q_k\Op (a_1,\dots,a_k)$.
\qed

\comment
It not true that $Q_k (a_1,\dots, a_k)^{-1} Q_k{\Op} (a_1,\dots, a_k)$
is always a product of commutators; we can easily find $a_1$ and $a_2$
\st $Q_2 (a_1,a_2) =
\(\smallmatrix 1 & 1 \\ 0 & 1 \endsmallmatrix\)$ and $Q_2\Op (a_1,a_2) =
\I$. The former is not a commutator in $\GL (2,\Z)$ (reduce modulo two).
\endcomment

It requires only an elementary induction argument that if each
of $Q_k(a_1, \dots, a_k)$, $Q_{k-1}(a_2,\dots, a_n)$, $Q_{k-2}(a_3,
\dots, a_k)$, \dots, $Q_2 (a_{k-1},a_k)$, and $Q_1 (a_k)$ is invertible
(a fairly hefty assumption), then 
$$Q_k (a_1,\dots, a_k)^{-1} Q_k{\Op}
(a_1,\dots, a_k)$$  is a product of commutators in the ring
generated by $\brcs{a_1, a_2, \dots, a_k}$. 

If
$a_1$ or $a_2$ is invertible, then $Q_2 (a_1, a_2)$ is conjugate to $Q_2
(a_2, a_1) = Q_2\Op (a_1,a_2)$, and so $Q_2^{-1} Q_2\Op$ is a single
commutator. However, if neither $a_1$ nor $a_2$ is invertible, the
conjugacy need not exist (as in $\Mn _N F$), but it usually seems to
happen that $Q_2^{-1} Q_2\Op$ is a commutator.

For the ring $R= \Mn _N \C$, $\spec Q_2 (a_1,a_2) \setminus \brcs{0} = \spec Q_2\Op(a_1,a_2)\setminus \brcs{0}$,
even counting multiplicities (algebraically, not geometrically). However,
nothing like that is true for $k=3$; with $N =2$, we can find
$a_1,a_2,a_3$ for which $Q_2 (a_1,a_2,a_3) $ is the identity matrix, but
the trace of $Q_2\Op (a_1,a_2,a_3) $ can be arbitrary. So it seems that
the determinant is the only symmetric function that is preserved under
$Q_k \mapsto Q_k\Op$. 

Define $\text{E}(2,R)$ to be the subgroup of $\gl (2,R)$ generated by all matrices of the following forms:
$$\eqalign{
t_a &:= \( \matrix  1 & a \\ 0 & 1 \\\endmatrix\) \quad \text{for all $a \in R$} \cr
m_{r,s} &:= \( \matrix  r & 0 \\ 0 & s\\\endmatrix\) \quad \text{for all $r,s\in \gl(1,R)$} \cr
j &:= \( \matrix  0 & 1 \\ 1 & 0\\\endmatrix\) . \cr
}$$
And $\text{PE}(2,R)$
 will denote the projective elementary group $\text{E}(2,R)/Z(R)$, where $Z(R)$ nonstandardly denotes the group of invertible elements in the centre of $R$. Brackets ($[ \ ]$) will denote equivalence class in $\text{PE}(2,R)$ (although these will often be suppressed---for example, using $t_a$ to denote $[t_a]$---for readability). 
 
 \Lem Proposition \threig. For $k \geq 2$, let $a = (a(1), \dots , a(k))$. Then 
 $$
[ \P_k (a)] = [\m r,s]
 $$
 for some invertible $r,s$ in $R$ if and only if $Q_{k-1} = 0 = P_{k}$, and when this holds, $Q_k (a) = Q_{k-2}(a(1), \dots, a(k-2)$ and is invertible, and 
 $r = P_{k-1} $ and $s= Q_{k-2}^{-1}$. 
 
 \Pf All of this follows from the definition of $\P_k (a)$, its invertibility, and the recurrence relations of the $Q$s. \qed
 
 \Lem Proposition \thrnin. Assume $\G(R)$ exists. Then the assignment
 $$
 \Arrow \phi; \brcs{T_a}_{a \in R} \cup \brcs{\M r,s}_{(r,s) \in \gl (1,R) \times \gl(1,R)\op} \cup \brcs{\J}. \text{PE}(2,R)
 $$
 given by 
 $$\eqalign{
 T_a & \mapsto \left[  \( \matrix  1 & a \\ 0 & 1 \\ \endmatrix\)\right] \cr
 \M r,s & \mapsto \left[ \( \matrix  r & 0 \\ 0 & s^{-1} \\ \endmatrix \) \right] \cr 
 \J & \mapsto  \left[ \( \matrix  0 & 1 \\ 1 & 0 \\ \endmatrix\) \right] \cr 
}$$
extends to an isomorphism of groups $\G(R) \to \text{PE}(2,R)$.

\Rmk Note the $s^{-1}$; this is necessary because multiplication in the second coordinate is in $R\op$.

\Pf Let $w(1), w(2), \dots, w(k)$ be elements of 
$\brcs{T_a}_{a \in R} \cup \brcs{\M r,s}_{(r,s) \in \gl (1,R) \times \gl(1,R)\op} \cup \brcs{\J}$. We first show that if the product $w(1) w(2) \dots w(k)$ is the identity in $\G(R)$, then the product $\phi(w(1)) \phi (w(2)) \cdots \phi(w(k))$ is the identity in $\text{PE}(2,R)$. This would imply that $\phi$ extends to a well defined group homomorphism.

Set $\T =  \brcs{T_a}_{a \in R}$ and $\Cal M =  \brcs{\M r,s}_{(r,s) \in \gl (1,R) \times \gl(1,R)\op} $; these are subgroups of $\G(R)$.  Set  $\H$ to be semidirect product of $\Cal M $ acting on  $T$ in the obvious way (as we have seen before). It is immediate that $\phi$ extends to group homomorphisms with domain $\T$ and $\Cal M$ respectively. Now we shw it extends to a group homomorphism with domain $\H$.

We have
$$\eqalign{
\M r,s T_a & = T_{ras} \M r,s \qquad \text{and verify directly that $\phi(\M r,s) \phi( T_a)  = \phi(T_{ras}) \phi( \M r,s$}) \cr
\M r,s \J & = \J \M s^{-1},{r^{-1}} \qquad \text{and verify directly that $\phi(\M r,s) \phi(\J)  = \phi(\J) \phi( \M s^{-1},{r^{-1}}$}) \cr 
}\tag *$$

The corresponding equalities with upper case letters replaced by equivalence classes of the lower class letters (elements of $\text{PE}(2,R)$ as is easy to verify), and we shall view this as part of (*)

Suppose that all $w(i)$ belong to $\T \cup \H$; we want to show that if $w(1) \dots w(k)$ is the identity of $\G (R)$, then so is $\phi(w(1) \dots \phi (w(k))$. If for two consecutive values of $i$, say $w(i)$ and $w(i+1)$, both belong to $\T$, then we can shorten the length of the sequence by taking their product; similarly, if both $w(i), w(i+1)$ belong to $\Cal M$, we can again shorten the length. So we can assume that $w(i)$ alternate in their membership, $T_{a(1)}, \M r(2),{s(2}, \dots$ (or beginning with a multiplication operator rather than a translation). 

By (*), we can move all the $\M r(i),{s(i)}$ to the right of the string, and thus end up with $T_a \M r,s = 1$. But this is impossible unless $a = 0$ and $r = s^{-1} = \lambda$ for some central invertible $\lambda$. This implies $\phi (T_a) \phi (\M r,s) = [\I]$. Since the elements $[t_b],[ \m t,u]$ also satisfy the corresponding relations in (*), we have that $\phi (w(1) \dots \phi(w(2) \dots \phi (w(n)) = t_a \m r,s$, and this equals $[\I]$, this portion of the proof is completed.

In particular, $\phi$ extends to a homomorphism from $\H$, which we call $\tilde \phi$.

Now assume the general case, allowing some of the  $w(i) $ to be  $\J$. If for a consecutive substring $w(i), w(i+1), \dots, w(j)$ all belong to $\H$, then since $\tilde \phi$
 is a homomorphism on $\H$,   we have $\phi (w(i)) \phi (w(i+1) \dots \phi (w(j)) = \tilde \phi (w(i) w(i+1) \dots w(j))$, and moreover, the product is of the form $T_a \M r,s$. If $w(j+1) = \J$, then from the second line of (*), we can move the $\Cal M$ term the right of $\J$. Then 
 $$\eqalign{
 \phi (w(i)) \phi (w(i+1) \dots \phi (w(j)) \phi (\J) & = \tilde \phi (T_a) \tilde \phi(  \M r,s) )\phi (J) \quad  \text{and by (*) for lower case letters}\cr
 & = \tilde \phi (T_a)\phi(\J)  \phi(  \M s^{-1},{r^{-1}} )\phi (J). \cr
 }$$
 We repeat this if there is another $\J$ term to the right of $w(j+1)$, and may continue until there is at most one $\Cal M$ term, and it appears to at the right end. That is, we are reduced to the case that odd-indexed terms are of the form $T_a$, and even ones $\J$, except the last term is an $\Cal M$. We might have to incorporate an additional  $T_0$ term to get the parity to be correct. We thus have a formula of the form 
 $$
 \J T_{a(1)} \J T_{a(2)} \dots \J T_{a(m)} = \M r^{-1},{s^{-1}} 
 $$
 for suitable $a = (a(i))$  and $r,s$ (not the same ones as in the preceding paragraph. Thus $P_m(a) = \lambda \diag(r^{-1},s)$ for a central invertible $\lambda$ (Corollary \oneone). Now $\P_m (a)$ is a product of terms of the form $\(\smallmatrix 0 & 1 \\ 1 & a(i) \\ \endsmallmatrix\)$; by direct computation, the equivalence class of this (in the projective group) is $\phi (\J) \phi (T_{a(i)}$, and the rest follows.

Thus $\phi$ extends to a group homomorphism. If $\gamma :=  \J T_{a(1)} \J T_{a(2)} \dots \J T_{a(k)}  \M r,s 
$ (an arbitrary element of $\G (R)$) is sent to the identity, then again by Proposition 
, $\P_k (a) = \lambda \diag (r,s^{-1}$. Then, for $x$ in the appropriate domain,
$$
\gamma \circ (\M r,s)^{-1} : x \mapsto (P_{k-1}x + Q_{k-1})(P_k + Q_k x)^{-1} = r^{-1} x s^{-1}.
$$
Thus $\gamma$ is the identity, so the extended $\phi$ is one to one; it is onto by definition. 
\qed

When $R$ has one in the stable range, as is well-known, $\text{E}(2,R) = \gl (2,R)$. Moreover, we have a small bound on the number of elementary matrices required to factorize an element: If $g:= 
\( \smallmatrix a & b \\ c& d \\ \endsmallmatrix \)$ is in $\gl (2,R)$, then $Ra + Rb = R$, so there exists $v$ in $R$ \st $a+bv$ is invertible. Thus 
$\( \smallmatrix a & b \\ c& d \\ \endsmallmatrix\) 
\( \smallmatrix 1 & 0\\ v& 1 \\ \endsmallmatrix\)$ has an invertible element in the $(1,1)$ position, and only a few more elementary matrices are required. At most two applications of the transposition matrix $\( \smallmatrix 0 & 1\\ 1& 0\\ \endsmallmatrix\)$
  are necessary. This last property motivates the notion of order ($\ord$) in Sections 4--6. 
  

 
 \SecT 4 Lengths of words
 
 Here $R$ is a general ring. Recall the definitions of $\m r,s$,$ t_a$, and $j$ from section 3. These generate $\text{E}(2,R)$, and so we redefine them to be their images in $\text{PE}(2,R)$. We also define $e_a := j t_a$.  These generators lead to a notion of length, which in turn, is related to conditions parallel (but not equivalent) to stable range. 
 
 We record relations among the various elements of $\text{PE}(2,R)$. The letters $r,r',s,s'$ refer only to invertible elements. 
 
 \item{(i)} $\m r,s  \m r',{s'} = \m {rr'},{ss'}$ 
 \item{(ii)} $e_a e_0 e_b  = e_{a+b}$
   \item{(iii)} $\m r,s e_a (\m r,s)^{-1} = e_{ras^{-1}}\m sr^{-1},{rs^{-1}}$
   \item{(iv)} (Proposition \threig) $e_{a(k)} \cdots e_{a(1)} = \m r,{s}$ if there exists a central invertible $\lambda$ \st the matrix equation
   $$
   \P_k (a) = \lambda \(\matrix r & 0 \\ 0 & s \\ \endmatrix \)
   $$ 
   holds, where $a = (a(1), a(2), \dots, a(k)$. 
   
   There may be some redundancy in this list; in particular, (iii) might follow from the other conditions.
   
   Condition (i) says that the natural map 
   $$\gl(1,R) \times \gl(1,R)/\brcs{(\lambda,\lambda)} \to \text{PE}(2,R)
   $$ 
 obtained from  $(r,s) \mapsto\m r,s $ is a group homomorphism. 
  It follows that every element of $\text{PE}(2,R)$ can be written as a string (possibly empty) of the form $e_{a(k)} e_{a(k-1)} \cdots e_{a(1)} \m r,s$. We are interested in the shortest possible string representing an element. 
  
  The identities $e_{a}^{-1} = e_0 e_{-a} e_{0}$ and $e_0^2 = 1$ follow from the relations above. Hence if we have a string of $e$s and their inverses, we can always convert it to a string without inverses (using $e_{a}^{-1} = e_0 e_{-a} e_{0}$), without two or more consecutive $e_0$s (as $e_0^2 = 1$), and then arrange that $e_0$ can only appear at the left or right end (using $e_a e_0 e_b = e_{a+b}$).
  
  Thus every element of $\pe$ can be expressed in the form $e_{a(k)} e_{a(k-1)} \cdots e_{a(1)} \m r,s$ where $a(j) = 0$ implies that $j$ is either $k$ or $1$.
  
  Define $\S (a) = e_{a(k)} e_{a(k-1)} \cdots e_{a(1)}$ if $a = (a(1), \dots, a(k))$. This is (the equivalence class) of $\P_k (a)$. The following is a partial converse to an earlier result; subject to only one invertibility hypothesis, for any sequence $a(1), a(2), \dots , a(k-2)$, there will exist $a(k-1), a(k), r, s$ (the latter two invertible) \st $\S (a(1), a(2), \dots , a(k-2), a(k-1),a(k)) = \m r,s$. We will then use this to define and analyze lengths of words.
  
  \Lem Proposition \fouone. Suppose that $k \geq 3$ and $a(1), a(2), \dots , a(k-2)$ is a list of $k-2$ elements of $R$ \st $Q_{k-1}(a(1), a(2), \dots , a(k-2))$  is invertible. Then there exist $a(k-1), a(k)$ in $R$ and $r,s$ in $\gl(1,R)$ \st $\S (a(1), a(2), \dots , a(k-2), a(k-1),a(k)) = \m r,s$, and these four elements are uniquely determined by $a(1), a(2), \dots , a(k-2)$.

\Pf By Corollary \thrthr, $P_{k-3} - Q_{k-3}Q_{k-2}^{-1}P_{k-2}$ (all applied to $a = a(1), a(2), \dots , a(k-2)$) or the truncated version, deleting $a(k-2)$) is invertible. 
Setting $a(k - 1) = -Q_{k-3} Q_{k-2}^{-1}$, we have that $P_{k-1} = P_{k-3} - P_{k-3} - Q_{k-3}Q_{k-2}^{-1}P_{k-2} := x$ is thus invertible. So we define $a(k) = -P_{k-2} x^{-1}$. It is now routine to verify that $P_k = Q_{k-1} = 0$ (now applied with the full list). By Proposition \threig, $\S_k = \m x,{y^{-1}}$ with $y = Q_{k-1}^{-1}$. 

Uniqueness (with $r,s$ defined modulo the action of central invertibles) is straightforward.
\qed

We can now analyze the equations $\S_k (a) = \m r,s$ for small values of $k$ 
  and use these to reduce lengths of strings. We are not aiming for a normal form (this is likely impossible here), but simply to reduce the complexity of the strings. 
  
  Conventionally, $a = (a(1), \dots, a(k))$. When $k \leq 2$,  $\S_k (a) = \m r,s$ entails $a (1) = a(2) = 0$, which gives nothing.
  
  When $k = 3$, in order to get 
  $\S_3 (a) = \m r,s$, we require $Q_1 (a(1))$ be invertible, that is, $a(1): = b$ itself is invertible. When this holds, we have an equation of the form 
  $$\eqalign{
  e_{-b^{-1}} e_{b} e_{-b^{-1}} &= \m b,{-b^{-1}}, \text{which rewrites to } \cr
  jt_b j &= t_{b^{-1}} j t_{-b} \m b,{-b^{-1}}
  }$$
  The left side has two of the obnoxious $j$s while the right side has only one, at a cost of a multiplier at the right, which, inside strings, can be moved to the right end anyway. So if among the $a(i)$s in any string other then the two ends, there is an invertible element, then we can reduce the number of $j$s that appear. This reduction will be used in the proof of the simplicity results.
  
  A similar phenomenon occurs with $k= 4$; the condition requires that $Q_2 (a(1), a(2)) = 1 + a(2) a(1)$ be invertible (which entails invertibility $1 + a(1) a(2)$), and in that case we obtain the equation
  $$
  j t_{-a(2)(1 + a(1) a(2)} j t_{-a(1)(1 + a(2) a(1)^{-1}} j t_{a(2)} j t_{a(1)} = \m {1+ a(1)a(2))^{-1}},{1 + a(2)a(1)} . 
  $$ 
  Conjugating with $j$ and rewriting yields 
  $$
  j t_{a(2)} j t_{a(1)} j = t_{a(1)(1 + a(2) a(1))^{-1}}j t_{a(2)(1 + a(1) a(2)} \m {1+ a(2)a(1)},{(1 + a(1) a(2))^{-1}}.
  $$
  There are three $j$s on the left, but only two on the right, and the $a$s are almost arbitrary; but we still require something to be invertible. 
  
  However, if $k=5$ and $1$ is in the stable range of $R$ (this is absolutely crucial), we do obtain an effective reduction procedure, so that at most two $j$ are required. This was actually done earlier  (end of section 3), but we proceed formally. 
  
  Select $a(1), a(2)$ arbitrarily. Since $Ra(1) + R(1+ a(2)a(1)) = R$ and $R$ has $1$ in the stable range, there exists $a$ in $R$ \st $a(1) + a\cdot  (1+ a(2) a(1))$ is invertible. Set $a(3) = a$. Then Proposition\, \threig\  gives a prescription for $a(4), a(5)$ and invertible $x,y$ \st $\S(a(1), \dots, a(5)) = \m x,y^{-1}$. After conjugating this last expression with $j$ and rewriting, we obtain 
  $$
  j t_{a(2)} j t_{a(1)}j = t_{-a(5)} j t_{-a(4)} j t_{-a} \m {y^{-1}}, x.
  \tag *$$
  Now the number of $j$s on the left is three and the number on the right is two, and $a(1),a(2)$ are completely arbitrary. This means any string $e_{a(1))} \dots e_{a(k)} \m r,s$ can be reduced to a string with just two $e$s (the multipliers being relegated to the right edge). 
  
  To quantify this, we introduce a notion of length for \pe, where $R$ is again arbitrary (no stable range assumption). Attached to each element $g$ of $\pe$ is a set $\Ord (g)$. If $g$ is  written as $e_{a(k)} e_{a(k-1)} \dots e_{a(1)} \m r,s$ where $a(i) = 0$ entails $i \in \brcs{1,k}$, then we define an element of $\Ord (g)$ as follows.
  $$
   \cases  0 & \text{if $k=0$}\\
  1^{-} & \text{if $k = 1$ and  $a(1) = 0$} \\
  k-2 & \text{if $k > 1$ and $a(1) = a(k) = 0$}\\
  k - \frac 32 & \text{$a(k) = 0$ and $a(1) \neq 0$}\\
 k^{-}  & \text{if $a (k) \neq 0$ and $a(1) = 0$} \\
 k & \text{$a(1), a(k)$ are both nonzero}\\
  \endcases
  $$
The reasons for this rather fine division into cases will become (more) apparent when we discuss normalizers of elements. 

Finally, we let $\ord (g)$ (the actual length of $g$) to be the minimum of all the elements in $\Ord (g)$, with $k- 1/2 < k^{-1} < k$. Thus $\ord (\m x,y) = 0$, $\ord (t_a) = 1/2$ unless $a = 0$ (as $t_a = e_0 e_a$), and $\ord (j) = 1^{-}$. 

Bizarrely, the examples in Lemma \thrsix\
 are in $\pe$. The image of $A$ there factors as $ e_{-y}e_xe_-y \m 1,{-1}$; thus $\ord (g) \leq 3$ (and probably is $3$, as the ring is very far from having stable range $1$), and $\ord (g^{-1}) \leq 2\slfrac12$. 
 
 \Lem Lemma \foutwo. Suppose $R$ has one in the stable range. Then for all $g \in \pe$, $\ord (g) \leq 2\slfrac12$.
 
 \Pf Our reduction procedure reduces any string to one with two or fewer $e$s, and thus (when expanded to involve $j$ and $t$s) contributes a term that is at most $2\slfrac12$ to $\Ord (g)$.
 \qed
 
 This was originally observed in [H]
 for matrices over the complexes (which of course have $1$ in their stable range) in terms of the two-sided fractional linear transformations studied there (the emphasis was on fixed points). 
 
 We will see that the converse is true (Proposition \sevsix
) and this will lead to further invariants (coming from bounds on $\ord (g)$). 

It is also possible to study indecomposable strings (those that cannot be reduced in length), but this leads too far afield.

\SecT 5 Normal subgroups

We prove that under modest conditions, $\petwo$ is perfect and equals the commutator subgroup of $\pe$. 

Define $\peone$ and $\petwo$ to be the subgroups of $\pe$ generated by $\Set{e_a}{a \in R}$ and $\Set{e_a e_b}{a,b \in R}$ respectively. If $R = \Mn _n F$ with $F$ a field, then $\petwo = \text{PSL}(2n,F)$ and $\peone$ is the image in $\text{PGL}(2n,F)$ of matrices with determinant $\pm 1$. 

Since $e_0 e_a = t_a$, we have that all $t_a$ belong to $\petwo$. Since $\petwo$ and $\brcs{\m r,s}$ generate $\pe$, to verify that $\peone$ or $\petwo$ is normal in $\pe$, it is sufficient to show that they are stable under conjugation with $\m r,s$. 

Testing this for $\petwo$, we have $\m r,s t_a \m r^{-1},{s^{-1}} = t_{ras}$; from $e_a e_b = j t_a j t_b$, conjugating with $\m r,s$ will similarly yield an element of $\petwo$.

For $\peone$, we observe that it is generated by $\Set {t_a}{a \in R} \cup \brcs{j}$. We require a little lemma.

  \Lem Lemma \fivfiv.  Let $R$ be a ring, and let $z$ be an invertible element. Then $\m z,{z} = e_{z^{-1}} e_1 e_{z-1} e_{-z^{-1}}$. 
  
  \Pf A straightforward computation of $P_3 (a)$, $P_4 (a), Q_3 (a), Q_4 (a)$ where $a = (z^{-1}, 1, z-1, -z^{-1})$. \qed

Then $\m r,{s^{-1}} j \m r^{-1},{s} = j \m_{rs},{(rs)^{-1}}$. With $z = rs$, we see that the conjugate of $j$ is a product of $e$s. 

Thus both $\peone$ and $\petwo$ are normal in $\pe$. 

\Lem Lemma \fivsix. Let $R$ be a ring. 
\item{(a)} $\petwo$ is of index one or two in $\peone$.
\item{(b)} If $R$ satisfies any of the hypotheses in Lemma \thrfiv, and there exists a central element $\lambda$ \st $\lambda^2 = -1$, then $\petwo = \peone$.

\Pf (a) We show $e_0 \, \petwo \cup \petwo = \peone$. An element of $\peone$ that can be written as a string of an even number of $e$s is already in $\petwo$, so consider a product of an odd number, $g = e_{a(1)} e_{a(2)} \cdots e_{a(2k+1}$. Since $t_{a(1)} = e_0 e_{a(1)}$, we can rewrite $g $ as $e_0 \cdot t_{a(1)} \cdot e_{a(2)} \cdots e_{a(2k+1)}$, which is a product of $e_0$ with a product of an even number of $e$s.

\noindent (b) By Lemma \throne, $\(\smallmatrix \lambda & 0 \\ 1&   \lambda  \\ \endsmallmatrix\)$ is a commutator in $\text{E}(2,R)$, and from the matrix equations displayed there, we see that it is in the commutator subgroup, without using any diagonal matrices; more is true: when translated to $e$s, we see that it is expressed as an even number of them. Hence $[\lambda e_0]$ belongs to $\petwo$, but since $\lambda$ is central, this is just $e_0$ (as an element of $\petwo$). 
\qed

The converse of 
is easy to prove in some cases, e.g., $R = \Mn _n F$ ($F$ a field) or other rings with a determinant-like function. 

We show that $\petwo$ is the commutator subgroup of $\pe$ and is perfect under very modest circumstances (of course, this type of result is what one obtains in lower K-theory, but working with the $e$s is efficient and instructive). Our first step is to obtain (crude) bounds on the number of $e$s necessary to express $\m r,s$ as a product thereof. 

\Lem Proposition \fivsev. Suppose that $r,s$ are invertible elements of $R$, and there exist (multiplicative) commutators, $u_1, u_2, \dots, u_k$ in $\gl (1,R)$ together with a central invertible $\lambda$ \st $s^{-1}r = \lambda^2 u_1 u_2 \cdots u_k$. Then there exist $a(1), a(2), \dots, a(8k+4)$ in $R$ \st 
$$
\m r,{s} = e_{a(1)}e_{a(2)} \cdots e_{a(8k+4)}.
$$
In particular, $\m r,{s}$ belongs to $\petwo$. 

\Pf  First let $u = x^{-1}y^{-1} xy$ be a commutator of invertibles $x,y$ in $R$. Set $a = u^{-1}x$, and observe that 
$$
\m a,{a} \m y,{x^{-1}yx} = \m u,{x^{-1}yxa} = \m u,1 .
$$
 Since $a$ is invertible, $\m a,a$ is a product of four $e$s (Lemma \fivfiv).
 
Next, we observe that everything of the form $\m y,{x^{-1}y^{-1}x} $ can also be expressed as a product of four $e$s--set $a_1 = x$ and $a_2 = x^{-1}(y^{-1}-1)$, so that $a_2 a_1 = x^{-1}y^{-1}x -1$ and $a_1 a_2 = y^{-1}$. Thus $y = (1+a_1a_2)^{-1}$ and $x^{-1}yx  = (1 + a_2 a_1)^{-1}$; then $x$. By Proposition \fouone, we can find $a_3, a_4$ \st  $P_3 (a_1,a_2, a_3) = 1+ a_1 a_2$ and $Q_2 = 1+ a_2 a_1$, and $\m {1+ a_1 a_2}, {1 + a_2 a_1} = \S (a_1,a_2, a_3, a_4))$. Thus $\m y,{x^{-1}y^{x}}$ is a product of four $e$s. 

Hence $\m u,1$ is a product of $8$ $e$s for any commutator $u$. Thus if $v$ is a product of $k$ commutators, then $\m v,1$
is a product of $8k$ $e$s. Since the $\lambda $ has no effect on the multiplication operators ($\m \lambda, \lambda$ is the identity), we can assume $\lambda = 1$. Finally, we have $\m r,s = \m rs^{-1},1 \m s,s$, and the left factor is a product of $8k$ $e$s, while the right is a product of four of them.
\qed

\comment

\SecT 4 Back to the concrete case
 
As in section\,1, we assume that $R$ satisfies ($\ddag$) for this section.
The following, a partial converse to Lemma\, \oneone,  says that subject to only one invertibility hypothesis, we
can choose any sequence
$a_1,
\dots, a_{k-2}$ of $k-2$ elements and then there will exist $a_{k-1}$,
$a_{k}$,
$x$ and $y$ \st $\S_k (a_1,\dots, a_k)= \M x,y$. Then we analyze the lengths of words
in $\G(R)$.

\Lem Proposition \fouthr. Suppose that $k \geq 3$ and $\brcs{a_1, \dots,
a_{k-2}}$ is a set of $k-2$ elements \st  $Q_{k-2} (a_1, \dots,
a_{k-2})$ is invertible in
$R$. Then there exist  $a_k$, $a_{k-1}$, and $\M x,y$
\st $\S_k (a_1, \dots, a_k) = \M x,y$, and these are uniquely determined
by $\brcs{a_1, a_2, \dots, a_{k-2}}$. 

\Pf By \thrthr, $P_{k-3} - Q_{k-3} Q_{k-2}^{-1} P_{k-2}$ is
invertible. Setting $a_{k-1} = - Q_{k-3}Q_{k-2}^{-1}$, we have that
$P_{k-1} = P_{k-3} - Q_{k-3} Q_{k-2}^{-1} P_{k-2}:= x$ is thus
invertible, and we can thus set $a_k = - P_{k-2} x^{-1}$, and it is now
routine to verify that $P_k = Q_{k-1} = 0$, and thus $\S_k (a_1, \dots,
a_k) = \M x,y$ with $y = Q_{k-2}^{-1}$.
\qed

We can analyze the equations $\S_k = \M x,y$ for small values of $k$ using Lemma\,\onetwo. When
$k
\leq 2$, the only way to obtain an equation of the form $\S_k =\M x,y$ is
trivially, e.g., $a_1 = a_2 = 0$. When $k = 3$, we require that
$Q_{1}(a_1,a_2,a_3)$ be invertible, that is, $a_1$ is invertible, and
then the equation is (with $a = -a_1^{-1}$)
$$
\J T_{-a^{-1}} \J T_{a} \J T_{-a^{-1}}  =  \M -a,a.
$$ 
This leads to  the equation $\J T_a \J = T_{a^{-1}} \J T_{-a} \M -a,a$
which can be used to reduce the number of $\J$s appearing in words (when
$a$ is invertible). We will see a more effective method of doing this
when we discuss the case that $k=5$.

When $k=4$, the condition becomes invertibility of $Q_2(a_1,a_2) = 1 + a_2 a_1$, and the
equation is
$$
\J T_{-a_2 (1+a_1 a_2)}\J T_{-a_{1}(1+a_2 a_1)^{-1}} \J T_{a_2} \J
T_{a_1}  =  \M {(1+a_1 a_2)^{-1}},{(1+a_2 a_1)^{-1}} .
$$ 
If $1 + a_2 a_1$ is invertible, it leads to a reduction by two in the
number of $\J$s. Conjugate the expression with $\J$, and we obtain 
$\J T_{a_2} \J T_{a_1} \J = T_{a_1 (1+a_2 a_1)^{-1}}\J T_{a_{2}(1+a_1
a_2)}\M {1+a_2 a_1},{1+a_1 a_1} $. However, this is not universally
effective, as the constraint (that $1+ a_2a_1$ is invertible) is not always satisfied.

With $k=5$, however, we do obtain a reduction procedure. Select any
$a_1$, $a_2$ in $R$. If we assume that $R$ has $1$ in the stable range,
then there exists $a$ in $R$ \st $a_1 + a (1+ a_2 a_1)$ is invertible
(since $R(a_1) + R(1+ a_2 a_1)$ is the improper left ideal). Set $a_3 =
a$. Then Lemma \onetwo\ gives the prescription for the $a_4$,
$a_5$, $x$ and $y$ \st $\S_5 (a_1,a_2, a_3,a_4, a_5) = \M x,y$. After
conjugating this with $\J$, the resulting expression can  be rewritten as
$$
\J T_{a_2} \J T_{a_1} \J = T_{-a_5} \J T_{-a_4} \J T_{-a}\M
{y^{-1}},{x^{-1}}, \tag{**}
$$ and this time, the two parameters, $a_1$ and $a_2$, are completely
arbitrary, while at the same time, the number of $\J$ terms drops from
three to two.

This suggests the following  definition. Let $T_i$ represent $T_{a_i}$.
We define a set attached to each element $g$ of $\G$, $\Ord (g)$, as
follows (implicitly, these representations allow the $x$ and $y$ to vary
over the invertible elements)
$$
\text{if } g = \cases \M x,y &\text{then
$0 \in \Ord (g)$} \\
T_1\M x,y &\text{then
$\frac 12 \in \Ord (g)$} \\
\J T_{k-1}\dots \J T_1 \J \M x,y & \text{then
$k^{-} \in \Ord (g)$} \\
\J T_k \J T_{k-1}\dots \J T_1 \M x,y & \text{then $k
\in
\Ord (g)$} \\
 T_k \J T_{k-1}\dots \J T_1 \J \M x,y & \text{then $k \in
\Ord (g)$}  \\ 
  T_{k+1} \J T_k \J T_{k-1}\dots \J T_1  \M x,y &
\text{then
$k +
\frac12 \in \Ord (g)$} \\
 \endcases
$$ 

Finally, define $\ord (g)$ to be the minimum of all the elements of $\Ord
(g)$, where $0 < \frac 12 < 1^{-} < 1 < 1\frac12 < 2^{-} < \dots$. Thus
$\ord (\M x,y) = 0$, $\ord (T_A) = \frac 12$ (provided $A \neq 0$), $\ord
(\J) = 1^{-}$.

Suppose that $g$ has a representation corresponding to Order exceeding
$2\frac12$, i.e., the string $\J T_a \J T_b \J$ appears in this
representation. By (*), we can replace this string by one with just
two $\J$s, and move the right/left multiplication operator to the far
right. The outcome is that the new representation has  fewer $\J$s.
We deduce that $\ord (g) \leq 2\frac12$ for all elements $g$---in other
words, every element of $g$ can be written in (at least) one of the forms
$T_a \J T_b \J T_c \M x,y$, $T_a \J T_b \M x,y$, $T_a \M x,y$ for some
choice of elements of $R$ and $\GL R$. [When we discuss the abstract version of
this construction, we will find that this reduction order $2\slfrac12$ or less
is a consequence of $1$ in the stable range.]

Call a word {\it indecomposable\/} if it cannot be  reduced to a word of
lesser order. For example, $\M x,y$ (order zero), $T_a$ if $a \neq 0$
(order one half), $T_a \M x,y$ if $a \neq 0$ (order one half), $\J$ 
(order $1^{-}$) are obviously indecomposable. The indecomposable words
are easy to recognize, but uniqueness fails, when the order is 2$\frac12$.

 If the word $g = \J T_a \J T_b \M x,y$ were  not
indecomposable, there would be a word of one of the forms, $\J T \J \M
*,{*} $ $T \J T \M *,{*}
$,
$\J T
\M *,{*}
$,
$\M *,{*}
$ (with the appropriate subscripts), equalling $g$. The first quickly
reduces to an equation involving at most two $\J$s. Each of the others
would yield an equation where the left side has three
$\J$s, and the right side is a multiplication operator. It easily follows
that this can only happen if $a$ is invertible (in which case $g$ reduces
to a shorter word) or either $a$ or $b$ is zero (both resulting in 
shorter words).   We deduce that $\J T_a \J T_b \M x,y$ is indecomposable
if and only if
$a$ is neither zero nor invertible and $b$ is not zero.

Shorter words are even easier to check for indecomposability; if $g =\J
T_a
\J  \M x,y$ (order 2${}^{-}$), indecomposability is equivalent to $a$
being neither zero nor invertible. 

With words of order 2$\frac12$, the situation is slightly more
complicated, because we can ask if it can be reduced to a word of order 2
or 2${}^{-}$. We can deal with both of these  simultaneously---suppose
$T_a \J T_b  \J T_c = \J T_d \J T_e \M r,s$ (we have absorbed the $\M
x,y$ into the right hand side; this plays no further role; the case of
equality to an element of order 2${}^{-}$ arises when $e = 0$). This
converts to the equation $ T_{-e}\J T_{-d} \J T_a \J T_b  \J T_c  =  \M
r,s $. We perform the substitution arising from 
$\M r, s T_{-c}
 = T_{-rcs} \M r,s$ and then conjugate the resulting expression with $\J$.
We obtain 
  ($k=4$)
$$
\J T_{-e + rcs}\J T_{-d} \J T_a \J T_b   = \M {s^{-1}},{r^{-1}}.
$$ It  follows that $1 + ba$ is invertible, and conversely, if $1+ba$ is
invertible, we can construct the rest of the terms to solve the
equations. Equating $g$ to still shorter forms yields the necessary
conditions that $a$ is not invertible, and $c$ is not zero. (The
conditions that neither $a$ nor $b$ be zero are subsumed by
noninvertibility of $1 + ba$.) 

We deduce that $T_a \J T_b  \J T_c  \M x,y$ is indecomposable if and only
if none of $b$, $1+ba$, and $1+ cb$ is invertible.

Having obtained the indecomposable forms, we now look at uniqueness. If
$g$ is represented by an indecomposable word of order 1$\frac12$ or less,
it is easy to see that the representation is unique. With order $2^{-}$,
given the indecomposable $\J T_a \J$ (with $a\neq 0 $), set it to $\J T_b
\J \M r,s$, we obtain $\J T_{b-a} \J = \M r,s$, which forces $b = a$. So
we have uniqueness here as well. 

There are two forms for order $2$. The first, $\J T_a \J T_b = \J T_c
\J T_d \M r,s$ yields $T_{-d} \J T_{a-c} \J T_b = \M r,s$ which forces
the middle term to be the identity, and it follows that $a=c$ and $b=d$,
so uniqueness occurs here as well. However, there is a second form, 
$\J T_a \J T_b =  T_c \J T_d \J \M r,s$. This yields $\J T_{-d} \J T_{-c}
\J T_a \J T_b = \M r,s$. This forces invertibility of  $1 + ba$---which
however, is not precluded by indecomposability of $\J T_a \J T_b$. In
fact, when $1+ ba$ is invertible, we can solve for $c$, $d$, $r$, and $s$
so that the equation holds. The outcome is that if $1+ba$ is not
invertible, then the indecomposable of order $2$ is unique (among those
of order two representing $g$), but if $1+ba$ is invertible, then
uniqueness fails, although with an element of the form $ T_c
\J T_d \J \M r,s$. In fact, in this case, there are exactly two forms for
the same element.

For order $2\frac 12$, the situation is again somewhat complicated.  The
equation $T_a \J T_b  \J T_c =  T_d \J T_e  \J T_f \M r,s$ yields
$T_{-e}\J T_{-f}\J T_{a-d}
\J T_b  \J T_c = \M r,s$. Moving the $T_c$ to the other side, we obtain
$T_{rcs -e}\J T_{-f}\J T_{a-d} \J T_b  \J = \M r,s$, and conjugating with
$\J$ yields
$$
\J T_{rcs -e}\J T_{-f}\J T_{a-d} \J T_b = \M s^{-1},{r^{-1}}.
$$ We deduce that $1 + b(a-d)$ is invertible, but also observe that if
this holds, then we can certainly solve for $e$ and $f$ (having been given
$a$, $b$, $c$, and assumed a value for $d$, and thus having $r$ and $s$)
so that the last equation holds, and of course, we may work backwards,
and thus obtain a huge family (parameterized by $d$ \st $1 + b(a-d)$ is
invertible) of expressions of order $2\slfrac12$. We can even obtain
multiple expressions without any multiplication operator at all, if
 $R = \Mn N \C$---since
$b$ is not invertible (else the original expression would not be
indecomposable), there exist (lots of) choices for $z$ \st $z b= b z = 0$;
set $d = a-z$, so that $s = r = 1$.

Now using the general results on $\P_k$ developed above in sections 2 and
3, we can define an abstract version of $\G (R)$ (and  two of its
normal subgroups) for any unital ring.

\SecT 5 Abstract $\G(R)$

We can define a version of $\G (R)$ for {\it any\/} ring, via the relations 
obtained in $\G$.

Let $R$ be a unital ring; define the free group $F \equiv F(R)$ with
generators
$$
\Set{E_a}{a\in R} \cup \Set {M_{r,s}}{r,s \in \GL R}.
$$
Now let $W \equiv W(R)$ be the normal subgroup of $F$ generated by the
following words:
\item{(i)} $M_{ r,s} M_{  r',s'} M_{  rr',s's}^{-1}$;
\item{(ii)} $E_0 E_a E_0 E_bE_0 E_{-(a+b)}$;
\item{(iii)}$M_{r,s} E_a M_{s^{-1},r^{-1}}^{-1} E_{s^{-1}a r^{-1}}^{-1}$; 
\item{(iv)} $ E_{a(k)} \cdots E_{a(1)}M_{r,s}^{-1}$, if
the size two matrix equation  
$$\P_k (a) =
\lambda \( \smallmatrix r & 0\\  0 & s^{-1} \\\endsmallmatrix\)$$
\item{} (where $a = (a(i))$)
holds for some central invertible element,
$\lambda$, of $R$.

\vskip4pt \noindent Parts (i) and (ii) reflect the multiplication and translation operators,
respectively; part (iii)  refers to the relationship between them. There may be some redundancy
in this list (in particular, (ii) might follow from the others). Let $e_a$ and $m_{r,s}$ be the respective images of
the generators modulo $W$, and define $\Ag (R)$ to be $F(R)/W(R)$. If
$\G(R)$ exists,  the assignment $E_a
\mapsto \J T_a$, $M_{r,s} \mapsto \M r,s$ induces a group homomorphism
$\Ag (R)
\to \G (R)$, as follows from the remark after the definition of $\P_k$. It is an isomorphism (\fivtwo). (Note that $e_0
\mapsto \J$.)

In general, (i) implies that the assignment $\GL R \times \GL
R\Op/\brcs{(\lambda,\lambda^{-1})} \to \Ag (R)$ induced by $(r,s) \mapsto
m_{r,s}$ is a group homomorphism. We also notice that the identity $e_a^{-1} =
e_0 e_{-a} e_0$ follows from the element of $W$, $E_a E_0 E_{-a} E_0$  (iv \& i) (and also follows from (ii \& i)), and
moreover, so does $e_0^2 = 1$. In addition, $e_a e_0 e_b = e_{a+b}$. Hence if we
have an arbitrary string of $e$s and their inverses, we can first convert it to a string
of $e$s without any inverses (using $e_a^{-1} =
e_0 e_{-a} e_0$), and then remove any $e_0$ that appear in the interior of the string
(using $e_a e_0 e_b = e_{a+b}$). This means that every word in the $e$s
and their inverses is equal to a word of the form $e_{a(k)} e_{a(k-1)}
\cdots e_{a(1)}$ where
$a(j) = 0$  entails that either $j = 1$ or $j=k$ (or both)---in other words, $e_0$ can
only appear at the ends of the string. This means that every element of $\Ag (R)$ can
be written in the form $e_{a(k)} e_{a(k-1)} \cdots e_{a(1)} m_{r,s}$ where only the
first and last $e$s can be $e_0$ (condition (0) allows us to move the $m$s to the
right, just as we did in the concrete version).

A particular consequence of the relations is that
if $\P_k ((a(i))) = \lambda r \oplus \lambda s^{-1}$, then $e_{a(k)}\cdots e_{a(1)}  =
m_{r,s}$. What is not clear at this point is whether the converse holds, because $W$ is
only the normal subgroup generated by the relations. This is equivalent to
deciding whether a word of the form $e_{a(k)}\cdots e_{a(1)} m_{r,s}$ is the
identity; it is not
at all clear how to decide this.

We also define $\Ag_1 (R)$ and $\Ag_2 (R)$ to be the subgroups of $\Ag (R) $
generated, respectively, by $\brcs{e_a}_{a\in R}$ and $\brcs{e_a
e_b}_{a,b \in R}$. (An alternative definition would define these two
groups by their generators and relations---just remove the $M$s from  (iv)
and delete (i \& iii) entirely---but then we would have to show that these are
subgroups, which is not clear.) Since $\Ag_2 (R)$ and $\brcs{m_{r,s}}$
generate $\Ag (R)$, to verify that $\Ag_i (R)$ is normal in $\Ag (R)$, it
is sufficient to show that each is stable under conjugation with every
$m_{r,s}$. For $i=2$, stability is obvious, as
$m_{r,s} e_a e_b m_{r,s}^{-1} = e_{s^{-1}ar^{-1}} e_{rbs}$. For $i=1$,
$m_{r,s} e_a  m_{r,s}^{-1} = e_{s^{-1}ar^{-1}} m_{z,z}$ where $z =
(sr)^{-1}$. However,  $m_{z,z}$ is a product of
four $e_x$s, specifically, $m_{z,z} = e_{-(z-1)}e_{-1} e_{z-1}e_1$.

We observe that with the appropriate translation to
this abstract setting, Proposition\,\fouthr\ applies. Now we have some routine results on behaviour under ring homomorphisms that
preserve the centre. The  {\it invertibles lifting\/} condition of
\fivone(b) means that if
$s$ is an invertible of $S$ in the image of $\phi$, then there exists an
invertible,
$r$, in
$R$
\st
$\phi(r) = s$. This occurs if (for example)
$R$ has one in the stable range, or if the kernel of $\phi$ is contained in the Jacobson
radical of $R$.

\Lem Lemma \fivone. Let $\Arrow \phi ; R. S$ be a unital ring homomorphism
sending the centre of $R$ into the centre of $S$. Then the map $F(R) \to
F (S)$ given by $E_a
\mapsto E_{\phi(a)}$ and $M_{r,s} \mapsto M_{\phi(r),\phi(s)}$ induces
well-defined group homomorphisms
$\Arrow \Ag(\phi);\Ag (R).\Ag (S)$, $\Arrow
\Ag_1(\phi);\Ag_1 (R).\Ag_1 (S)$, and $\Arrow
\Ag_2(\phi);\Ag_2 (R).\Ag_2 (S)$. {\par}
\item{(a)} If $\phi$ is onto, then $\Ag_j (\phi)$ are onto for $j=1$ and $j=2$.
\item{(b)} If $\phi$ is onto and invertibles lift, then $\Ag (\phi)$ is onto.

\Pf Note that the map $F(R) \to \Ag (S)$ determined by the assignments of the
generators is automatically a homomorphism (since $F(R)$ is free), so once we show
the induced mapping is well-defined, it will   be a group
homomorphism from the quotient of $F(R)$. It thus suffices to show that each of the
generators of $W$ that appear in (0) and (i) in the definition of $\Ag (R)$ is mapped
to $W(S)$. But this is a consequence of the elementary fact that $\phi(P_k ((a(i)))) =
P_k (\phi(a(i)))$ (and similarly with the $Q$s), that is, $\phi$ respects the size two
matrices $\P_k$ (the condition on the centre is essential for this). We thus obtain
an induced map $\Ag (R) \to \Ag (S)$, necessarily a homomorphism of groups. Since
$e_a \mapsto e_{\phi(a)}$, that $\phi$ restricts to maps $\Ag_j (R) \to \Ag_j (S)$ for
$j = 1$ and $j=2$. Moreover, the generators of $\Ag_j (S)$ are obviously in the
image of $\Ag_j (R)$, hence (a) holds.

(b) Since $\Ag (S)$ is generated by $\Ag_1 (S)$ and $\brcs{m_{t,u}}$ (where $t,u$
are invertible in $S$), it suffices to be able to lift $t$ and $u$ to invertible elements
$r$ and $s$ respectively, so that $m_{r,s}$ is an element of $\Ag (R) $ and it is
mapped to $m_{\phi(r),\phi(s)} = m_{t,u}$.
\qed

Of course, what we really want is that if $\phi$ is one to one, then $\Ag (\phi)$
should also be one to one. The difficulty is the same word problem mentioned
above---just because $e_{a(k)} \cdots e_{a(1)} = m_{r,s^{-1}}$, it does not necessarily
follow that that the corresponding relation on the $a(i)$ and $r,s$ hold, that is,
$\P_k(a(i)) = \lambda r \oplus \lambda s^{-1}$---there may be other relations which are
simply consequences of these. Fortunately, there is a way around this word
problem. It still leaves some difficulties, however.

\Lem Proposition \fivtwo. If $\G(R)$ exists, then the natural map $F(R)
\to
\G(R)$ induces isomorphisms $\Ag (R) \to \G(R)$, $\Ag_j (R) \to \G_j (R)$
given by $e_a \mapsto
\J T_{a}$ and $m_{r,s} \mapsto \M r,s$.

\Pf We already observed that elements of $W$ are sent to the identity in $\G(R)$, so
there is an induced group homomorphism $\Ag(R) \to \G(R)$, which also restricts to
group homomorphisms from the $\Ag_j (R)$ to their counterparts. Each is clearly
onto, so it remains to show that the map is one to one. Suppose that an element,
$e_{a(k)}e_{a(k-1)}\cdots e_{1} m_{r^{-1},s}$
of $\Ag (R) $, is mapped to the identity in $\G(R)$. This means that
$$
\J T_{a(k)} \J T_{a(k-1)} \cdots \J T_{a(1)} \M r^{-1},s = \M 1,1.
$$
The left factor  is $\S_k (a(i))$, which we know has the form $x \mapsto (P_{k-1}x +
Q_{k-1}) (P_k x + Q_k)^{-1}$ (where we have abbreviated $P_k (a(i))$ to $P_k$,
etc), and we know that for this to be $\M r,s^{-1} $, we
require that $\P_k (a(i)) = \lambda r \oplus  \lambda s^{-1}$. However, this is
one of the relations in $W$, so that the  word, $E_{a(k)}E_{a(k-1)}\cdots
E_{1} M_{r,s}$, belongs to $W$, and thus $e_{a(k)}e_{a(k-1)}\cdots e_{1}
m_{r^{-1},s} = 1$. Obviously, the restrictions to the subgroups are one to one as well.
\qed

\Lem Proposition \fivthr. Suppose that $\Arrow \phi; T.R$
is a unital ring homomorphism sending the centre of $T$ into the centre of $R$.
Suppose moreover that $\G (R)$ exists. If $\phi$ is one to one, then so is
$\Ag (\phi)$.

\Pf By the preceding proposition, it suffices to prove that the map determined by
$e_t \mapsto \J T_{\phi (t)}$, $m_{u,v} \mapsto \M{\phi(u)},{\phi(v)}$
(compose the induced map with the isomorphism) is one to one. If $g:=
e_{a(k)}e_{a(k-1)}\cdots e_{1} m_{u,v}$ is in the kernel, then we obtain
the equation
$$
\J T_{\phi(a(k))} \J T_{\phi(a(k-1))} \cdots \J T_{\phi(a(1))} \M \phi(u)^{-1},{\phi(v)} =
\M 1,1.
$$
This implies that $\P_k (\phi(a(i))) = \lambda \phi(u) \oplus  \lambda
\phi(v)^{-1}$ for some central invertible $\lambda$ in $R$. Since $\P_k$ commutes
with $\phi$, we deduce that $\lambda \oplus \lambda$ is expressible as a polynomial
(with noncommuting coefficients, of course) in $\phi(a(i))$, $\phi (u)$, and $\phi(v)$
(the matrix equation is $\lambda \I = \phi(\P_k  (a(i)) (\phi (u)^{-1} \oplus \phi (v))$), and
this means that there exists $\mu$ in $T$ \st $\phi (\mu) = \lambda$. Since $\phi$ is
one to one, this $\mu$ is unique; it is invertible, since $\mu =
\P_k(a(i)) (u \oplus v^{-1})$.

For $t$ in $T$, we have that $\phi(t\mu - \mu t) = \phi (t) \lambda - \lambda \phi(t) =
0$ (since $\lambda$ is in the centre of $R$). As $\phi$ is one to one, it follows that
$t\mu = \mu t$, i.e., $\mu$ is in the centre of $T$. We thus have $\phi (\P_k  (a(i))-
\mu u^{-1}
\oplus \mu v^{-1}) = 0$ (extending $\phi$ to the corresponding map on size two
matrix rings, also one to one). Hence $\P_k  (a(i))-
\mu u^{-1} \oplus \mu v^{-1} = 0$, so that $g = 1$.
\qed

In particular, if $D$ is any subring of $ \C$, then $T := \Mn m D \subset
R:= \Mn m
\C$, and the hypotheses apply. In particular, $D$ can be any countably  generated
(or even generated by a set of first  uncountable cardinality)  domain of characteristic
zero. If $D$ is a Dedekind domain and $m > 3$, it will turn out that the
invariant
$\Ag (\Mn m D)/\Ag_2 (\Mn m D)$ is rather interesting---it is
$D^{\text{inv}}/(D^{\text{inv}})^{2m}$, that is, the unit group of $D$ modulo
those that are $2m$ powers.

On the other hand, we do not know that (for example) $R$ weakly or stably finite
implies that $\Ag (R)$ is not trivial, or, in the other direction, that really strongly
infinite ($R$ not a domain, but for all $a$ in $R$, there exist $r$ and $s$ \st
$ras = 1$) implies something drastic about $\Ag (R) $.

However, we can prove the expected results concerning perfectness or
minimal normal subgroup property of
$\Ag_2 (R)$. There is an interesting connection between stable range
conditions and the minimal number of terms required to express elements
of the group, e.g., if a type of stable range condition holds, then there
is a bound on the $k$ required to express every element in
$\Ag (R)$ as $e_{a(k)} e_{a(k-1)}\cdots e_{a(1)} m_{r^{-1},s}$ (actually, we have a
more precise statement, which is made a little more complicated by the fact that we are
using $e_a$ corresponding to $\J T_a$ in the concrete case). 
\endcomment

The following  mimics the proof that $\(\smallmatrix 1 & a \\
0 & 1 \\ \endsmallmatrix\)$ is a commutator (\throne).

\Lem Proposition \sixtwo. Let $R$ be a ring \st for all $a$ in $R$, there exist $b,r$ in $R$ with $r$ invertible \st $rbr- b = a$. Then $\petwo$ is pertect and equals the commutator subgroup of $\pe$. In particular, this applies if there exists central $k$ \st $k(k^2-1)$ is invertible. 

\Rmk Sufficient is that Lemma \thrfiv\ apply; in particular, this holds if  $R$ contains a central field with at least four elements.

\Pf  We note that $\m r,r t_b (\m r,r)^{-1} t_b^{-1}= t_{rbr-b} = t_a$, and thus $t_a = e_0 e_a$ is in $\pe'$; and since $\m r,r$ is in $\petwo $ (Lemma\,\fivfiv),  thus $t_a$ in in $\petwo'$ (every $t_c \in \petwo$). Thus so is $j t_d j = e_c e_0$ for every $d$ in $R$. Since $e_c e_0 e_0 e_a = e_c e_a$, all the generators are commutators, and so $\petwo$ is perfect. 

Since $\pe/\petwo$ is abelian, $\petwo$ contains $\pe'$, and thus equals it.

If $k$ is central and $k(k^2-1)$ is invertible, then $k$ and $k^2 -1$ are invertible. Set $r = k$ and $b = (k^2-1)^{-1}a$, so that $rb r - b = a$.
\qed

\def\SR{\Cal S\Cal R}

\SecT 6 Stable range and orders of elements

We show a connection between properties of the $Q$s (resembling stable range conditions) and maximal lengths of elements (Proposition \sevsev). This permits us to define a sequence of stable range-like conditions; these are called $\SR(n)$. The first one, $\SR (1)$, actually is equivalent to one in the stable range, but the higher ones are not stable range conditions. For example, $\Z$ fails all $\SR (n) $.  

First we have an elementary result about solution of equations
involving $Q_k$s. The last $c_l$ does not need to
be invertible.
This lemma also applies when the $Q_i$ are replaced by $Q_i\Op$.

\Lem Lemma \sevone. Let $c_1, c_2, \dots, c_l$ be $l$  elements in
$R$ \st $c_i$ is invertible for $1 \leq i \leq l-1$. Then there exist
$x(1), x(2),
\dots, x(l)
$ in
$R$
\st the equations
$$\eqalign{
Q_1 (x(1)) & = c_1 \cr
Q_2 (x(1),x(2)) & = c_2 \cr
\vdots \qquad& \quad\vdots \cr
Q_l (x(1), x(2), \dots, x(l)) & = c_l \cr
}$$
hold.

\Pf Fix the $c_i$, and set $c_0 = 1$ (since $Q_0 =1$, this is consistent
with the equations); we proceed by induction. Since
$Q_1 (a) = a$, we set
$x(1) = c_1$. Assume $x(1), x(2), \dots, x(k)$ ($k < l$) have been found
so that  $Q_i (x(1), \dots, x(i)) = c_i$ for $1 \leq i \leq k$. The
$k+1$st equation can be rewritten as 
$$
x(k+1) Q_k(x(1), \dots, x(k)) + Q_{k-1} (x(1), \dots, x(k-1)) = c_{k+1};
$$
this simplifies to $x(k+1) c_k + c_{k-1} = c_{k+1}$. Since $k < l$,
$c_{k} $ is invertible, and thus we can set $x(k+1) = (c_{k+1} -
c_{k-1})c_k^{-1}$, completing the induction.
\qed

\comment
We now define  conditions for  a Banach algebra or more
generally a normed (or merely metrized) unital ring, $R$,to satisfy. Let
$U_n (B)$ be the set of unimodalar (left) $n$-tuples; that is, 
$$U_n (B) = \Set{(r_1, r_2, \dots, r_n)}{\sum Br_i = B}.$$

Let $v = (a(i))$ a tuple (of some specified size) of elements of $R$, and
abbreviate $Q_k (a(1),a(2), \dots , a(k))$ to $Q_k(v)$ (if $v$ is a
$k$-tuple).If $v$ is a $k$tuple, let $S_{k} (v)$ (or
$S(v)$ if no ambiguity will result) be the $(k-1)$tuple   obtained by
deleting the last (rightmost) entry. For $n \geq 1$ and a normed algebra
$R$, we say that it has $n$ in its topological Q-stable range ($\tQsr$)
if the set 
$$\eqalign{
&\Set {(Q_{n} (v),Q_{n-1}(S_n(v))}{v \in R^{n}} 
\text{is contained in
the closure of }\cr
&\Set {r(Q_{n-1} (S_nw),Q_{n}(w))}{w \in R^{n} \text{ and } r \in \GL
R}. 
}$$
That is, each pair of the form $(Q_{n}(v),Q_{n-1} (Sv))$ can be written
as a limit of a sequence of ordered pairs $(r_m Q_{n-1} (Sw_m), r_m
Q_{n}(w_m))_{m\to
\infty}$ where the $ r_m$ are invertible and the $w_m$ are  elements of
$R^{n}$. (Neither  sequence $(r_m)$ nor $(w_m)$ is required to converge.)
The important point is that the order is reversed. 

When $n=1$, the condition is simply that for all $a$ in $R$, $(a,1)$ lies
in the closure of 
$$\Set{(r,rb)}{r\in \GL R \text{ and } b\in R} .
$$
 This
forces the invertibles to be dense (project onto the first coordinate),
which condition is tsr$(R) = 1$ ($1$ is in the topological stable
range). On the other hand, if the invertibles are dense, then we may find
invertibles $r_j \to a$ and set $b_j = r_j^{-1}$, so that $(r_j, r_j b_j)
\to (a,1)$. Hence $\tQsr(R) = 1$ is equivalent to tsr$(R) = 1$.

It is possible (and desirable) that tsr$(R) \leq n$ implies $\tQsr(R)
\leq n$; however, I was not able to prove this.

\Lem Proposition . Suppose that $R$ is a normed algebra satisfying right
$\tQsr (R) \leq n$ and $(\ddag)$. For all
$n+1$-tuples of elements of
$R$,
$v= (a(1), a(2), \dots, a(n+1))$, there exists  $w = (a(n+2), \dots
a(2n+1))$ \st $Q_{2n+1}(v,w)$ is invertible.

\Pf Abbreviate $Q_{n+1} (v) $and $Q_{n} (Sv) $ to $Q_{n+1}$ and $Q_{n}$
respectively. Let $x = (x(2n+1), \dots, x(n+1))$ be an $n$tuple of
elements of $R$ whose values are to be determined. We have that
$Q_{2n+1}(v,x) = Q_n(x)\Op Q_{n+1} (v) + Q_{n-1}\Op(Sx)Q_n(v)$.

By the relations among the $P$s and $Q$s, we have that 
$$1 = (-1)^n P_n\Op
(n) (Sv) Q_{n+1} (v) + (-1)^{n+1} P_{n+1}\Op (v) Q_n(v).
$$ Hence, with $c
= P_n\Op(-Sv)$ and $d = P_{n+1}\Op(-v)$ (noting that $P_k\Op (-a) =
(-1)^{k-1}P_k (a)$ for any $k$tuple $a$), we have that $-1= cQ_{n+1} (v)+
d Q_n(Sv)$. Now $P_k (a) = Q_{k-1}(a(2),a(3),\dots, a(k))$, and  a
similar relation applies to the opposites. Hence, there exists $u$ in
$R^n$ \st $(d,c) = (Q_{n}\Op(u),Q_{n-1}\Op(Su))$; of course, $Q_k\Op(a) =
Q_k\Op (a\Op)$ (the last indicates $a$ written in reverse order). Hence
$\tQsr(R) \leq n$ means there existsequences $(r_j)$ of elements of $\GL
R$ and $(w_n)$ of elements of $R^n$ \st $r_j Q_{n-1}\Op (Sw_j) \to
Q_n\Op(u) = d$ and $r_j Q_{n}\Op (w_j) \to Q_{n-1}\Op(u)=c$. Hence $r_j
Q_n (w_j)\Op Q_{n+1}(v) + r_j Q_{n-1}\Op (Sw_j)Q_n(v) \to -1 $. Hence
for all sufficiently large $j$, the expression $r_j
Q_n (w_j)\Op Q_{n+1}(v) + r_j Q_{n-1}\Op (Sw_j)Q_n(v) = r_j
Q_{2n+1}(v,w_j)$ is invertible. Since every $r_j$ is invertible, so are
all but finitely many $Q_{2n+1}(v,w_j)$.
\qed

In analogy with the concrete case, we define the Order of an element of
$\Ag (R)$ to be the set of  numbers associated to  strings that represent
the element. Since $e_0$ corresponds to $\J$ (in the concrete case) and
$e_0 e_a$ corresponds to $T_a$, the following definition is natural.

If $g =  e_{a(k)} e_{a(k-1)} \cdots e_{a(1)} m_{r,s}$ where $a(j) = 0$
implies $j \in \brcs{1,k}$ (we can always reduce to this case, because of
the relation $e_a e_0 e_b = e_{a+b}$), then we define an element of
Order($g$) to be
$$\cases 
0 & \text {if $k = 0$}\\
1^{-}& \text{if $k=1 $ and $a(1) = 0$;}\\
k-2 & \text{if $a(k) = a(1) = 0$ and $k > 1$.}\\
k- \frac32 & \text{if $a(k) = 0$ but $a(1) \neq 0$;}\\
k^{-}& \text{if $a(k) \neq 0$ but $a(1) = 0$;}\\
k& \text{if $a(k),a(1) \neq 0$;}\\
\endcases
$$
\noindent Finally we define $\ord (g)$ to be the infimum of the elements
of Order($g$). It is an easy exercise to verify that this agrees with the
notion of order in the concrete case, and that   order is subadditive on
products. Let $S$ denote the operator which takes finite sequences of any
length and  deletes the last entry. 
\endcomment

Recall the definition of $\ord (g)$ for $g $ in $\pe$ (Section 4). 

\Lem Proposition \sevtwo. Suppose that $R$ is a ring. For a fixed positive
integer
$n$, the following are equivalent  {\par}
\item{(i)}for all $a = (a(1),a(2),\dots, a(n+1))$ in $R^{n+1}$, there
exists
$b = (b(1),b(2), \dots,b(n))$ in $R^{n}$ \st $Q_{2n+1}(a,b)$ is
invertible.{\par} 
\noindent \item{(ii)}  $\ord (g) \leq n + \frac 32$ for all $g$
in
$\petwo$. 

\Pf (i) implies (ii). We find a reduction rule $e_{a(n+1)}e_{a(n)}\dots
e_{a(1)} e_{a(0)}
\mapsto e_0 w m_{r,s}$ (for any $a(0), a(1), \dots,$ in $R$) where $w$ is a suitable
word in the 
$e$s. This will eventually reduce everything to an
element of   order $n + \slfrac 32$ or less.

Select $a(1)$, $a(2)$, \dots, $a(n+1)$ arbitrarily. By hypothesis, there
exist
$a(n+2),
\dots, a(2n+1)$
\st
$Q_{2n+1}(a(1), \dots, a(2n+1))$ is invertible. By \fouthr, there exist
$a(2n+2)$ and $a(2n+3)$ together with invertible $r$ and $s$ \st
$\P_{2n+3}(a(1), \dots, a(2n+3)) = r\oplus s$. It follows that
$e_{a(2n+3)}e_{a(2n+2)}\cdots e_{a(1)} = m_{r,s}$. This rewrites to 
$$\eqalign{
e_{a(n+1)}e_{a(n)}\cdots e_{a(1)}& =
\(e_{a(2n+3)}e_{a(2n+2)}\cdots e_{a(n+2)}\)^{-1}m_{r,s}   \cr
& = e_0e_{-a(n+2)} e_{-a(n+3)} \cdots e_{-a(2n+3)}e_0 m_{r,s};
\text{\ hence}\cr 
e_{a(n+1)}e_{a(n)}\cdots e_{a(1)}e_a & = e_0e_{-a(n+2)} e_{-a(n+3)}
\cdots e_{-a(2n+3)}e_0 m_{r,s}e_a \cr
& = e_0e_{-a(n+2)} e_{-a(n+3)}
\cdots e_{-a(2n+3)+  s^{-1}ar^{-1}} m_{r,s}. \cr
 }$$
If $a(n+1)$ and $a$ are both nonzero, the Order of the left side is
$n+2$; the Order of the right hand side is at most $n+ \frac 32$. A word of Order $(n+2)^{-}$ will have $a =
0$, and so the reduction also applies to this. It is
easy to see that if $w$ is any word of Order $n+2$ or more, by iteration
of this substitution, we will eventually obtain a word of  length at most
$n + \frac 32$. Note however, that  words of Order $n+ \frac 32$ are
unaffected by this reduction. 

\noindent (ii $\implies$ i). Pick $n+2$ elements, $a(0), a(1), \dots,
a(n+1)$ in $R$, and set $w = e_{-a(0)} e_{-a(1)}\dots e_{-a(n+1)}$, so
that $e_0 w e_0 e_{a(n+1)}\dots e_{a(0)} = 1$. 
Set $a = (a(0),a(1),\dots, a(n))$ (not with $a(n+1)$ at the end). We
will find $b $ in $R^n$ \st $Q_{2n+1}(a,b)$ is
invertible. By (ii), 
$\ord (w) \leq n+\frac 32$. 

First, suppose $\ord (w) = n+\frac32$. Then
there exist $b(i)$  ($i= 0, \dots, n+1$), and invertible $r$ and $s$ \st
$w = e_0 e_{b(n+1)}
\dots e_{b(0)} m_{r,s}$. Substituting, we obtain $1 = e_{b_{n+1}} \cdots
e_{b(0)} m_{r,s} e_0 e_{a(n+1)}\dots e_{a(0)}$. Moving the $m_{r,s}$ to
the left of the long expression yields
$$\eqalign{
m_{r',s'}  &=  e_{b'(n+1)} \cdots
e_{b'(0)}  e_0 e_{a(n+1)}\dots e_{a(0)} \cr
& =  e_{b'(n+1)} \cdots e_{b'(1)}
e_{b'(0)+a(n+1)} e_{a(n)}\dots e_{a(0)}.
}$$
Set   $b' = (b'(0)+a(n+1), b'(1), \dots,
b'(n+1))$. The  right side is represented by
$\P_{2n+3} (a,b')$, the left by $\(\smallmatrix r' & 0 \\ 0 &
(s')^{-1}\endsmallmatrix\) $. Hence $Q_{2n+3} (a,b')$ is invertible and
$Q_{2n+2} (a,Sb') = 0$, from which it follows that $Q_{2n+1}(a,S^2 b') =
Q_{2n+3} (a,b')$ and is thus invertible. Set $b = S^2 b' = (b'(0)+a(n+1), b'(1), \dots
b'(n-1)))$.

If $\ord (w) = n+1$, then either $w$ can be written the form 
\item{($\alpha$)} $w
= e_{b(n+1)} \cdots e_{b(1)} m_{r,s}$ or 
\item{($\beta$)} $w = e_0 e_{b(n+1)}
\cdots e_{b(1)} e_0 m_{r,s}$. 

\noindent
 In the first case, after passing the
$m_{r,s}$ to the left, and then to the other side we obtain an expression
of the form 
$m_{r',s'} e_0 =  e_{b'(n+1)} \cdots e_{b'(1)+ a(n+1)} e_{a(n)} \cdots
e_{a(0)}$. With $b' = (b'(1)+a(n+1), b'(2),\dots, b'(n+1))$  and $a =
(a(0),a(1),\dots, a(n))$ the right side is represented by
$\P_{2n+2}(a,b')$, and the left by $\(\smallmatrix r' & 0 \\ 0 &
(s')^{-1}\endsmallmatrix\) \(\smallmatrix 0 & 1 \\ 1 &
0\endsmallmatrix\) = \(\smallmatrix 0 & r'  \\ 
(s')^{-1} & 0 \\ \endsmallmatrix\)$. Hence $Q_{2n+1} (a,Sb') = r'$ which
is invertible; set $b = Sb'$. 

If $(\beta)$ applies, the corresponding expression is $m_{r',s'} =
e_{b'(n+1)} \cdots e_{b'(1)}e_{a(n+1)}\cdots e_{a(1)}$, leading to the
matrix equation $\P_{2n+3}(a,b') = \(\smallmatrix r' & 0 \\ 0 &
(s')^{-1}\endsmallmatrix\)$. Thus $Q_{2n+3} (a,b') = (s')^{-1}$ is
invertible and $Q_{2n+2} (a,Sb') = 0$, so that $Q_{2n+1}(a,S^2b')$ is
invertible; set $b = S^2 b'$. 

If $\ord (w) = (n+1)^-$, then $w = e_{b(n+1)}\cdots e_{b(2)}e_0
m_{r,s}$, and the equation becomes $m_{r',s'} e_0 =  e_{b'(n+1)} \cdots
e_{b'(2)} e_{a(n+1)}  \cdots e_{a(0)}$, and thus $\(\smallmatrix 0 & r'  \\ 
(s')^{-1} & 0 \\ \endsmallmatrix\) = \P_{2n+2}(a,b')$ where $b' =
(a(n+1),b'(2),\dots, b'(n+1))$. Then $Q_{2n+1}(a,Sb') = r'$ is
invertible, so we may set $b = Sb'$. 

If $\ord (w) < (n+1)^-$ (that is, $\ord (w) \leq n + \frac12$), the
arguments are even simpler.
\qed

In view of this, we define property $\SR (n)$ for a ring 
$R$: $n-1 + \frac 32<\max _{g \in \pe} \ord (g) \leq n + \frac 32$. We will see (Proposition \sevsix) that $\SR (1)$ is the same as stable range one. Of course, $\SR (\infty)$ means the failure of all $\SR (n)$. 

Now we require some Fibonacci-like identities. We recall that $Q_0 =
1$. In the following, the indices in the arguments of the $Q_k\Op$ are
descending, while those of the $Q_k$ are ascending. 

\Lem Lemma \sevthr. For any ring $R$, integers $N > m >0$, and any selection of
$b(1), \dots, b(N)$, the following is true:
$$\eqalign{
Q_{N} (b(1),b(2), \dots, b(N)) &= Q_{m}\Op (b(N), \dots,
b(N+1-m))Q_{N-m}(b(1), \dots, b(N-m)) +\cr
& \qquad Q_{m-1}\Op (b(N), \dots,
b(N+2-m))Q_{N-m-1}(b(1), \dots, b(N-m-1)).\cr
}$$

\Pf We may  assume that $R$ is the free ring on $b(1), b(2),
\dots, b(N)$. Every
$Q_k$ (and
$Q_k\Op$) consists of
$f(k)$ distinct words---each with multiplicity one---in the arguments,
where
$f(k)$ is the $k$th Fibonacci number (beginning $f(0) = 1 = f(1)$). Hence
the number of words on the left side is $f(N)$, and the number of words
on the right is $f(m)f(N-m) + f(m-1)f(N-m-1)$---which by one of the
thousands of Fibonacci identities is $f(N)$. Hence it suffices to show
that every word that appears on the right does so with multiplicity one
and appears in the left (thereby exhausting the left side). 

Now an arbitary $Q_k (b(1), \dots, b(k))$is the sum of all words
$b(i_1) b(i_2) b(i_3) \dots$of the following form:
\item{(i)}$i_1$ and the terminal $i$ have the same parity as $k$
\item{(ii)} $i_1 > i_2 > \dots$
\item{(iii)} $i_j - i_{j+1}$ is odd
\item{(iv)}the empty word, $1$, appears if and only if $k$ is even.

\noindent It is routine to establish this. Of course,  $Q_k\Op$
consists of the words of $Q_k$ written in reverse order. It follows from
the conditions that the length of each word has the same parity as $k$.

Select a word $w = b(i_1) b(i_2) b(i_3) \dotsb(s)$ that appears in the 
product $Q_{m}\Op (b(N), \dots,
b(N+1-m))Q_{N-m}(b(1), \dots, b(N-m))$; because the sets of variables in
the arguments are disjoint, there exists $u$\st $i_u \geq N+1-m >
i_{u+1}$(if the latter exists; if it does not, the second factor is the
empty word).This yields a factorization $w = w_1 w_2$, and now we 
check that with $w_1$ appearing $Q_{m}\Op (b(N), \dots,
b(N+1-m))$ and $w_2$ appearing in $Q_{N-m}(b(1), \dots, b(N-m))$, the
various conditions that have to hold for a word to appear in $Q_{N}$ do.
The only one that presents any problem is that $i_u - i_{u+1}$must be
odd---but this follows from the facts that $u$ and $m$ have the same
parity, that the length of $w_2$has the same parity as $N-m$, and that
within each of $w_1$ and $w_2$, the difference  of successive indices is
odd. The first condition on the parity of the first entry follows from
the reverse  ordering on the arguments in the opposite $Q$.

The same considerations apply to the other product, 
$$Q_{m-1}\Op (b(N), \dots,
b(N+2-m))Q_{N-m-1}(b(1), \dots, b(N-m-1)).$$
 We need only show that no
word in the sum of the products appears more than once. Certainly, the
contributions from each of the summands is disjoint, since those of the
first are of the form $w_1 w_2$ with the length of $w_1$ of the same
parity as $m$, while those of the second have the length of $w_1$ of the
same parity as $m-1$. And disjointness of the sets of arguments yield
that no word in one of the products has multiplicity greater than one in
the individual product.
\qed

\Lem Proposition \sevsix. For a ring $R$, the following are equivalent.
{\par}\item{(i)} For all $x,y$ in $R$ there exists a positive integer $l$ together with
$b = (b(1),\dots,b(l))$ in $R^l$ \st both of $Q_{l+2}(x,y,b)$ and $Q_{l-1}(b(2),\dots,b(l))$ are invertible;
\item{(ii)} $R$ has stable range $1$;
\item{(iii)} for all $x,y$ in $R$ there exists an element $c$ of $R$ \st $Q_{3}(x,y,c)$ is invertible.{\par}
\noindent In particular, $\ord (g) \leq 2\slfrac12$ for all $g$ in $\pe$ if and only if $R$ has stable range $1$.

\Pf (i) implies (ii). Suppose $Ra +Rc = R $; then there exist $s$ and $t$ in $R$ \st $sa - tc = -1$; then $Q_2(a,s) = 1 + sa = tc$. By (i), there exists an $l$-tuple $b$ so that $Q_{l+2} (a,s,b)$ and $Q_{l-1}(b(2),\dots,b(l)) = Q_{l-1}\Op (b(l),b(l-1),\dots, b(1))$ are invertible. By the previous lemma with $ N = l+2$ and $m = l$, we can write  
$$\eqalign{
Q_{l+2} (a,s,b) &= Q_{l-1}\Op(b(l),b(l-1),\dots,b(1))Q_1(a,s,\dots) + Q_l\Op(b(l),\dots,b(1)))Q_2(a,s,\dots)\cr
& = Q_{l-1}\Op(b(l),b(l-1),\dots,b(1))a
+ Q_l\Op(b(l),\dots,b(1))tc\cr
&= ua + vc  \cr
}$$
where $u$ is invertible; hence $a + u^{-1}vc = u^{-1} Q_{l+2}(a,s,b)$ is invertible, yielding stable range $1$.

\noindent (ii) implies (iii). Since $Ra + R(1+ba) = R$, and $R$ has stable range $1$, there exists $c$ in $R$ \st $a + c(1+ba)$ is invertible. But this expression is  $Q_3(a,b,c)$.

\noindent (iii) implies (i). Set $l = 1$ and observe that $Q_0$ is identically $1$.

The final assertion follows immediately from \sevtwo\ and (ii) if and only if (iii).
\qed

In part (i), the $l$ is permitted to vary with the choice of $x$ and $y$. Note how close condition \sevtwo(i) is to implying  stable range  $1$: if the $b$ can be chosen so that additionally $Q_{2n-1}(a(3),\dots, a(n+1),b)$ is invertible (but we do not require that the $n$ be uniform in $b$, that is, it can vary, depending on the first two coordinates of $a$), then stable range 1 applies.

If $R = \Z$, then $2$ is in the stable range of $R$, but $\Z$ fails all $\SR(n)$, as there is no bound on the number of transpositions  required to factorize  $\(\smallmatrix 0 & 1 \\ 1 & 1\\ \endsmallmatrix\)^k $ as elementary matrices (as $k$ varies). 

\comment
A condition sufficient to guarantee that (i) above hold, is the following
($n$ is fixed).

\itemitem{} for all $a$ in $R^{n}$ there exist $b$ in $R^{n}$, $c$ in
$R$, and invertible $r$ in $R$ \st 
$$
(P_{n+1}(a,c), Q_{n+1}(a,c)) = r(Q_{n-1}(b),P_{n-1}(b)).
$$

\noindent (In the topological setting, it would be enough to verify that the set consisting of
the possible right side is dense in the set consisting of the possible left side.)
The difficulty in attempting to verify this (even for $n=2$) is that the
$P$s and $Q$s are reversed. However, it does (vaguely) resemble a stable
range condition.

Results on simplicity of $\Ag_2 (R)$ are still in their nascent stage. We
require that $R$ be simple and generated as a ring by the commutator
subgroup and the central elements---these are not unreasonable
conditions. However, the proof at the moment requires that the order of
every element of $\Ag (R)$ be less than $3$, as occurs if $R$ has
stable range $1$.
The difficulty resides with products of three (or more) $e$s, as we
discuss afterwards.
\endcomment

Now we consider simplicity of  $\petwo$. Let $\Nn (g)$ denote the smallest normal subgroup of $\pe$
containing $g$. First a couple of elementary results.

\Lem Lemma \sevfou. Suppose that $g$ is an element of $\pe$.
{\par}
\itemitem{(a)} If for some integer $k > 0$, $\Ord (g)$ contains $k +
\frac 12$, then
$g$ is conjugate via an element of $ \petwo$ to  $g'$ for which
$\Ord (g')$ contains
$l$ for some
$l \leq k$.{\par} 
\itemitem{(b)} If $\Ord (g)$ contains $k^{-}$ for some $k \geq 3$,
then
$g$ is conjugate via an element of $\petwo$ to $g'$ with $k-2$ in $\Ord
(g')$; if
$\Ord (g)$ contains $2^{-}$, then $g$ is conjugate to $g'$ having
$\frac 12$ in
$\Ord (g')$.

\Pf (a) We note that $g = e_0 e_{a(1)} \cdot e_{a(2)}\dots
e_{a(k+1)}m_{r,s}$ is conjugate to $g':= e_{a(2)}\dots
e_{a(k+1)}m_{r,s}\cdot e_0 e_{a(1)}$. This expands to $ e_{a(2)}\dots
e_{a(k+1)} e_0 e_{ra(1)s}m_{r,s}$, which equals $e_{a(2)}\dots
e_{a(k+1)+ ra(1)s} m_{r,s}$, and thus $\Ord (g')$ contains $k$.

\noindent (b) Write $g = e_{a(k)}  e_{a(k-1)} \cdot e_{a(k-2)}\dots e_{a(2)} e_0
m_{r,s}$; this is conjugate  to $g' = e_{a(k-2)}\dots e_{a(2)} e_0
m_{r,s}e_{a(k)}  e_{a(k-1)}$, which equals $e_{a(k-2)}\dots e_{a(2) +
s^{-1}a(k)r^{-1}}   e_{a(k-1)} m_{r,s}$. If $k=2$, then $g = e_a \cdot
e_0 m_{r,s}$, which is conjugate to $e_0 m_{r,s} e_a = e_0
e_{s^{-1}ar^{-1}} m_{s^{-1},r^{-1}}$.
\qed

\Lem Proposition \sevfiv. Suppose that  $R$ is a simple  ring that  is generated as a ring by its invertibles, and its centre has at least four elements. If $g$ is an element $\pe$ with $\ord (g) \leq 2\slfrac12$, then 
 $\petwo \subseteq \Nn (g)$.

\Pf  The centre of a simple ring is a field, $K$, so if $k$ in $K$ is not $0,\pm1$, then $k(k^2-1)$ is not zero, so must be invertible.

 By  \sevfou, we may assume
that $g$ has Order zero, $\frac12$, $1^{-}$, $1$, or $2$.

\noindent  (a) {\it $g = e_0 e_a$ for some nonzero $a$ in $R$.}
Conjugating with $m_{r_i,s_i}$, we obtain $g_i = e_0 e_{r_i a s_i}$ in
$\Nn (g)$, and so the product $h=  g_1 g_2 \dots g_l$ belongs to $\Nn
(g)$; induction on $e_a e_0 e_b = e_{a+b}$ yields that $h = e_0 e_{x}$
where $x = \sum r_i a s_i$. Since $R$ is simple, every element of $R$ can
be expressed as a finite sum, $\sum t_j a u_j$
for some elements $t_j$, $u_j$ in $R$; since $R$ is generated by its
invertibles, each $t_j$ and $u_j$ is a sum of  invertibles, hence an
arbitrary element of $R$ can be put in the form $\sum r_i a s_i$ with
$r_i$ and $s_i$ invertible. Hence $\Set{e_0 e_b}{b \in R} \subseteq \Nn
(g)$. Obviously this implies that $\Set{e_b e_0}{b \in R} \subseteq \Nn
(g)$  as well.

Since $\petwo $ is generated by  $\Set {e_0 e_b, e_c e_0}{b,c \in R}$,
we have that $\petwo \subseteq \Nn (g)$.

\noindent (b) {\it $\ord (g) = 0;$ that is, $g = m_{r,s}$ for some
invertible
$r,s$ in $R$ \st
$m_{r,s}$ is not the identity.} Conjugate $g$ with $e_0 e_a$, to obtain
$e_0 e_{a - ras} m_{r,s}$; multiplying by $m_{r,s}^{-1}$ on the right, we
see that $e_0 e_{a-ras} \in \Nn (g)$ for all choices of $a$, and since
$m_{r,s}$ is not the identity, there exists $a$ \st $a \neq ras$, so that
$e_0 e_{a-ras}$ is not trivial. By (a), $\petwo  \subseteq \Nn (e_0
e_{a-ras}) \subseteq \Nn (g)$.

\noindent (c) {$\ord (g) = \frac12$.} The cases not covered in (a) occur
when $g = e_0 e_a m_{x,y}$ where $a $ is not zero, and $m_{x,y}$ is not
the identity. Conjugate $g$ with $e_0 e_b$; the outcome is $e_0
e_{a+b-xby} m_{x,y}$; multiplying on the right by the inverse of $g$, we
obtain $e_0 e_{b-xby}$ in $\Nn (g)$. Since $m_{x,y}$ is not the identity,
we may find $b$ \st $b \neq xby$, and therefore (a) applies.

\noindent (d)  {\it $g = e_0$.} Conjugate with $m_{a^{-1},a^{-1}}$ for
any invertible $a$ in $R$; the outcome is $e_0 m_{a,a}$, so that
$m_{a,a}$ belongs to $\Nn (g)$, and now (b) applies.

\noindent (e) {\it $\ord (g) = 1^{-}$.} Excluding the cases covered in
(d), we can write $g = e_0 m_{x,y}$ where $m_{x,y}$ is not the identity.
Conjugate with $m_{r,s}$ (for any pair of invertible elements of $R$);
the result is $e_0 m_{s^{-1},r^{-1}} m_{x,y} m_{r^{-1},s^{-1}}$ in $\Nn
(g)$. Multiply by $g^{-1} = m_{x^{-1},y^{-1}} e_0$ from the left; we
obtain $g_1:= m_{x^{-1}s^{-1} x r^{-1}, s^{-1}yr^{-1}y^{-1}}$ in $\Nn
(g)$. If $g_1$ is not the identity, case (b) applies; so we have to show
that there is a choice for $r,s$ so that $g_1$ is not trivial.

If  $g_1$ were the identity for a particular choice of $r$ and $s$, then
  there would exist central invertible $\lambda$
(depending on $r$ and $s$ in $R$) \st $x^{-1}s^{-1} x r^{-1} = \lambda$
and $ s^{-1}yr^{-1}y^{-1} = \lambda^{-1}$ (since the $m$s are in the image of $\GL (R) \times \GL
(R\Op)/\brcs{(z,z^{-1})}$). Multiplying them, we obtain
$$
x^{-1}s^{-1}x r^{-1} = yry^{-1} s.
$$
If $g_1$ were the identity for every choice of invertible $r$ and $s$,
then the displayed equation  holds for every such choice. If $s =1$, then
we deduce $r^{-1} = yry^{-1} $ for all invertible $r$ in $R$. If $s =2$
(and $2$ is invertible), then $r^{-1} = 4 yry^{-1}$, which is now
impossible.

\noindent (f) {\it $\ord (g) = 1$.}  We may write $g = e_a m_{x,y}$ with
$a
\neq 0$. There exists (since the invertibles
generate $R$ as a ring)  invertible $r$ \st $ra \neq a$. Form $h= m_{r,1}
g m_{r,1}^{-1}$, which expands to $e_{r^{-1}a}m_{u,v}$ for some invertible
$u$ and $v$. Then
$h^{-1}g = m_{u,v}^{-1} e_0 e_{-r^{-1}a} e_0 e_a m_{x,y}$, which
simplifies to $m_{u,v}^{-1} e_0 e_{a-r^{-1}a} m_{x,y}$. Passing the left
$m$ term gives us an element of the form $e_0 e_b m_{x',y'}$, where $b$
is not zero. This gives an element of Order $\frac12$ in $\Nn (g)$, which
now reduces to an earlier case.

\noindent (g) {\it $\ord (g) = 2$.} This is a slightly more complicated
version of the argument for $\ord(g)= 1$. Set
$h_0 = g = e_a e_b m_{t,u}$. Neither
$a$ nor $b$ can be zero as the order of $g$ is two. Let
$\brcs{r_i,s_i}_{i=1}^m$be  invertibles \st $\sum s_i^{-1}a
r_i^{-1} = a-1$. For each $i$, define $g_i$ to be $m_{r_i,s_i}e_a e_b
m_{t,u}m_{r_i,s_i}^{-1}$; this simplifies to $e_{s_i^{-1}a
r_i^{-1}}e_{r_ib s_i}m_{r_i t r_i^{-1},s_i^{-1}u s_i}$. Then (for this
computation, $\sim$ denotes conjugate to)
$$\eqalign{
g_1^{-1}h_0 & = m_{r_i t r_i^{-1},s_i^{-1}u s_i}^{-1} e_0 e_{-r_1 b s_1}
e_{-s_1^{-1}a r_1^{-1}} e_0 e_a e_b m_{t,u}\cr
& = m_{r_i t r_i^{-1},s_i^{-1}u s_i}^{-1} e_0 e_{-r_1 b s_1} \cdot
e_{a-s_1^{-1}a r_1^{-1}}  e_b m_{t,u}\cr
& \sim e_{a-s_1^{-1}a r_1^{-1}}  e_b m_{t,u} m_{r_i t r_i^{-1},s_i^{-1}u
s_i}^{-1} e_0 e_{-r_1 b s_1} \cr
& = e_{a-s_1^{-1}a r_1^{-1}} e_{b_1}m_{t_1, u_1}:= h_1.
 }$$

Now a similar computation, this time for $g_2^{-1}h_1$, together with the
conjugation which puts the first three terms at  the end yields $h_2 $ of
the form $e_{a - s_1^{-1}a r_1^{-1}- s_2^{-1} a r_2^{-1}} e_{b_2}m_{t_2,
u_2}$. Inductively, we obtain $h_m = e_x e_{b_m} m_{t_2,
u_2}$ where $x = a - \sum s_i^{-1} a r_i^{-1}= 1$, and of course each
$h_j$ is in  $\Nn (g)$. Thus $\Nn (g)$ contains an element of the form
$e_1 e_y m_{w,z}$  for appropriate $w,y,z$.

Now the identity $e_{-1} e_1 e_{-1} = m_{-1,1}$ yields that for any $y$,
$e_1 e_y = e_0 e_1 e_{-1-y}$, so that $h_m = e_0 e_1 e_{-1-y}m_{-w,z}$,
which has Order $\frac 32 < 2$, and we reduce to Order 1 by Lemma \sevfou.
\qed

\comment
If we attempt to extend the arguments of this result, say to $\Ord (g)
= 3$ or $4$ (which would  be required when either there is no upper
bound on the orders of elements, or when the upper bound is $3\frac12$
or $4\frac12$ respectively), something curious happens. Write $g = e_a
e_b \cdot e_c m_{x,y}$.   This is conjugate to $ e_c m_{x,y}e_a e_b=
e_c e_{y^{-1}ax^{-1}} e_{xby}m_{r,s}$; now conjugate this with
$m_{-1,1}$, and then with $e_0$ to obtain $h = e_0e_{-c}
e_{-y^{-1}ax^{-1}} e_{-xby}m_{x,y}e_0$. Now $g$ is
 conjugate to $k = m_{x,y}e_a e_b e_c$, and 
$$\eqalign{
kh &= m_{x,y}e_a e_b e_c e_0 e_{-c} e_{-y^{-1}ax^{-1}}
e_{-xby}m_{x,y}e_0 \cr
&  = m_{x,y}e_a e_b e_0 e_{-y^{-1}ax^{-1}}
e_{-xby}m_{x,y}e_0\cr
& = m_{x,y}e_a e_{b-y^{-1}ax^{-1}}e_{-xby}m_{x,y}e_0 \cr
& \sim e_a e_{b-y^{-1}ax^{-1}}e_{-xby} e_0 m_{y^{-1}x,yx^{-1}}.
}$$
A final conjugation and rewriting will put the $e_0$ in front, and we
would obtain an element of Order $2\frac12$, which subsequently reduces
to Order 2 or less---{\it except\/} if $a = xby$ and $y = x$, when the
expression on the right is the identity (i.e., $h = k^{-1}$). By a cyclic
permutation on $a$, $b$, $c$, and conjugation with suitable
$m_{t,u}$, we obtain a reduction to the case that $g = e_a^3$.

If $a=0$ or $a$ is invertible, we are done. If not, we consider the orbit
of $g$ under elements of the form $m_{r,s}$; we obtain elements of the
form $e_{s^{-1}a r^{-1}}e_{ras} e_{s^{-1}a r^{-1}}
m_{(rs)^{-1},(rs)^{-1}}$ in $\Nn (g)$. It is not clear how to proceed
further, except when $R$ happens to be von Neumann regular, when it is
possible to reduce to  Order two  or less (from Order three---however, 
it is not clear how to reduce from Order four to Order three even when
$R$ is von Neumann regular).

A similar dead end arises when $g = e_a e_b e_c e_d$; the process above
reduces to the form $e_be_{c-a} e_{-b} e_{a-c}$ in place of $e_a^3$.

\endcomment  

Now we recall the intersection of translates conditions from Section 1. We only require \gui 3, which asserts that for all $a,b$ in $R$, there exists invertible $u$ \st both $u+ b $ and $u+c$ are also invertible. For examples (and non-examples), see the two Appendices. 

\Lem Proposition \sevele. Suppose that $R$ is a simple ring whose centre contains at least four elements, and in addition, $R$ is generated (as a ring) by its invertibles. If $R$ either has $1$ in the stable range, or satisfies \gui 3, then $\petwo$ is contained in every normal subgroup of $\pe$. If additionally, $R$ is generated additively by $\gl (1,R)'$ and the central invertibles, then $\petwo$ is simple.

\Rmk In the presence of \gui 3, the assumption that $R$
 is generated by the invertibles is redundant. 
 
 \Pf If $R$ has $1$ in the stable range, then $\ord (g) \leq 2\slfrac 12$ for all $g$ in $\pe$ (Proposition \sevsix); by Proposition \sevfiv, $\petwo \subseteq \Nn (g)$, thus the result.
 
Now suppose that  $R$ satisfies \gui 3. It is enough to show that every $g$  can be conjugated to an element has $\ord$ value at most $2\slfrac12$. To do this, we exploit the reduction in the number of $j$s afforded by 
$$
j t_a j = t_{a^{-1}} j t_{-a} \m a,{a^{-1}}
$$
when $a$ is invertible (see the discussion of reductions after Proposition \fouone). 

Let $g $ be an element of $\petwo$, and suppose $\ord (g) \in \brcs{k^{-1},k, k + \slfrac12}$ for some positive integer $k$. If $k \leq 2$, we are done; so assume $k \geq 3$. By conjugating, we may assume $\ord (g) = k$ (Lemma \sevfou). Thus $g = e_{a(1)} \cdots e_{a(k)}\m r,s$ where the $a(i)$ are nonzero elements of $R$.

There exists an invertible element $a$
\st both $a + a(1)$ and $a + (a(1) + a(2))$ are invertible (that is, applying \gui {3} to the pair $a(1), a(1) + a(2)$). Let $b = a + a(1)$ and $c = a + a(1) + a(2)$, so that each of $a,b,c$ is invertible. Let $d = -b(1+ a(2) b^{-1})$; then $d = -bcb^{-1}$, so that $d$ is also invertible. 

We will conjugate $g$ with $jt_a j$, and show that the resulting element has $\ord$ value less than $k^{-}$, by showing it has at most $k-1$ $j$s in one of its factorizations. 

The following computation is not as horrible as it looks.
$$\eqalign{
(jt_aj) g (jt_a j)^{-1} & = j t_a j j t_{a(1)}j t_{a(2)}\cdots j t_{a(k)} \m r,s j t_a j\cr
& = (j t_{a + a(1) } j) t_{a(2)}\cdots j t_{a(k)} j\m {s^{-1}},r  t_a j \cr
& = t_{b^{-1}}j t_{b} \m b,{-b^{-1}} t_{a(2)}\cdots j t_{a(k)} j t_{-s^{-1}ar}\m {s^{-1}},r   j \cr
& = t_{b^{-1}}j t_{b} t_{-ba(2)b^{-1}}j \m -b,{b} t_{a(3)}\cdots j t_{a(k)} j t_{-s^{-1}ar} j\m *,* \cr
& = t_{b^{-1}}j t_{d}j  \m -b,b t_{a(3)}\cdots j t_{a(k)} j t_{-s^{-1}ar} j\m *,* \cr
& = t_{b^{-1}}j t_{d}j  \m -b,b j t_{a(3)}\cdots j t_{a(k)}  t_{-r^{-1}a^{-1}s} j t_{s^{-1}ar}\m *,* \cr
& = t_{b^{-1}}t_{d^{-1}}j t_{-d} \m d,{-d^{-1}}\m -b,b  jt_{a(3)}\cdots j t_{a(k)}  t_{-r^{-1}a^{-1}s} j t_{s^{-1}ar}\m *,* \cr
}$$
The various $\m *,*$ terms can be moved to the right without increasing the number of $j$ terms. The resulting expression has had two $j$ terms removed from the initial segment, and one $j$ term added to the end, resulting $k-1$ $j$s. We can continue this until the number is $2$ or less; at that point, the $\ord$ value of the conjugate is at most $2\slfrac12$, and we are done.

Now assume that $R$ is generated by $\GL {R}'$ and the invertible
central elements;  we show that the conjugations used in the constructions above can be
implemented by elements of $\petwo$, where necessary. For  (a) of Proposition \sevfiv, it
is enough that the elements $m_{r_i,s_i}$ belong to $\petwo$, and
sufficient is that $r_i s_i^{-1}$ belong to $\GL R'$, which can be
arranged by the current hypotheses. For (b), conjugation is already
implemented by an element of $\petwo$, as is the case in (c). Case (d)
only applies directly when $e_0$ belongs to $\petwo$, that is, when
$\peone = \petwo$; we observe that $m_{x,x}$ is already  $\S_4
(\dots)$, so belongs to $\petwo$. Case (e) uses only conjugations
with elements of the form $m_{x,y}$, and the arguments permit us to
select the $x$ and $y$ in $\GL R'$ when the current hypotheses apply.
The same applies in the $\Ord (g) = 1$ or $2$ cases.

Finally, the conjugations in this proof to reduce the $\ord$ value 
are implemented by elements of the form $j t_a j = e_a e_0$, which belong to $\petwo$.
\qed 

\comment
\SecT 7 Abelianized $\pe$

We examine the map $\Arrow \Psi_R;\GL R . \pe /\petwo $ induced by
$\Arrow \psi_R; \GL R . \pe $ in section 5.

\Lem Proposition \eigone.
The inclusion $\Arrow \psi_R; \GL R . \pe $ induces  an onto map
$$
\Arrow \Psi_R;\GL R/(Z(R)^2\cdot [\GL R,\GL R]) . \Ag(R)/\petwo  .
$$

\Pf We already know that the image of $(Z(R)^2\cdot [\GL R,\GL R])$ under
$\psi_R$ belongs to $\petwo $.
\qed

Hence necessary and sufficient for $\Psi_R$ to be an isomorphism of
(abelian) groups is that for all invertible $r$ in $R$, $m_{r,1}$
belonging to $\petwo $ entails that $r = \lambda^2 u$  where $\lambda $
is a central invertible, and $u$ is in $\GL R'$. Sufficient for this is
the following condition, which we discussed earlier.

\Lem Lemma \eigtwo. Suppose  for all $k$ and $a = (a(i))$ in $R^{1 \times 2k}$ \st
$Q_{2k} (a)$ is invertible, that $Q_{2k}^{-1} (a) Q_{2k}\Op (a)$
is in $Z(R)^2 \cdot \GL R'$. Then
$$
\Arrow \Psi_R;{\GL R \left/\(Z(R)^2\cdot [\GL R,\GL R]
\right.\)}. \pe / \petwo .
$$
is an isomorphism of abelian groups.

\Pf Suppose $s$ is an invertible element of $\petwo $; then there exist
invertible
$r_i$ and central $ \lambda$ in $R$  \st $(s\lambda)^{-1} = Q_{2k}(r_1,
\dots, r_{2k})$ and $\lambda = Q_{2k}\Op (r_i)$. Thus $s\lambda^2 =
Q_{2k}^{-1} Q_{2k}\Op$, so that $s $ belongs to $Z(R)^2\cdot 
[\GL R,\GL R]$. Hence the kernel of the composite map $\GL R \to \pe
\to
\pe/\petwo $ is contained in  $Z(R)^2\cdot [\GL R,\GL R]$. However,
this group is always contained in the kernel.
\qed

We know that $Q_{k}^{-1}Q_k\Op \oplus 1$ belongs to $\GL {2,R}'$. If $T$
is a unital ring with
$k$ in its  stable range, then the map $\GL {m,T} \to \GL {m+1,T}$ ($v
\mapsto v \oplus 1$) induces isomorphisms on the abelianizations if $m
\geq 2+ k/m$, and thus so does the map $\GL {m,T} \to \GL {2m,T}$ given
by $v\mapsto v \oplus \I$. Hence if
$R = \Mn m T$ for such an $m$, we deduce $Q_k^{-1}Q_{k}\Op$ is in $\GL
R'$ whenever $Q_k$ is invertible. Thus in these cases, we see that
$\Psi_R$ is an isomorphism. 

A Dedekind domain, $D$, has $2$ in the stable
range, so that if $R = \Mn _m D$ with $m \geq 3$, then $\pe /\petwo $
is isomorphic to $\GL R^{\text{ab}} /\text{Image} (D^*)^2$. In case
$\text{SK}_1(D) = 0$ (as occurs with number fields), $\GL R^{\text{ab}}$
is just $D^*$ (the unit group, via the determinant map), so $\pe /\Ag_2
(R) \iso D^*/(D^*)^{2m}$. If $D = \Z$, then $\pe /\petwo  \iso \Z_2$
regardless of $m > 2$.
\endcomment
\comment
\SecT 9 Banach algebras and friends

In this section, we show via a different style of proof that  if $R$ is a
Banach algebra or a suitable dense subalgebra with one in the stable
range, then
$\G(R)/\G_2 (R)$ is naturally isomorphic to $\GL R/Z(R)^2 \cdot[\GL R,\GL
R]$. 

For a Banach space $B$, we  recall the definition of the derivative of a
function
$\Arrow
\phi;U.B$ defined on an open subset, $U$, of $B$. For $x$ in $U$, $\phi$
has derivative
$D\phi (x)$ (a linear map from $B$ to itself) if there exists a bounded
linear map $\Arrow L \equiv L_x; B.B$ \st for all $\epsilon > 0$, there
exists $\delta$ \st for all $a$ in $B$ \st $\| a \| < \delta$, 
$$
\|\phi (x + a) - \phi(x) - L(a) \| < \epsilon \| a \|,
$$
and $D\phi (x) (a)  $ is defined to be $L_x (a)$.
As is well known, the chain rule holds for differentiable functions, that
is, $D (\phi \psi) (x)(a) = D\phi (\psi(x)) (D\psi (x)(a))$. If
$R= B$ is a Banach algebra with invertibles dense, the functions $\J$,
$T_b$ ($b$ in $R = B$) and $\M r,s$ (for $r$, $s$ in $\GL R$) are
differentiable (and defined on open sets), with respective derivatives,
$L_x = \M -x^{-1},{x^{-1}}, \M 1,1, \M r,s$ (of course, the last two do
not depend on the choice of $x$). Elements of $\G (R)$, being compositions
of such maps, are therefore differentiable (on the dense open subset in
their domains), and the formula for their derivative follows from the
chain rule. 

\Lem Proposition \ninone. Let $R$ be a Banach algebra whose set of
invertibles is  dense, and suppose that $r$ and $s$ are invertible in
$R$.{\par}
\item{(a)} If $\M r,s$ belongs to $\G_1 (R)$, then there exist invertible
$\lambda$ in the centre of $R$ and $g$ in $\GL {R}'$ \st $rs^{-1}$ is one
of $\pm \lambda^2 g$.{\par}
\item{(b)} If $\M r,s$ belongs to $\G_2 (R)$, then there exist invertible
$\lambda$ in the centre of $R$ and $g$ in $\GL {R}'$ \st $rs^{-1} =
\lambda^2 g$.{\par}
\item{(c)} The group $\G (R)/\G_2 (R)$ is naturally isomorphic to $\GL
R/(Z(R)^2 \cdot \GL R')$.

\Pf Let $a(1), \dots , a(k) $ be $k$ elements of $R$, and set $g_j = \J
T_{a(j)}$, so that $\S_i \equiv \S_i(a(1),a(2), \dots a(i)) = g_i g_{i-1}
\dots g_1$. We observe that for $x$ in the common domain of all the
relevant maps (and such exist, since the individual domains contain dense
open sets, and we are taking an intersection of finitely many), that
$Dg_1 (x) = \M -{(x + a(1))^{-1}},{(x + a(1))^{-1}}$. Inductively, we
obtain
$$\eqalign{
D \S_k (x)(b) &= D (g_k \S_{k-1})(x)(b) = D g_k (\S_{k-1}(x))(D
\S_{k-1}(x)(b)) \cr
& =  \M -{(\S_{k-1}(x) + a(k))^{-1}},{(\S_{k-1}(x) + a(k))^{-1}} D g_k
(\S_{k-1}(x))(D
\S_{k-1}(x)(b)) \cr
& =  \M -{\S_{k}(x)},{\S_{k}(x)} (D \S_{k-1}(x)(b) )\cr
&= \M -{\S_{k}(x)},{\S_{k}(x)} \M -{\S_{k-1}(x)},{\S_{k-1}(x)} \cdots \M
-{\S_{1}(x)},{\S_{1}(x)} (b) .
}$$
Set $v_i (x) = \S_i (x)$, and abbreviate this to $v_i$; by construction,
these are invertible---the $x$ is in the common domain of all the $g_i$
and everything else appearing here. We deduce
$D\S_k (x) = \M {(-1)^k v_k v_{k-1} \cdots v_1},{ v_1 v_{2} \cdots v_k}
$. If
$\M r,s$ belongs to $\G_1(R)$, then there exist a choice of $a(i)$ \st
$\M r,s = \S_k$ (we suppress the parentheses following $\S_k$---it should
be $\S_k (a(1), a(2), \dots, a(k))$, but this will cause confusion with
what we have denoted $\S_k (x)$, the effect of $\S_k$ applied to the
element $x$ in its domain). Pick $x$ in the common domain of all the
(finite set of) relevant transformations, and take derivatives; we deduce
that 
$$
\M r,s =   \M {(-1)^k v_k v_{k-1} \cdots v_1},{ v_1 v_{2} \cdots v_k} 
$$
(as linear functions on $R$). Hence there exists a central invertible
$\lambda$ (depending on $x$, as the $v_i$ do) \st 
$$
r = (-1)^k \lambda v_k v_{k-1} \cdots v_1 \qquad s = \lambda^{-1}  v_1
v_{2}
\cdots v_k.
$$
Thus $rs^{-1} = (-1)^k\lambda^2 v_k v_{k-1} \cdots v_1 v_k^{-1}
v_{k-1}^{-1}
\cdots v_1^{-1} $. The image of $g = v_k v_{k-1} \cdots v_1 v_k^{-1}
v_{k-1}^{-1} \cdots v_1^{-1}$ is trivial in any abelian group, so $g$ is
in
$\GL R'$.

If $\M r,s$ is in $\G_2$, then we can assume that the $k$ appearing above
is even, and the result follows. 

Part (c) is an immediate consequence of \eigtwo. 
\qed

\Lem Proposition \nintwo. Suppose that $R_0$ is a dense unital subring of
a Banach algebra $R$ \st $\GL {R_0}$ is relatively dense in $R_0$. If $r$
and $s$ are invertibles in $R_0$ \st $\M r,s$ belongs to $\G_2 (R)$, then
$rs^{-1} \in Z(R)^2 \cdot \GL {R_0}'$.

\Pf As $R$ is the completion of $R_0$, it follows that invertibles are
dense in $R$, so that $\G (R)$ and $\G_2 (R)$ are defined. By hypothesis,
there exists even $k$ and $a(i) $ in $R_0$ \st $\M r,s = \S_k (a(i))$ as
elements of $\G (R_0)$; density entails that this equality hold as
elements of $\G (R)$. By the proof of the preceding result, for all $x$
in the common domain in $R_0$ of all the operators involved---and now we
can assume that $x$ is in $R_0$, there exist the $v_i (x)$, and the fact
that
$x$ is in the common domain in $R_0$ means that each $v_i (x)$ not only
belongs to $R_0$  but is also invertible. Now the rest of the argument is
the same as in the preceding case. 
\qed

For example, this applies if $R = \Mn m C(X)$ where $X$ is compact and of
(topological) dimension one or less, and $R_0$ is a dense subalgebra \st
elements of
$R_0$ with determinant invertible in $C(X)$ are invertible in $R_0$.
However, these are subsumed by stability results.

\endcomment

\comment

Other stuff? relative orbits of $T_a$; that is, if $g$ is an element of
$\G(R)$ \st $g T_a g^{-1} = T_b$, then there exist invertibles $r$ and
$s$ \st $b = sar$ (easy; argument works in any Banach algebra with tsr
=1, and also for dense subalgebras for which invertibles are relatively
dense). 

If $\ord (g) \leq 2$ and $\M r,s = g \M t,u g^{-1}$, then there exists
$h$ in the subgroup generated by $\brcs{\M x,y} \cup \brcs{\J}$ \st 
$\M r,s = h \M t,u h^{-1}$; this fails if $\ord (g) = 2\frac12$ (first
part less easy, second part trivial).

Can occur that $T_a$ is conjugate to a nontrivial $\M r,s$, but not very
interesting.

Wanted nice description of relative orbit of $\J T_a \M r,s$; this seems
out of reach at the moment.

Conjugacies by elements of $\Ag (R)$ do not preserve fixed point sets
(even when $R = \Mn n \C$); what happens with $\Ag_2 (R)$ conjugacies?

\endcomment
\let\M = \Mn 

\SecT Appendix A. Intersections of translates of $\gl (1,R)$ 

In this paper, the following class of properties was used. Let $R$ be a ring. For a cardinal $\alpha >1$ (usually $\alpha =3$), we say {\it $R$ satisfies \gui {\alpha} } if for every subset of $S \subset R$ with $|S| = \alpha$, we have
$$
\bigcap_{s\in S} (\gl (1,R) - s) \neq \emptyset \tag {($\alpha$)}
$$
Here the notation $\gl (1,R) - s$ indicates the {\it translate\/} of $\gl(1,R)$ by $s$, that is $\Set{u-s}{u \in \gl(1,R)}$ (so if $s = 0$, the translate is simply $\gl(1,R)$ itself). [Warning: old-fashioned notation (which hopefully will never be used again) for $U \setminus A$ (for sets $U$ and $A$) is $U - A$.] We can obviously restrict to the case that one of the elements of $S$ is $0$, because of the following operation.

 For $u,v$ invertible in $R$ and $v_1, \dots, v_k$ elements of $R$, the transformation $R \to R$ given by $a \mapsto uav$ induces a bijection 
$$
\gl(1,R) \cap  \(\cap_{i=1}^k (\gl(1,R) - v_i \) \to 
\gl(1,R) \cap  \(\cap_{i=1}^k (\gl(1,R) - uv_iv \) .
$$

Throughout, an {\it invertible\/} element of a ring means two-sided invertible. As a noun, {\it invertible\/} is interchangeable with {\it unit,} but the former is preferred. 

So if $\alpha = n$, condition \gui {n} is equivalent to, for all sets of $n-1$ elements ${s_i}$ of $R$, there exists an invertible $u$ \st each of $u+s_i$ is invertible. (We can   allow some of the $s_i$ to be zero, but these can be discarded.)

If $n =2$, \gui {2} is just the condition that for all $s $ in $R$, there exists an invertible $u$ \st $u + s$ is also invertible, and this can be rephrased as, every element of $R$ is the sum of two invertibles (or equivalently, every element is a difference of invertibles). This condition was called {\it 2-good\/} in [V], where a number of interesting properties were proved. Some properties of 2-good rings extend easily to rings satisfying \gui {k} for general $k$; others extend, but not easily; and still others don't extend at all.

In the same paper, the property {\it $k$-good\/} was introduced: every element is a sum of $k$ invertibles. This is a weakening of $2$-good; on the other hand if $k > 2$, then \gui {k} is a {\it strengthening\/} of \gui{2} (and indeed, $\Mn_2 \Z$ satisfies \gui {2}, but not \gui 3). 

Other  related conditions in the literature include {\it twin-good\/}: for all $a$ in the ring, there exists an invertible $u$ 
\st both $u \pm a$ are invertible [SS]; this is weaker than \gui 3. A generalization of this is discussed in [KKN], {\it $n$-tuplet good\/}: for any $a$ and  $n$ central invertibles $u_i$, there exists an invertible $u$ \st all of $au_i + u$ are invertible. This is weaker than \gui {n+1}. 

Let $V = V(S)$ be the set $\Set{u-v}{u,v \in \gl(1,S)}$, the set of differences of units in $S$ (we could just as well have written $u+v$ of course, but then the {\it set of sums of units} is ambiguous). The condition that $S$ satisfies \gui {2} is equivalent to $S = V(S)$. 

We also (tentatively) define \gui 1: for all nonzero $a \in S$ \st $aS \neq S$,  $aS \cap V \neq \brcs{0}$.  A plausible but different alternative  definition for \gui{1} would be that the ring is generated (as a ring) by its invertibles (that is, every element is a sum of invertibles); another alternative would be to not define \gui {1} at all. As an exercise, if $R$ is a commutative ring all of whose proper images are finite, and having an invertible of infinite order, then $R$ satisfies this version of \gui 1. (This will be proved later in Appendix D, Lemma \Atwf, along with surprising negative results.)

The motivations for studying \gui {k} come from the discussion of intersections of translates of $\gl (1,R)$ in section 1, and crucially, the results in Section 6, about simplicity of $\petwo$ requiring either $1$ in the stable range or \gui 3. Because of this, the emphasis here is on \gui 3. 

Most results of the form, $R$ satisfies (or doesn't satisfy) \gui {n}, are very simple. But some, dealing with finite rings (reducing to matrix rings over finite fields) are not entirely trivial. 

First, some straightforward examples and properties.

\Lem Lemma \Aone. For fixed $\alpha$, the property \gui {\alpha} is preserved by 
\item{(i)} factor rings
\item{(ii)} products
\item{(iii)} direct limits

\Pf Obvious. \qed

Let $J(R)$ denote the Jacobson radical of the ring $R$.

\Lem Lemma \Atwo. If $R$ is a ring \st $R/J(R)$ satisfies \gui {\alpha}, then so does $R$. 

\Pf For $r$ in $R$, $r$ is invertible if and only if $r +J(R)$ is invertible in $R/J(R)$. \qed

\Lem Lemma \Athr. Let $R$ be a ring, and let $u,v$ be invertibles in $R$. For any  subset $S \subset R$, the map $x \mapsto uxv$ induces a bijection
$\cap_{s \in S} (\gl(1,R)-s) \to \cap_{s \in S} (\gl(1,R)-usv)$.

\Pf Obvious. \qed

So if $S$ contains zero, and one of its other elements is invertible, we can replace that element by the identity (of course the other, nonzero, elements are changed). 

\Lem Examples \Afou. 
\item{(i)} Matrix rings over a division ring with  infinite centre satisfy \gui {n} for all  integers $n$ [this will be improved; see Corollary \Asev]
\item{(ii)} (von Neumann) regular algebras over an uncountable field and with no uncountable direct sums of nonzero right ideals satisfy \gui {\aleph_0}
\item{(iia)} Stably finite simple regular rings that are algebras over an uncountable field satisfy \gui {\aleph_0}
\item{(iii)} Banach algebras satisfy a condition  stronger than  \gui {n} for all $n$
\item{(iiia)} If $S $ is a dense subring of a Banach algebra $R$ \st $\gl(1,S) = S \cap \gl(1,R)$, then $S$ satisfies a condition  stronger than  \gui {n} for all $n$.
\item{(iv)} Algebraic algebras over an infinite field satisfy \gui {n} for all  $n$.

\Pf (i) Set $R = \Mn_n D$, where $D$ is a division ring with infinite central subfield $K$. In $R$, an element is invertible if its (right) annihilator is zero. 
For $s $ in $R$, define 
$$
V_s: = \Set {\lambda \in K}{\lambda + s \text{ is not invertible}}.
$$
 Since  the sum of right modules 
$\sum_{\lambda \in V_s} (\lambda + s)^r $ is a direct sum (the same argument as in the standard argument that the number of eigenvalues of a matrix over a field is at most $n$), and since $R$ contains no nontrivial direct sums with more than $n$ summands, it follows that $|V_s| \leq n$. Hence the intersection of any finite set of $V_s^c$ (the complements) is nonempty. 

\noindent (ii) Similar argument as in (i), except that here $V_s$ is at most countable. 

\noindent (iia) Such rings admit a faithful rank function, hence cannot contain a nontrivial uncountable direct sum of right ideals. So (ii) applies.

\noindent (iii)  Let $T$ be any set bounded in norm, say by $c$. If $\lambda$ is a complex number \st $|\lambda | > c$, then $\lambda + c$ is invertible. 

\noindent (iiia) This is not entirely trivial, because we have not assumed $S$ is an algebra. However, the integers are contained in $S$, and the assumption on invertibles implies $S$ is an algebra over the rationals. Now the argument in (iii) (using the assumption about invertibles again) works with any sufficiently large rational number, for a bounded set $T$.

\noindent (iv) If $a$ in $R$ is algebraic over a central field $K$, then $a$ is invertible if and only if it is a nonzero divisor in the (unital) algebra it generates, $K[a]$. For $s \in R$, the finite dimensional algebra $K[s]$ admits a  (unital) faithful representation to some some matrix ring over $K$, and its image thus has only finitely many eigenvalues. In particular, for all but finitely many elements of $K$, $\lambda + s$ has invertible image, and thus is not a zero divisor in $K[s]$, hence is invertible therein. So $V_s$ (see (i)) is finite. 
\qed

\comment
\Lem Corollary. Let $R$ be an artinian $K$-algebra, where $K$ is infinite. Then $R$ satisfies \gui {n} for all integers $n$. 
\endcomment

The most important case for our purposes is \gui 3. Trivially, both finite fields $F_2$ and $F_3$  fail to satisfy \gui 3.

\Lem Example \Afiv. $S = \Mn_2 F_2$ does not satisfy \gui 3.

\Pf Set $b = \(\smallmatrix 1 & 0 \\ 0 & 0 \\ \endsmallmatrix \)$ and $c = \(\smallmatrix 0 & 1 \\ 0 & 0 \\ \endsmallmatrix \)$. Since the group of invertibles consists of only six matrices, it takes about a minute (less than the time to write this up) to verify that \gui {3} fails in this example. \qed

 It also follows that $\M_2 \Z$ fails to satisfy \gui 3.
 
 \comment
\Lem Lemma . Let $(X,\mu)$ be a probability space, and let $U,Y,Z$ be (measurable) subsets of $X$. Then 
$$\eqalign{
1- \mu (U \cap Y \cap Z) &\leq (1- \mu(U) + (1-\mu(Y) +(1- \mu (Z)) - \(\mu(U^c \cap Y^c) + \mu (U^c \cap Z^c)\right. \cr  
& \qquad+ \left.\mu (Y^c \cap Z^c) - \mu (U^c \cap Y^c \cap Z^c)\) \cr
& \leq (1- \mu(U) + (1-\mu(Y) +(1- \mu (Z) - \(\mu(U^c \cap Y^c) + \mu (U^c \cap Z^c)\) .\cr
}$$

\Pf Inclusion-exclusion says
$$
\mu (U^c \cup Y^c \cup Z^c) = \mu (U^c) + \mu (Y^c) + \mu(Z^c) -\(\mu(U^c \cap Y^c) + \mu (U^c \cap Z^c) + \mu (Y^c \cap Z^c)\) +\mu (U^c \cap Y^c \cap Z^c),
 $$
 and the statements follow. 
 \qed 
The  following  will be superseded by a more complete result (Proposition \Aele). But the proof is interesting.
\endcomment

\comment
\Lem Proposition \Asev. Let  $F = F_q$ be the finite field with $q$ elements. Then for all $n$, $\Mn_n F$ satisfies \gui {q-1}. In particular, if $q \geq 4$, then $\Mn_n F$ satisfies \gui {3} for all $n$. 

\Pf We show that $A: = |\gl(n,F)|/|\Mn_n F| > 1 - 1/(q-1)$, from which the result follows easily. Define $f_n(q) = \prod_{i=1}^n (1- q^{-i})$; this is  $A = \mu (\gl(n,F))$ where $\mu$ is the normalized counting measure on $\Mn_n F$. Let $f (q)= \prod _{i=1}^{\infty} (1- q^{-i})$, the  infinite product (convergence is routine to verify).
Then $q \mapsto f(q)$ is monotone increasing, and a basic approximation argument shows $f (q) > 1- 1/(q-1)$. Thus  $f_n (q) > f(q) > 1 - 1/(q-1)$ for all $n$ and $q \geq 4$. 

Let $b$ be an element of $R$; we observe that each of the  translates $\gl(1,R) - b$ has the same cardinality as $\gl(1,R) = \gl (n,F)$. For $b_1, b_2, \dots , b_{n-1}$ in $\Mn_n F$, the intersection $\gl(n,F) \cap \( \cap_{i=1}^{n-1}
\gl(n,F)- b_i)\)$ is an intersection of $n$ subsets of $\Mn_n F$, each having more than $1- 1/(q-1)$ of the elements of $\Mn_n F$. The intersection is nonempty by Lemma \Asix.\qed
\endcomment

Frequently, it turns out that proofs of results for ${2}$-good (property \gui 2) rings carry over to those satisfying \gui k, for some or all $k$. This occurs for the easier results above, e.g., \Afou, which appear in [V] for \gui 2. However, if we look at the argument in [V, Prop 8] to show that if $R$ satisfies \gui 2, then so do all matrix rings $\Mn_n R$, we see that it does not carry over to, for example, \gui 3. However, the proof of a slightly more general result, [WR, 3.6(i)], {\it can\/} be extended to general \gui k, with minor modification.

Let $e$ be an idempotent in a ring $R$. We say that $a$ in the corner ring $eRe$ is {\it relatively invertible\/} if there exists $\overline{a}$ in $eRe$ \st $a\overline{a} = e = \overline {a}a$. For $b$ in $(1-e)R(1-e)$, we use the notation $\underline{b}$ for the relative inverse of $b $ in $(1-e)R (1-e)$, that is, $b\underline{b} = 1-e = \underline {b}b$.

\Lem Proposition \Aeig.  Let $R$ be a ring, and suppose that $e$ is an idempotent \st both $eRe$ and $(1-e)R(1-e)$ satisfy \gui {k} for some integer $k\geq 2$ (that is, regarding $eRe$ and $(1-e)R(1-e)$ as unital rings with respective identity elements $e$ and $1-e$). Then $R$ satisfies \gui k. 

\Pf (After the proof of [WR, 3.6(i)]) Let $a(1), a(2), \dots, a(k-1)$ be $k$ elements of $R$. We will find invertible $u$ in $R$ \st $u + a(i)$ are all invertible. Using the Peirce decomposition, we write 
$$
a(i) = e a(i) e + ea(i) (1-e) + (1-e) a(i) e + (1-e)a(i) (1-e).
$$
By hypothesis, in the ring $eRe$, there exist $u = eue$, together with $\overline u$ \st $u \overline{u} = \overline {u} u = e$ and each of $ ea(i) e -u:= u_i$ are relatively invertible in $eRe$ (that is, there exist $\overline{u_i}$ in $eRe$ \st $\overline {u_i} u_i = e = u_i \overline{u_i}$).

Define $c_i = (1-e) a(i) (1-e) - (1-e) a(i)e \overline{u_i}ea(i) (1-e)$; this belongs to $(1-e) R (1-e)$. By hypothesis, there exist   $v$ and corresponding $\underline {v} $  in $(1-e) R(1-e)$ with $v \underline {v} = 1-e = v\underline {v}$, \st for all $i$,
$$
c_i - v:= v_i
$$ 
is relatively invertible, that is, $\underline {v_i} $ exist. 

We have 
$$\eqalign{
a(i) &= e(u + u_i) e + ea(i)(1-e) + (1-e)a(i) e + \(v + v_i + (1-e) a(i)e \overline{u_i} e a(i) (1-e) \) \cr
& = (u + v) + \(eu_i e  + v_i + (1-e) a(i)e \overline{u_i} e a(i) (1-e) \). \cr
}$$
The first parenthesized term on the second line, $u+v$, is clearly invertible in $R$ and independent of $i$. Call parenthesized term to its right, $A(i)$, so I don't have to repeatedly copy and paste it. It suffices to show $A(i)$ is invertible in $R$. 

To this end, consider
$$\eqalign{
\( e - (1-e)a(i) \overline{u_i} + 1-e \) &A(i) \(e - \overline{u_i}e a(i)(1-e) + 1-e\)\cr
& = \( u_i + ea(i)(1-e) + v_i\) \(e - \overline{u_i}e a(i)(1-e) + 1-e\)\cr
& = u_i + v_i \cr
}$$
Each $u_i + v_i $ is invertible in $R$; the  parenthesized terms on the first line are of the form $1 + n$ where $n^2 =0$, so are invertible. It follows that  $A(i)$ is invertible [and if we were talking about $2 \times 2$ matrices, this would express $A(i)$ as a product of three elementary matrices].
\qed

Variants on the argument are given in Appendix C, where they are used to shorten the proof of the main result there. 

So we deduce immediately the results corresponding to [V, Proposition 8],  for \gui k. This also yields a  proof that for $F_q$, the finite field with $q$ elements, $\Mn_n F_q$ satisfies \gui {q-1}, but we  do better than that in Appendix B. 

\comment
The following will be subsumed by a more general result, xxx, when $R$ is right or left artinian. 

\Lem Corollary. If $R$ is a finite ring, then $R$ satisfies \gui {3} iff it has no factors of the form  $F_2, M_2 F_2$, or $F_3$.

\Pf $R$ satisfies \gui {3} if and only if its factor by the Jacobson radical satisfies it. But the latter is a finite semisimple ring, and thus is a product of various size matrix rings over finite fields, and a product satisfies\gui {3} if and only if each of the constitutents does. \qed 

\endcomment

\Lem Corollary \Anin. Let $R$ be a ring. Let $k \geq 2$ be an integer. 
\item{(a)} If  $R$ satisfies \gui {k}, then so do all matrix rings, $\Mn_n R$.
\item{(b)} Let $W_k(R) = W_k$ be 
$
\Set{n \in \N}{\Mn_n R \text{ satisfies \gui k}}.
$
\noindent Then either $W_k$ is empty or is closed under addition.

If $W_k(R)$ is not empty,   it is thus a numerical semigroup (ignoring the non-participation of zero---the latter is superfluous anyway) if it contains a subset with greatest common divisor $1$.  For example, $W_3 (F_2) = \N \setminus \brcs{1,2}$ and $W_3 (F_3) = \N \setminus \brcs{1}$; this also equals $W_4 (F_3)$, but by an ad hoc  proof.

It is possible to write down a prosaic proof of \Anin(b) (that is, without using Proposition \Aeig) in the case that $R$ is right or left Hermite, but this has become moot.

\Lem Corollary \Atff. Let $R$ be a right or left artinian ring. Then $R$ satisfies \gui {n} for all $n$ if and only if it has no finite images.

\Pf Let $S$ be a finite ring; there is a trivial upper bound on $\Set{n}{\text{$S$ satisfies \gui n}}$---$S$ cannot satisfy \gui {|S| + 1} (a better bound will be given later). By Lemma \Aone(i), any $R$ with a finite factor cannot satisfy \gui {n} for all sufficiently large $n$.

Now assume $R$ has no finite images. We reduce to $\overline{R} = R/J(R)$ by Lemma \Atwo, and further reduce to $R = \Mn _k D$ where $D$ is a division ring, which is necessarily infinite. Thus $D$ trivially satisfies \gui {n} for all $n$, so by Proposition \Aeig, all matrix rings do. 
\qed

We use the standard notation, $F_q$, to denote the finite field with $q$ elements. We will show (in  Appendix B) that $\Mn _n F_q$ satisfies \gui {q} for all $n \geq 2$, and for $q \geq 3$ and $n$ even, it satisfies \gui  {q+1} (Propositions \Aele\ and \Bsev). In  Appendix C, we will show that $\Mn _n F_2$ satisfies \gui {3} for all $n \geq 3$.

\Lem Corollary \Atwe. Suppose that $R$ is a right or left artinian ring. Then $R$ satisfies \gui {3} if and only if it has no homomorphic images isomorphic to any one of $F_2$, $\M _2 F_2$, or $F_3$.

\Pf By Lemma \Atwo, we can factor out the Jacobson radical, so $\overline{R} \iso \prod_i \M_{n(i)} D_i$ where $D_i$ are division rings. If $D$ is  infinite, it thus  satisfies \gui {k} for all $k$, and from Corollary \Anin, so do all of its matrix rings. If $D_i$ is finite, then it is a finite field.  $\overline{R}$ satisfies \gui {3} if and only if all $\Mn_{n(i)}D_i$ do, and the result then follows from Corollary \Aele\  and Proposition \Bone. \qed

\comment
Of course, this last generalizes, e.g., if $R$ is right or left artinian, then $R$ satisfies \gui {k} for all $k \geq 2$ if and only if it has no finite images (it is easy to see that $\Mn_n F_q$ cannot satisfy \gui {n(q-1) +1}, which we will address later). 
\endcomment

In the other direction, $\Z$ and $\Mn_2 \Z$ obviously fail to satisfy \gui {3} (since they have factor rings that don't). And $\Z$ fails to satisfy \gui {2}, and even \gui 1. If we eliminate the offending primes $2$ and $3$, we might wonder whether $\Z[1/6] $ satisfies \gui {2}. The answer is no. On the other hand, even $\Z[1/2] $ satisfies \gui 1. 

It is true however, that if $R$ is a (commutative) principal ideal domain, then $\M_2 R$ satisfies \gui {2}. More generally, if $R$ is an elementary divisor ring, then $\Mn_n R$ satisfies \gui {2} for all $n > 1$. 

\Lem Proposition \Athi. ([V, Proposition 6]). Let $R$ an elementary divisor ring. Then for any $n \geq 2$, $\Mn_n R$ satisfies \gui 2.

We can (marginally) extend this. 
The relevant results we need about Hermite rings are [Ka; Theorem 3.5], asserting that for right Hermite rings,  if $A$ is in $\Mn_n S$ then there exists $U  \in \gl(n,S)$
 \st $AU$ is triangular---upper or lower, we can choose which we want, and [Ka; Theorem 3.6], that matrix rings over right Hermite are again right Hermite. 

A stronger condition, being an elementary divisor ring, is also defined in [Ka, p 465]. The definition is that for every rectangular matrix $A$ there exist invertible matrices of the appropriate sizes \st $PAQ$ is diagonal, and an additional divisibility condition holds (we invariably do not use the additional divisibility condition). 
For commutative domains, this is equivalent to (a) every finitely generated ideal is principal, and (b) for every triple $(a_1,a_2,a_3)$ \st $\sum a_i S = S$, there exist $x,y$ \st $xa S + (xb + yc)S = S$, [Ka, Theorem 5.2]. 

So we make a slightly weakened definition, a ring is {\it weakly divisorial\/} if for every rectangular matrix $A$ there exist invertible matrices of the appropriate sizes \st $PAQ$ is diagonal.

\Lem Lemma \Aftn. Suppose $S$ is a weakly divisorial ring. Then for every $n >1$, $\Mn_n S$ satisfies \gui 2.

\Pf Given $B \in \Mn_n S$, there exists  $U,V \in \gl(n,S)$
\st $UBV:= T$ is diagonal, and thus it suffices to find an invertible matrix $W$ \st $W + T$ is invertible. Let $W$ be the matrix that is given as the cyclic permutation matrix (ones in $(n,1), (1,2), \dots (n-1,n)$ positions), plus the matrix with $-d_{1,1}$ in the $(1,1)$ position and no other nonzero entries. It is an easy exercise to show that $W$ and $W + D$ are both invertible, by expressing each as a product of elementary matrices from the right, and from the left. 
\qed

Here is a conjecture which seems reasonable, but which I haven't been able to prove. It would yield at least that $\M_4 \Z$ satisfies \gui 3. 

\Lem Conjecture \Affn. Suppose that $S$ is an elementary divisor ring satisfying \gui 2. Then for all $n \geq 2$, $R = \Mn_n S$ satisfies \gui 3.

Now we deal with non-examples. 
Recall, for a ring $S$, the definition of $V(S)$, the set of  differences of units,
$$
V(S) = V = \Set{u-v}{u,v \in \gl (1,S)}. 
$$

\Lem Proposition \Asxn. [V, Proposition 10] Let $R$ be a ring  containing an $n$-generated  right ideal $J = \sum_{i=1}^n a_i R$ that is not $n-1$-generated. If $J \cap V(R) = \brcs{ 0}$, then $\Mn_n R$ fails \gui 2.

This applies to $R = \Z[x]$ (take $J$ to be the ideal generated by $\brcs{3^{n-i} x^i}_{i=1}^n$) and  $F[x,y]$, and for these rings for no $n$ does $\Mn_n R$ satisfy \gui {2} (let alone \gui 3). Along these lines, we can find small rings with similar properties.

\Lem Lemma
\Asvn. Let $S$ be a commutative domain with field of fractions $K$. Suppose there exist nonzero elements $a,b,c$ in $S$ with the following properties
\item{(i)} the ideal $J: = <a,b>$ (generated by $\{a,b\}$) is not principal;
\item{(ii)} $ c S \cap V = \brcs{0}$ and $J \cap V = \brcs{0}$.{\par}
\noindent Then $ \Mn_3 S$ does not satisfy \gui 3.
If additionally, for each pair $r,s$ of elements of $S$ \st $ar = bs$, each of $r$ and $s$ belongs to $J$, then $\Mn_4 S$ does not satisfy \gui 3. 

\Pf Define 
$$
B = \( \matrix 0 & a & b\\
0 & 0 & 0 \\ 
0 & 0 & 0 \\\endmatrix\) \qquad C = \( \matrix c & 0 & 0\\
0 & 0 & 0 \\ 
0 & 0 & 0 \\\endmatrix\).  
$$
Suppose there exists $U $ in $\gl(3,S)$ \st both $U+B$ and $U+C$ are invertible. 

Let $d_i$ be the determinant of the $2 \times 2$ minor obtained from $U^{(1,i)}$, that is, obtained from deleting the first row and $i$th column of $U$. Then $< d_1, d_2, d_3 >$ is $S$. 

We compute the determinants of $U+ B$ and $U+C$; these are elementary.
$$\eqalign{
\det (U + B) & = \det U - ad_2 + b d_3\cr
\det (U + C) & = \det U + c d_1 \cr
}$$
Thus $cd_1$ is a difference of units (belongs to $V$), and so is $-ad_2 + bd_3$. By hypothesis (ii), both must be zero, that is, $c d_1 = 0$ and $-ad_2 + b d_3 = 0$. The first yields that $d_1 = 0$, and the second, that $a d_2 = b d_3$. 

Now we claim that the ideal generated by $\brcs{d_1, d_2}$ is not principal: it is isomorphic to $J$, via the observation that $a < d_2,d_3> = d_3 < a,b>$. Hence $< d_1, d_2, d_3> = < d_2,d_3>$ cannot be $S$. 

Now suppose the supplementary hypothesis holds, and redefine $B$ and $C$ to be the size four matrices with first rows $\(\matrix a & b & 0 & 0 \\ \endmatrix \)$ and $\(\matrix   0 & 0 &a & b\\ \endmatrix \)$ respectively, with all their other rows consisting of zeros. Define $d_i$ to be the size three minors obtained by expanding along the top row of the candidate invertible $U$ (\st $U + B$ and $U+C$ are invertible). 

Then $\det (U + B) = \det U + a d_1 - bd_2$ and $\det (U+C) = \det U + ad_3 - bd_4$. Again by hypothesis (ii), $a d_1 = b d_2$ and $a d_3 = b d_4$. The supplemental hypothesis says that each $d_i$ belongs to $J$, so that $< d_1, d_2, d_3, d_4>$ cannot be improper; but this contradicts invertibility of $U$. 
\qed 

Now we give an example satisfying all the conditions. 

\Lem Example \Aegn. An order $S$ in a number field   \st $\M_3 S$ and $\M_4 S$ fail to satisfy \gui 3.

\Pf Let $S = \Z[3i]$, a subring of $\Z[i]$, the Gaussian integers. We can rewrite $S = \Z + 3\Z[i]$. Set $J = 3 \Z[i]$, which we view as an ideal of $S$. It cannot be principal (for otherwise its endomorphism would be $S$, but its $S$-module endomorphism ring is obviously $\Z[i]$, which is not isomorphic $S$).

We can write $J = 3S + 3iS$, and set $a = 3i$, and $b = c = 3$. The unit group of $\Z[i]$ consists only of powers of $i$, and thus the unit group of $S$ is simply $\brcs{\pm 1}$. Therefore $V(S) = \brcs{0,\pm 2}$. Since all elements of $J$ are of the form $3m + 3in$ ($m,n$ integers), it follows immediately that $J \cap V = \brcs{0}$. Thus $\M_3 S$ fails \gui 3. 

Next, suppose that $3r = 3i s$ for some $r,s $ in $S$. Writing $r = A + 3Bi$ and $s = C + 3iD$ for integers $A,B,C,D$, we have $3A =-9D$ and $3B = C$. Thus $3$ divides $A$ and $C$, so that $r,s \in J$. The supplemental hypothesis is now satisfied.
\qed

This example is not Dedekind, and it is likely there are no examples among number fields, as the unit groups are big (except for rings of integers in quadratic imaginary extensions).
See section C for interesting behaviour of functions fields.

Sometimes we can give an upper bound on the $k$ for which $R$ satisfies \gui k. Unfortunately, these tend to be crude, based on simple ideas. Recall that $F_q$ denotes the field with $q$ elements. 

\Lem Lemma \Antn. If $R = \Mn_n F_q$, then $R$ does not satisfy
\gui {(q-1)n+1}.

\Pf Let $e_i$ denote the $n \times n$ matrix whose only nonzero entry is a $1$ in the $(1,i)$ position, and consider the set with $(q-1)n$ matrices, $\Set{z e_i}{z \in F_q\setminus \brcs{0}; 1 \leq i \leq n}$. 

Suppose that $U$ is in $\gl(n,R)$ and $U + z e_i$ is invertible for all $(q-1)n $ choices of $z$ and $i$. Let $U^{1,i}$ be the square $n-1$ minor obtained by deleting the first row and $i$th column, and set $d_i =(-1)^{i+1} \det U^{1,i}$. 
We easily see that 
$$
\det (U + z e_i) = \det U + z d_i.
$$
Now fix $i$; if any two of these are equal, we obtain $d_i = 0$. Otherwise, we obtain $q-1$ distinct nonzero elements of $k$. Since $\det U \neq 0$, subtracting it from all of the others, one of the differences, $z d_i$, must be zero, and so $d_i = 0$ in any case. Since this is true for all $i$, $\det U = 0$, a contradiction.
\qed

For example, this means that $\Mn_n \Z_2$ fails \gui {n+1}, and this may even be sharp (all we know at the moment is that it satisfies \gui {3} if $n\geq 3$; this is proved in Appendix B). For $\Mn _2 F_3$, we know it satisfies \gui 3; using inclusion-exclusion, it is possible to show that it satisfies \gui 4. This is sharp, as it doesn't satisfy \gui 5{} by Lemma \Antn. 

On the other hand, $\Mn_n F_7$ satisfies \gui 7, a far cry from the upper bound obtained here, \gui {6n}. 

Along the same lines, we have a result applicable to $\Z$. 

\Lem Lemma \Atwh. Let $S$ be a ring, suppose there exists  $a$ in $S$ \st $aS \cap V(S) = \brcs{0}$. 
Then $\Mn_n S$ fails to satisfy \gui {n+1}. 

\Pf [Same idea as the proof of [V, Proposition 10].] Let $e_i$ be the standard $i$th basis element for  $S^{1\times n}$, and set $E_i$ to be the $n\times n$ matrix whose top row is $e_i$ and all the rest of its entries zero. 

Let $U$ in $\gl(n,S)$ satisfy $U+aE_i$ invertible for $i =1,2,\dots, n$. Then $(U+aE_i)U^{-1}$ is $\I + a D_i$, where $D_i$ is the matrix whose top row is the $i$th row of $U^{-1}$ and all the remaining entries are zero. 

It is routine to verify that invertibility of $\I + a D_i$ implies that $1 + a v_{i,1}$ is invertible in $S$, where $v_{i,1}$ is the corresponding entry of $U^{-1}$. For each $i$, we thus have $av_{i,1} \in V$, and thus $av_{1,i} = 0$. Since $v_{i,1}$ run over the first column of an invertible matrix, this implies $a = 0$, a contradiction. 
\qed

In particular, if $\Mn _n S$ satisfies \gui {n+1}, then $S$ satisfies \gui 1.

\comment
\Lem Lemma \Atty. Let $R$ be a commutative ring, and let nonzero  $a$ in $R$ satisfy $aR \cap V = \brcs{0}$. Then $\Mn_n R$ does not satisfy \gui {n+1}. 

\Pf Set $A(i) = a e_i$ (same notation as in the previous proof), and suppose that for some $U $ in $\gl (n,R)$, all of $U + a(i)$ are invertible. As before, 
$$
\det (U + A(i)) = \det U + a d_i
$$  
(where, once again, $d_i$ are the signed determinants of the minors). This expresses each $a d_i$ as a difference of units, that is, $ad_i \in V \cap aR$, so $ad_i = 0$. However, $\det U$ is an $R$-linear combination of the $d_i$, that is, $\sum d_i R = R$. This forces $a = 0$, a contradiction. 
\qed

\endcomment

Thus for neither $R = \Z$ nor $F[x]$, can $\M_3 R$ satisfy \gui 4. But I still don't know whether $\M_3 \Z$  satisfies 
\gui 3, although extensive calculation suggests it does. All   proper factor rings of $\Mn _2 \Z[1/2]$
satisfy \gui 4{}  (because of the result, mentioned above, that $\Mn _2 F_3$ does), but I don't know whether $\Mn _2 \Z[1/2]$  even satisfies \gui 3.

For the problem of whether (e.g., for a commutative domain $S$ satisfying \gui 2) $\Mn_n S$ satisfies \gui {3} if $n\geq 2$, we have the following tantalizing result. It is along the same lines as the proof of Proposition \Aeig. 

\Lem Lemma \Atwn. Let $S$ be a right or left Hermite ring satisfying \gui 2, and suppose that $B$ in $\M_{n} S$ is of the form $F \oplus 0$ (where the $0$ is of size at least $1$). Then for each $C$ in $\Mn_n S$, there exists $U$ in $\gl(n,S)$ \st both $U+B$ and $U+C$ are invertible. 

\Pf There exist $V, W$ in $\gl(n-1,S)$ \st  $VFW$ is upper triangular; then $V \oplus 1$ and $W \oplus 1$ will render $B$ upper triangular, with $n$th column consisting of zeros. By applying an elementary transformation (cyclically permuting the columns) we can arrange that the resulting form is $B_0$ strictly upper triangular (with zeros along the main diagonal); we relabel it $B$. At the same time, all these transformations are applied to $C$, and the relabelled matrix will still be called $C$. 

Now define $U$ in $\gl (n,S)$ as follows. For each of the diagonal entries $c_{i,i}$ of $C$, there exist invertibles $u_i$ in $C$ \st $u_i + c_{i,i}$ is invertible (this is \gui {2} for $S$). 
$$
U_{i,j} = \cases u_i & \text{if $i = j$} \\
-c_{ij} & \text{if $i < j$}\\
0 & \text{if $i > n$}\\
\endcases
$$
Then $U + B$ and $U$ are upper triangular with invertibles along the main diagonal, so are invertible, and $U + C$ is lower triangular with all diagonal entries invertible, so is also invertible.
\qed

This means that if we want to test, say for  a commutative domain that is also an elementary divisor ring  satisfying \gui 2, whether $\Mn_n S$ satisfies \gui 3, the problem reduces to both $B$ and $C$ being of full rank. Unfortunately, I haven't been able to do anything with this.

This can be a viewed as a variation on the proof of Proposition \Aeig. Another variation occurs in Appendix B.  

\comment
As an ad hoc definition, we call a ring $S$ {\it very weakly Hermite\/} if for all integers $n$ and $B \in \Mn_n S$, there exist $U,V$ in $\gl {n,S}$ \st $UBV$ is  triangular. This includes all right or left Hermite rings, along with a lot of others, including some that are very strange (for example,  rings \st for {\it all\/} nonzero $a \in \Mn_n S$, there exist invertibles $u,v$ \st $uav$ is the identity matrix; the right or left maximal ring of quotients of a domain that is not Ore satisfies this). 
\endcomment

\comment
The following is plausible because of the way \gui {2} works (Proposition xxx). It would lead to a proof that $\M_{2n} \Z$ satisfies \gui {3} for all $n >1$

\Lem Conjecture \Atwt. Suppose that $S$ is an elementary divisor ring satisfying \gui {2}. Then for all $n \geq 2$, $\Mn_n S$ satisfies \gui {3}. 
 \endcomment

Another approach for an interesting class of examples is via the following question. Let $S$ be a set of primes of $\Z$, and denote by $\Z[S^{-1}]$ the ring obtained by inverting all the members of $S$. Under what  conditions  will $\Z[S^{-1}]$ satisfy \gui 2? (or even \gui 3?). 

An obvious necessary condition is that $2 \in S$ (otherwise, $\Z[S^{-1}] $ maps onto $ F_2$). In the other direction, if $2 \in S$ and $S$ is cofinite (in the set of primes), then $\Z[S^{-1}]$ is semilocal, so that it (and its matrix rings) satisfy \gui {k} if and only if the finite factor rings---corresponding to the primes not in $S$---do (thus $\Z[S^{-1}]$ satisfies \gui 2, and will satisfy \gui {k} if all the primes not in $S$ are greater than $k$; since the Jacobson radical behaves well \wrt matrix rings, we obtain corresponding results for $\Mn_n \Z[S^{-1}]$. 

To decide whether $\Z[S^{-1}]$ satisfies \gui 2, we reduce to the following question: given an integer $m$ can we represent it as sum of difference of two elements, each which are of the form $s_1 ^{m(1)} $
if we take $S = \brcs{2,3}$, then $\Z[1/6]$ fails \gui 2. This is a consequence (as my colleague Damien Roy pointed out) of results on $S$-unit equations, [S, chapter 5, section 2]. It is plausible that if $S = \brcs{2,3,5}$, then $\Z[S^{-1}]$ satisfies \gui 2. 

At the other extreme, if $S$ is cofinite in the set of primes, then $R = \Z[S^{-1}]$ is semilocal and $R/J(R) \iso \prod_{p \notin S} F_p$. In this case, if $2,3 \in S$, then $R$ satisfies \gui {\min_{p \notin S} \brcs{p-1}}.

 \comment
A class of non-examples. 

\Lem Lemma. Let $S$ be a ring, and suppose that there exist $a_1, a_2, \dots, a_n$  elements of $S$ \st the right ideal they generate, $J:= \sum a_i S$ cannot be generated by a set with less than $n$, and moreover  $J \cap V = \brcs{0}$.  Then $\Mn_n S$ does not satisfy \gui 2, and $\M_{n+1} S$ fails to satisfy \gui 3.

\Pf Let $B$ the $n \times n$ matrix whose top row is $\(\matrix a_1 & a_2 & \dots & a_n \\ \endmatrix\)$ and whose remaining entries are zero. Suppose $U$ and $U + B$ are both invertible in $\Mn_n S$ for some matrix $U$. 

Then $(U + A)U^{-1}$ is invertible; all of its rows except the first  are simply the $i$ basic row. Denote the top row $(1 + a_1', a_2', \dots , a_n' )$. It is an easy exercise to check that this matrix being invertible implies $1 + a_1'$ is invertible; this implies in turn, that $a_1' \in V$.

We also have that $\sum a_i' S = \sum a_i S$, and in particular, $a_1' \in J \cap V$, forcing $a_1' = 0$. This entails that $J$ has an $n-1$-element generating set, a contradiction. 

To see that $\Mn_{n+1} S $ fails \gui {3}, let $B$ be the previous version of $B$ with a column of zeros adjoined to the left, and $C$ be the matrix whose only nonzero entry is a $1$ in the $(1,1)$ entry. The details are routine. 
\qed

For example, if $S = \Z[x]$, then $V = \brcs{0,\pm2}$. The ideal $J = < 3^{n-1}, 3^{n-2}x, \dots, x^n >$ cannot be generated by fewer elements, and does not contain $2$. Thus $\Mn_n \Z [x]$ fails \gui {2} for all $n$.  

Similarly, if $F$ is a field, then with $S = F[x,y]$ (the polynomial algebra in two variables), $\Mn_n S$ fails \gui {2} for all $n$, and this case is even easier since $V = F$, all of whose nonzero members are units. 

\endcomment
In particular, $\M_2 \Z$ fails to satisfy \gui {3} (which we already knew), and $\M_3 \Z$ fails 
\gui {4} (this is also a consequence of $\M_3 \Z_2$ failing to satisfy it, since for $S = F_2$, $V(S) = \brcs{0}$.

For Dedekind domains, $V(S)$ can be interesting. For example, let $S = \Z[\gamma]$ where $\gamma$ is the golden ratio. If $s \in V\setminus \brcs{0}$, then the norm of $s$ is plus or minus a Lucas number ($\brcs{1,3,4,7,11,18, \dots}$), and moreover, every Lucas number is achieved by an element of $V$. This type of result extends to  other quadratic extensions that have an infinite unit group (that is, not the purely imaginary ones).

In particular,  since the only square Lucas numbers are $1$ and $4$, the only integers in $V\setminus \brcs{0}$ are  $\pm1,\pm2$.

\comment
Thm 4.1  S satisfies: every fg right ideal is principal. Let A be a square matrix whose second through mth columns are all zero. Then there exists invertible U st AU is lower triangular. 

Definition of Hermite: if every 1 x 2 matrix admits diagonal reduction (that is, (a b) Q = (d 0) for some invertible Q), then $R$ is right Hermite; 2 x1  gives the definition of left Hermite. 

Really only need  PAQ theorem to get triangular (upper or lower)
\endcomment

\def \jnf{\text{JNF}}

\SecT Appendix B

Here $F_q$ (or simply $F$) denotes the finite field with $q$ elements. We show that for all $n> 1$, $\M_n F$ satisfies \gui {q}, and that for even $n$ and $q \neq 2$, $\M_n F$ satisfies \gui{q+1}. 

For a field $D$, its group of nonzero elements is denoted $D^*$.

Let $(X,\mu)$ be a normalized measure space (that is, $\mu (X) = 1$; $\mu $ is a probability measure). For a finite set $S$, let $\mu$ denote the normalized counting measure (if $T \subset S$, define $\mu$ via $\mu (T) = |T|/|S|$). We will use several such in the proofs below. First,  we  some elementary consequences of inclusion-exclusion.

\Lem Lemma \Asix. Let $(X,\mu) $ be a probability space, and let $X_1, X_2, \dots, X_k$ be measurable subsets of  $X$. Any of the  following conditions, (a), (b), or (c),  is sufficient for $\mu(\cap X_i) > 0$.{\par} 
{\noindent (a)} $\sum_{i=1}^k (1-  \mu(X_i))  < 1$; 
{\par}
\noindent (b) there exists a positive integer  $s \leq k/2$ and $\delta \geq 0$ \st 
$$\eqalign{
\sum_{j=1}^s \mu(X^c_{2j-1} \cap X^c_{2j}) & = \delta \quad\quad \text{and}\cr
\sum_{i=1}^k (1- \mu (X_i)) &< 1 + \delta;
\cr
}$$
\noindent (c) there exist $i < j < k$ \st 
$$\eqalign{
\(\mu (X_i^c \cap X_j^c) + \mu (X_i^c \cap X_k^c) + \mu (X_j^c \cap X_k^c)\) -\mu (X_i^c \cap X_j^c \cap X_k^c) &= \eta\quad\quad \text{and}\cr
\sum_{i=1}^k (1- \mu (X_i)) &< 1 + \eta. \cr
}$$

\Pf (a) 
$$\eqalign{
\mu(\cap X_i) & = 1 - \mu (\cup X_i^c) \cr
& \geq 1 - \sum \mu (X_i^c) \cr
& =   1 - \sum (1- \mu (X_i)) > 0.
}$$

\noindent (b) For each $j$, we have $\mu (X_{2j-1}^c \cup X_{2j}^c) = \mu (X_{2j}^c) + \mu (X^c_{2j}) - \mu(X^c_{2j-1} \cap X^c_{2j})  $. Summing over $j = 1$ to $s$, we obtain 
$$
\sum_{j=1}^s \mu(X_{2j-1}^c \cup X_{2j}^c) = \sum_{i=1}^{2s} \mu (X_i^c) - \delta. 
$$
Set $Y_j = X_{2j-1}\cap X_{2j}$. Thus 
$$\eqalign{
\sum_{j=1}^s \mu (Y_j^c)  + \sum_{i > 2s} \mu (X_i^c) & = \sum_{
i=1}^k \mu (X_i^c) - \delta < 1.
}$$
By (a), $(\cap Y_j )\cap (\cap_{i > 2s} X_i)$ has positive measure; but this intersection is just $\cap_{i=1}^k X_i$. 

\noindent (c) Inclusion-exclusion says
$$
\mu (X_i^c \cup X_j^c \cup X_k^c) = \mu (X_i^c) + \mu (X_j^c) + \mu(X_k^c) -\(\mu(X_i^c \cap X_j^c) + \mu (X_i^c \cap X_k^c) + \mu (X_j^c \cap X_k^c)\) +\mu (X_i^c \cap X_j^c \cap X_k^c).
 $$
 Set $Z = X_i \cap X_j \cap X_k$. Then $Z \cap \(\cap_{s \notin \brcs{i,j,k}} X_s\) = \cap_t X_t$ and 
 $$\eqalign{
 (1- \mu(Z^c) + \sum_{s \neq i,j,k} \(1 - \mu (X_s)\) & = \sum_t (1- \mu (X_t)) - \eta < 1.
 }$$
 So (a) applies, and thus $\cap X_t = Z \cap \(\cap_{s \notin \brcs{i,j,k}} X_s\)$ has nonzero measure.
 \qed

\Lem Lemma \Aten. Let $F = F_q$, $n \geq 2$, $R = M_n F$, and $v$ be an element of $R$ with rank $r$. Let $D$ be a  maximal subfield of $\M_n F_q$ (so $D \iso F_{q^n}$). Then 
$$
\left|\Set{d \in D \setminus \brcs{0}}{d + v  \text{ is not invertible}}\right| \leq \frac{q^n-1 - (q^{n-r}-1)}{q-1}.
$$

\Pf Denote by $T_v$ the set in the display. $D$ is a field unitally embedded in $R$, having $q^n $ elements, and the nonzero ones are of the form  $e^i$, $0 \leq i \leq q^{n} -2$ for a suitable choice of $e$ in $D$.  For $d \in D\setminus \brcs{0}$, define 
$$
V_d := \ker (d+v),
$$
an $F$-vector subspace of $F^{n\times 1} = F^n$, and 
$$
S_v:= \Set{d \in D \setminus \brcs{0}}{d + v \text{ is not invertible}},
$$
and this of course is $\Set{d \in D \setminus \brcs{0}}{V_d \text{ is not zero}}$. 

We first observe that if $d \neq d'$, then $V_d \cap V_{d'} = \brcs{0}$: pick $x \in F^n $ in the intersection, and note that $(d+ v)x = 0 = (d'+v)x$; thus $(d-d')x = 0$; since  $d-d'$ is a nonzero element of the field $D$ (and thus an invertible matrix), it follows that $x = 0$. [Warning: because elements of $D$ need not commute with $v$, the sum, $\sum V_d$, need not be a direct sum.] 

In addition, $V_d \cap \ker v = \brcs{0}$: for $x$ in the intersection, we have $vx = 0$ and $(d+v)x = 0$, so $dx =0$, again implying $x = 0$. 

If $V_d$ is not trivial, then it contains at least one one-dimensional subspace, and all of these are disjoint (modulo zero) from $\ker v$ and each other. The number of one-dimensional subspaces of $F^n$ that intersect only in zero with $\ker v$ is 
$$
\frac{q^n-1}{q-1} - \frac{q^{n-r}-1}{q-1}.
$$
Hence this is an upper bound for $|T_v|$. 
\qed

 These bounds are independent of the choice of $D$. Sometimes a specific choice of $D$ will result in improvements, but  in many cases, the estimates are actually achieved.

\Lem Proposition \Aele. If $n \geq 2$, then $\Mn _n F_q $ satisfies \gui q. 

\Pf Let $v_1, v_2, \dots, v_{q-1}$ be nonzero elements of $R = \Mn _n F_q$. We wish to show that  there exists an element $d \in \gl(n,F_q)$ \st all of $d + v_i$ are invertible. 

First assume that at least one of the $v_i$ is invertible. Without loss of generality, $i =1$. There exist invertible $U$ and $V$ in $R$ \st $Uv_1 V = \I$, and we relabel the  $Uv_i V$ as $v_i$.

Pick a maximal subfield  $D$ of $R$; it has $q^n-1$ nonzero elements, one of which is $\I$. 
For each $i$, let $S_i = \Set{d \in D\setminus \brcs{0}}{d + v_i \text{ is invertible}}$, and let $\nu $ denote the normalized counting measure on $D^* = D \setminus \brcs{0}$.  

 It is clear that $S_1= D^* \setminus \brcs{-\I} $, so $|S_1| = q^n-2$, and thus  $\nu (S_1) = 1 - 1/(q^n-1)$. 

By Lemma \Aten, the remaining $S_i$ satisfy $|S_i| \geq q^n-1 - (q^n-1)/(q-1)$, with strict inequality in case $v_i$ is not invertible. In particular, for these $S_i$, $\nu (S_i) \geq 1 -1/(q-1)$
Thus 
$$\eqalign{
\sum_{i=1}^k (1- \nu (S_i) & \leq (1- \nu(S_1)) + \sum_{i=2}^{q-1}  \frac 1{q-1}\cr
& = \frac 1{q^n-1} + \frac{q-2}{q-1} < 1. \cr 
}$$
Thus  Lemma \Asix(a) applies. 

Now suppose that none of the $v_i$ are invertible, that is have rank $n-1$ or less. By Lemma\, \Aele,  $|S_i| \geq (q^n-1)(1- 1/(q-1)) +1$, and so $1- \nu (S_i) \leq 1/(q-1) - 1/(q^n-1)$. Hence 
$$\eqalign{
\sum \( 1- \nu(S_i)\)& \leq \frac {q-1}{q-1} - \frac{q-1}{q^n-1} < 1. 
}$$
And again, Lemma \Asix(a) applies. 
 \qed

When $q= 2$, this only yields that when $n \geq 2$, then $\Mn _n F_2$ satisfy \gui 2, which was already known [V, Proposition 6]. The result of    appendix C is that $\Mn _n F_2$ satisfies \gui {3} for all $ n\geq 3$. 

\Lem Proposition \Bsev. If $q \geq 4$, and $n$ is even, then $\M_n F_q$ satisfies \gui {q+1}.

Owing  to Corollary \Anin(a), it suffices to prove this when $n =2$. There are a few differences between the treatments of odd and even $q$.

We require a number of results. For $V$ in $\M_n F_q$, define $\T_V := \Set{g \in \Gl(n,F_q)}{g+ V \in \Gl( n,F_q)}$; this is of course just $\Gl(n,F_q) \cap (\Gl(n,F_q) - V)$. And we define the measure $\mu$ on $\Gl(n,F_q)$ via $\mu (X) = |X|/|\Gl(n,F_q)|$ as usual. The complement of $X \subset \Gl(n,F_q)$, denoted $X^c$, is \wrt  $\gl(n,F_q)$. The identity matrix is denoted $\I$.

\Lem Lemma \Bten. Suppose $q \geq 3$ and $n = 2$. Then 
\item{(i)} $\mu (\T_{\I}) = 1- (q^2-2)/(q^2-1)(q-1)$.
\item{(ii)} If $V$ is one of $\diag (r,s)$ or $\(\smallmatrix t& 1 \\ 0 & t\\ \endsmallmatrix\)$ with $r,s,t \neq 0$ and $r = 1$ implies $s \neq 1$, then 
$\mu (\TT \I ,V.)\geq (q^2-q +1)/(q^2-q)(q^2-1) $.
\item{(iii)} If $V$ is rank one, then $\mu (\TT \I,V.) \geq 1/(q^2-1) $. 
\item{(iv)} If $V$ is rank one, then $\mu (\T_V)  = 1- q/(q^2-1)$.
\item{(v)} If $V$ and $W$ are both rank one, and $\brcs{V,W}$ is linearly independent, then $\mu (\TT V ,W.)\geq 1/(q^2-1) $.

\Rmk In (v), if $V = \lambda W$ for some nonzero scalar $\lambda$ (and $V$ is rank one), then $\TT V,W.$ is empty, which is useless. 

 \Pf (i) The set $\T^c_{\I}$ consists of invertible matrices having $-1$ as an eigenvalue. Since $n=2$, the other eigenvalue must also belong to $F_q^*$. Hence all such matrices must be conjugate to one of
  $$
 \brcs{\(\smallmatrix -1 & 0 \\ 0 & \lambda \\ \endsmallmatrix\)}_{\lambda \in F^*} \cup  \brcs{\(\smallmatrix -1 & 1 \\ 0 & -1 \\ \endsmallmatrix\)}
 $$
 For $\lambda \neq -1$, the number of elements conjugate to $\diag (-1,\lambda)$ is $
 |\gl(2,F)|/|F^* \times F^*|$ (since $\gl (2,F)$ acts transitively by conjugation, and the centralizer is $F^* \times F^*$). When $\lambda = -1$, the matrix is scalar, so the number of conjugates is $1$. Finally, the number of conjugates of the last matrix is $ |\gl(2,F)|/|F^* \times F|$ (the centralizer consists of polynomials in the matrix with nonzero determinant). Since these conjugacy classes are disjoint, we have 
 $$\eqalign{
 \mu (\T^c_{\I}) & = \frac{q-2}{(q-1)^2} + \frac 1{|\gl(2,F)|} + \frac{1}{(q-1)q} \cr
 & = \frac{q^2-2}{(q^2-1)( q-1)}.\cr
 }$$
 
 \noindent (ii) If $r = s \neq 0,1$, then the $\TT \I,V.$ consists of the invertible matrices whose eigenvalues are $-1$ and $-r$, hence are conjugate to $\diag (-1,-r)$ (since $r \neq 1$). For each $r$, the number in the conjugacy class is $|\gl(2,F)|/|F^* \times F^*|$; taking the union over $r$, we obtain $\mu(\TT \I,V.) = (q-2)/(q-1)^2$.
 
 Now suppose $r \neq s$, and both are unequal to $0$ or $1$. Let $g =  \(\smallmatrix a & b \\ c &d \\ \endsmallmatrix\)$ belong to $\TT \I,V.$. Taking determinants 
 $$\eqalign{
ad + (a+d+1) &= bc \cr
ad + (as + dr + rs) &= bc, \cr
 }$$
 and together with $a + d + 1 \neq 0$ (so that $\det g \neq 0$), these are sufficient for $g$ to belong. 
 We deduce $a(s-1) + d(r-1) = 1-rs$. If $a = -1$, then $d =-r \neq 0$ (so $a + d + 1 \neq 0$). The first equation is just $(a+1)(d+1) = bc$, so $bc = 0$, and there are thus $2q-1$ matrices of this form (arising from $b= 0$ or $c = 0$, with $a=-1$ and $d = -r$). 
 Similarly, if $d=-1$, then $a = -s$ and there are $2q-1$ matrices, and this set is disjoint from the matrices with $a = -1$. 
 
 If $a,d \neq -1$, then $bc \neq 0$, so there are $(q-1) $ choices for $(b,c)$. As $a = (s-1)^{-1} (1-rs - d(r-1))$, $a$ is determined by $d$ (and vice versa), and we have to exclude the cases that $a+ d =-1$. There are thus $q-2$ possible values for $d$, leading to $(q-2)(q-1)$.
 
 The total is then $(q-2) (q-1) + 4q-2 = q^2 + q$. 
  
 Suppose $s= 1$ and $r\neq 1$  (the same method applies if $r =1$ and $s \neq 1$). Then $d = -1$, so $a \neq 0$, and  $a$ varies over $F^*$, so the number of $(a,d)$ is $q-1$. From $bc = (a+1) (d+1) = 0$, there are $2q-1$  choices $(b,c)$, making $(q -1)(2q-1)$ in all.

Now suppose $V = \(\smallmatrix r & 1 \\ 0 & r \\  \endsmallmatrix\)$ with $ r\neq 0$. The equations are
 $$\eqalign{
ad + (a+d+1) &= bc \cr
ad + ((a + d)r + r^2 - c) &= bc. \cr
 }$$
 If $r = 1$, then $c = 0$, and thus $b$ is arbitrary, and $(a+1)(d+1) =0$. If $a=-1$, then $d\neq 0$ and vice versa, yielding $q(2q-1)$. 
 
 If $r \neq 1$, then $c$ is determined by $a$ and $d$, and $a + d+1 \neq 0$. If $c = 0$, then $r^2-1 + (r-1) (a+d) = 0$, whence $a+d = -(r+1) \neq -1$ (as $r\neq 0$). There are $q$ choices for $b$ and $q$ choices for $d$, giving us $q^2$ matrices. If $ c \neq 0$, then $b$ is uniquely determined by $(a,d)$, and we have $q^2-q$ choices (the first factor from the set of possible $(a,d)$ with $a+d + 1 \neq 0$); the total is $q(2q-1)$.

 \noindent (iii) We may conjugate the pair $\I,V$ with an arbitrary element of $\gl(2,F)$, so reduce to the case that $V $ is one of $\(\smallmatrix \lambda &0 \\ 0 & 0 \\\endsmallmatrix\)$ or $\(\smallmatrix 0 & \lambda \\ 0 & 0 \\ \endsmallmatrix\)$ ($\lambda \neq 0$). 
 
 In the trace $\lambda$ case, the equations are 
 
 $$\eqalign{
ad + (a+d+1) &= bc \cr
ad + \lambda d &= bc. \cr
 }$$
 Thus $(\lambda-1)d = a+1$. If $\lambda =1$, then $a = -1$, forcing $d \neq 0$ (to ensure that $a+ d +1 \neq 0$, so that $ad-bc \neq 0$). It also forces $bc = 0$, yielding $2q-1$ choices for $(b,c)$, whence the set has $(2q-1)(q-1) $ elements. 
 
 If $\lambda \neq 1$, then $d$ is uniquely determined by $a$. If $a=-1$, then $d = 0$ which is not permitted (as $a + d + 1 \neq 0$). So there are $q-1$ choices for $(a,d)$---if $d = -1$, then $a = -\lambda \neq 0$. In the latter case, $bc = 0$, so there are $2q-1$ elements of this form. 
 
 Otherwise, neither $a$ nor $d$ is $-1$, so $bc \neq 0$, and thus there are $q-1$ choices for $(b,c)$, whence $(q-2)(q-1)$ elements of this form (the $q-2$ arises from $a \neq -1, -\lambda$). 
  The total is thus $(q-2)(q-1) + 2q-1 = q^2 -q + 1$. 
 
 Finally, let $V$ be one of the nilpotent matrices. The second equation is then 
 $$
 ad - \lambda c = bc,
 $$
 yielding $c = -\lambda^{-1}(a+d + 1)$, and thus  $c$ is required to be nonzero. Hence $b$ is determined by $(a,d)$, the only constraint being $a+ d+1 \neq 0$. There are thus $q^2 - q$ elements in the set. 
 
\noindent (iv) The $P,Q$ trick allows us to assume $V = \(\smallmatrix 1 &0 \\ 0 & 0 \\\endsmallmatrix\)$. In this case, the inequalities for $M + V$ being invertible are $ad-bc \neq 0$ and $ad-bc + d \neq 0$. The negation of the latter is $d(a+1) = bc$. 

If $a = -1$, then $bc = 0$ forcing $ad \neq 0$, that is $d \neq 0$. Thus there are $q-1$ choices for $d$, and $2q-1$ choices for $(b,c)$, yielding $(q-1)(2q-1)$ elements.

If $a \neq -1$, then $d$ is determined by $(b,c)$, except that $d \neq 0$. Thus $bc$ is not zero, so there are $(q-1)^2$ choices for $(b,c)$, and $q-1$ choices for $a$, whence the number of matrices is $(q-1)^3$. 

The total is $(q-1)^3 + (q-1)(2q-1) = q^3 -q^2 $. 
 
\noindent (v) The $P,Q$ trick allows us to assume $V = \(\smallmatrix 1 &0 \\ 0 & 0 \\\endsmallmatrix\)$. This is invariant under various  transformations: multiplication of the second row or column by a nonzero scalar, adding a multiple of the second row or column to the first, and transpose (which is not elementary, but is permitted).  

Beginning with arbitrary rank one $W$, we apply these operations and  reduce to only two possibilities. 
Suppose $W_{12} \neq 0$. Multiplying the second column by a scalar, we can assume $W_{12} = 1$. Adding a multiple of the second column to the first, we can assume $W_{11} = 0$. Since the matrix is rank one, we must have $W_{21} = 0$. If now $W_{22} =0$, we have reduced to (the new) $W_1 = \(\smallmatrix 0 &1 \\ 0 & 0 \\\endsmallmatrix\)$. On the other hand, if $W_{22}\neq 0$, then we can multiply the second row by a scalar and add a multiple of it to the first, resulting in $W_2 = \(\smallmatrix 0 &0 \\ 0 & 1 \\\endsmallmatrix\)$. 

Now suppose that at  the outset, we have $W_{12} =0$, but $W_{22} \neq 0$. Then the rank one condition forces $W_{11} = 0$, and now we do the same as before, yielding two choices, $W_2, W_1^T$; but the latter is the transpose of the first choice, so needn't be considered. 

Finally, if the second column consists of zeros, we can transpose it, and reduce to the case $W_{11} \neq 0$ and $W_{12} = 0$; but this is a multiple of $V$, which violates linear independence. 

So we are reduced to considering only $W $ equalling one of $W_1$ or $W_2$.

If $W = W_1$, we have the two equations,
$$\eqalign{
ad + d &= bc \cr
ad - c & = bc, \cr
}$$
in addition to $d \neq 0$ (to force $ad-bc \neq 0$). This yields $c = - d \neq 0$. Thus $d(a+1) = -d b$, so $b = -(a+1)$. There are $q-1$ choices for $d$, and $q$ choices for $a$, resulting in $(q-1)q$ matrices. 

If $W = W_2$, the equations are 
$$\eqalign{
ad + d &= bc \cr
ad + a & = bc. \cr
}$$
Then $a = d$ and we must have $d \neq 0$, as previously. Then $bc = a(a+1)$. If $a = -1$, then there are $2q-1$ choices for $(b,c)$. If $a \neq -1$, then $bc \neq 0$, so there are $q-1$ choices for $(b,c)$, and thus $(q-1)^2$ matrices. The total is 
$(q-1)^2 + 2q-1 = q^2$. 
\qed

We need a few simple results about finite fields of even order. For the field $F =F_q$ with $q$ even, let $B = B(F)$ be $\Set {h^2  + h}{h \in F}$. It is easy to see that $B$ is an additive subgroup of $F$, and is the image of the $F_2$-linear map $F \to F$, $x \mapsto x^2 + x$; moreover, as the kernel consists of $\brcs{0,1}$, the image (that is, $B$) is an additive subgroup of index two in $F$. In addition, every element of $F$ is a square (since $F^*$ is a group of even order), and of course, the square root is unique. Moreover, $B$ is stable \wrt taking squares or square roots.

\Lem Lemma \Beig. Let $F = F_q$ with even $q$, and suppose $a,b \in F_q$. {\par}\item{(a)} The polynomial $g$, defined by $g(x) = x^2 + ax + b$, is irreducible over $F$, if and only if $a, b \neq 0$ and $b/a^2 \notin B(F)$. {\par}
\item{(b)} Let $\rho \in F\setminus\brcs{0,1}$. There are $q/2$ elements $s \in F\setminus \brcs{0,1}$ \st $s\rho/(s^2+1)$ does not belong to $B(F)$.

\Pf (a) Since $g$ has degree two, irreducibility is equivalent to no roots in $F$. Necessity of $ab  \neq 0$ is clear. Rewrite $g(x) = a^2 f(x/a)$, where $f (x) = x^2 + x + b/a^2$. If this has a root in $F$, then $b/a^2 \in B$, and conversely. 

\noindent (b) Replace $s$ by $s+1$ (so the new $s$ is also not zero or one); then $s\rho/(s^2+1)$ is replaced by $(s +1)\rho/s^2$. There exists $\alpha$ \st $\alpha^2 = \rho$; obviously, $\alpha^2  + \alpha \neq 0$. Then
$$\eqalign{
\frac{(s +1)\rho}{s^2} & = \(\frac {\alpha}{s}\)^2 + \frac {\alpha}{s} + \frac{\alpha^2 + \alpha}{s} \cr
& \equiv   \frac{\alpha^2 + \alpha}{s}      \mod B.
}$$
As $s$ varies over $F_q\setminus\brcs{0,1}$, $1/s$ varies over the same set. If $s = 1$, the outcome is in $B$, and since $0$ is also in $B$, that leaves $q/2$ values for $s$ \st the fraction is not in $B$, since $|B| = q/2$.
\qed

\noindent {\it Proof of Proposition \Bsev.} Let $S:= \brcs{v(1), v(2), \dots , v(q)}$ be a set consisting of $q$ distinct elements of $\M_2 F_q$ for which we want to find $M := \(\smallmatrix a & b \\  c & d\\ \endsmallmatrix\) \in \gl(2,F)$ \st all $M+ v(i)$ are invertible.  

\noindent (a) {\it All the matrices $v(i)$ are rank one.} By Lemma \Aele(iv), $1- \mu(\T(v(i)) \leq q/(q^2-1)$, so that $\sum (1- \mu(\T(v(i))) \leq 1 + 1/(q^2-1)$. Since $q > q-1$, there exists $i$ \st $v(i)$ is not a scalar multiple of any of the others.

Hence there exists $ j \neq i$ \st $\brcs{v(i),v(j)}$ is linearly independent. In addition, there exists either a pair $i(',j')$  with distinct  entries which are also distinct from $\brcs{i,j}$ \st ${v(i'),v(j')}$
is linearly independent, or there exists $k$ unequal to $i$ or $j$, \st $\brcs{v(i),v(k)}$ is linearly independent (to see this, consider at the worst case: $v(t) = \lambda_t v$ for some rank one $v$ with $\lambda_t$ running over all of $F_q^*$, and one more element $v(q)$). In the former case, using the notation of Lemma \Asix, $\delta \geq 2/(q^2-1)$ by Lemma \Aele, so  \Asix(b) applies. In the latter case (with the $k$), $\eta \geq 2/(q^2-1)$ (since $\eta$ is at least as large as the sum of the measures of two of the intersections of the complements), and \Asix(c) applies. 

\noindent (b) {\it There exist $i \neq j$ and invertible $P,Q$ \st $Pv(i)Q = \I$ and either $Pv(j) Q$ is a scalar matrix, or its characteristic polynomial is irreducible over $F_q$. }

We may thus assume (after re-ordering, and replacing each $v(i)$ by $Pv(i)Q$, and relabelling) that $v(1) = \I$ and $v(2) = \delta$ the latter either a scalar matrix, or it has  irreducible characteristic polynomial. In either case, there is a maximal subfield of $\M_2 F_q$, $D$, containing $\delta$. Obviously, the cardinality of $D$ is $q^2$.  The normalized counting measure, $\nu$, is defined \wrt $D^*$, as are complements. 

Define $X = D^*$ with normalized measure $\nu$ given by  $\nu(Y) = |Y|/(q^2-1)$. For $V \in \M_2 F_q$, set $\II (V) = \Set{d \in D^*}{V + d \in \gl(2, F)}$, and set $Y_i = \II (v(i))$. Then $\nu (Y_1) = \mu (Y_2) = 1 - 1/(q^2-1)$, and by Lemma \Aten\ with $n=2$, $\nu(Y_i) \geq 1 - ((q^2-1)/(q-1))/(q^2-1) = 1 -1/(q-1)$. Thus 
$$
\sum_{i=1}^q (1-\nu(Y_i)) \leq \frac 2{q^2-1} + \frac {q-2}{q-1} = \frac q{q+1} < 1. 
$$
By Lemma \Asix(a), $\cap Y_i $ is not empty, and we are done with this case. 

From this point, in dealing with subsets $S =\brcs{v(1), \dots, v(q)}\subset \M_2 F_q$ with at least one of the elements being invertible, we may assume that for all $i \neq j$,

\item{(*)}If $v(i)$, $v(j)$ are invertible, then $v(i) ^{-1} v(j) $ is neither a scalar matrix, nor  has reducible characteristic polynomial.

\noindent Property (*) is invariant under left and right multiplications by invertibles, as is easy to check. So we may assume (if necessary) that $S = \brcs{\I, v(2), \dots, v(q)}$, still satisfying (*). So we assume that $S$ is a counter-example to the assertion of the proposition, and derive a contradiction by successively narrowing the possibilities. 

\noindent {(c)} {\it $S $ satisfies $(*)$, at least one of its elements is invertible, and at least two of its elements are rank one.} Let $D$ be any maximal subfield of $\M_2 F_q$, and let $\nu$ be the normalized measure thereon, and let $\II_i = \Set {d \in D^*}{d+v(i)} \in \gl(2,F)$. Since $v(1) = \I$, by Lemma \Aele,
$$\eqalign{
\sum (1- \nu (\II_i))& \leq  \frac 1{q^2-1} + \frac{q-3}{q-1} + \frac{2q}{q^2-1} \cr
&\leq \frac{1 + q^2-2q-3 + 2q}{q^2-1} = \frac{q^2-2}{q^2-1} < 1;\cr
}$$
thus Lemma \Asix(a) applies.

\noindent (c) {\it At most one element of $S$ is not invertible and $S$ satisfies $(*)$.} We observe that if $v(i)$ and $v(j)$ have a common nonzero eigenvalue, say $\lambda$, then $-\lambda I \in \II_i^c \cap \II_j^c$. Assume that $v(q)$ has rank one to begin with. For any choice of maximal subfield $D$, $\nu (\II_q) \leq q/(q^2-1)$. 

Now consider the eigenvalues of the rest of the invertible matrices: since $v(1) = \I$, it has only $1$ as an eigenvalue. If any of $v(2)$ through $v(q-1)$ have two distinct eigenvalues in $F_q^*$ (they always have at least one in $F_q^*$ by (*)), we would obtain an overlap, that is  ($v(i), v(j)$ have an eigenvalue in common; then  at least one $\nu(\II_i^c \cap \II_j^c) \geq 1/(q^2-1)$, and then the formula \Asix(a)  yields the result (the corresponding sum is less than one). 

Hence each of $v(2), \dots, v(q)$ has only one eigenvalue, and are nonscalar (by (*)); moreover, the eigenvalues must be distinct (otherwise an overlap occurs). If $v(q)$ has a nonzero eigenvalue, we are once again going to get an overlap---hence $v(q)$ is nilpotent. 

We say a matrix $V$ in $\gl (n,F)$ is {\it monopotent} [HL, p\,425] if it has exactly one  eigenvalue in the algebraic closure of $F$ (which of course implies the eigenvalue belongs to $F$), and the eigenvalue is nonzero. Alternatively, $V = \lambda \I  + N$ where $N$ is nilpotent and $\lambda$ is in $F^*$. Of course, when $n =2$, $N^2=0$, and $N $ can be written in the form  $\(\smallmatrix a\\  b\\ \endsmallmatrix\)\(\smallmatrix c & d \\ \endsmallmatrix\)$ subject to $ac + bd = 0$. We say $V$ is {\it properly monopotent\/} if it is monopotent but not a scalar multiple of the identity. 

Let $V$ be a properly monopotent matrix in $\M_2 F_q$. We associate to $V$ its {\it type\/}, a symbol in $F_q \cup \brcs{\infty}$, and its {\it form,} a pair whose first coordinate is its type, and whose second coordinate is in $F_q^*$, as follows ($\lambda $ is an nonzero element of $F_q$). \hfill\break
$
\lambda \(\I + \kappa \(\matrix \alpha & -\alpha^2 \\ 1 & -\alpha \\\endmatrix\) \)
$ has type $\alpha$ and form $(\alpha, \kappa)$; \hfill\break
$\lambda \(\I + \kappa \(\matrix 0 & 1 \\ 0 & 0 \\\endmatrix\) \)$ has type $\infty$ and form $(\infty,\kappa)$. 

Every properly monopotent matrix has exactly one form. Transpose interchanges the forms $(0,\kappa)$ with $(\infty,\kappa)$;  if $\alpha \in F_q$, transpose sends $(\alpha, \kappa)$ to $(-\alpha^{-1}, -\kappa \alpha^2)$. It is also easy to check that if $V$ and $V'$ are properly monopotent with forms $(\alpha,\kappa)$ and $(\alpha',\kappa')$ respectively, then $VW = WV$ if and only if $\alpha = \alpha'$; moreover, $V^{-1}W$ is a scalar matrix if and only if $(\alpha,\kappa) = (\alpha',\kappa')$. 

\Lem Lemma \Bnin. Suppose   that $V,V' \in \gl(2,F_q)$ are properly monopotent with eigenvalues $\lambda, \lambda'$ and forms $(\alpha, \kappa), (\alpha',\kappa')$ respectively. Suppose that $V^{-1}V'$ is  properly monopotent and $\alpha \neq \alpha'$.
{\par}\item{(a)}If $q$ is odd,  
\itemitem{(i)} if $\alpha,\alpha' \in F_q$, then  $\kappa \kappa' (\alpha-\alpha')^2 = 4$  and the eigenvalue of $V^{-1 }V' $ is $- \lambda'/\lambda$;
\itemitem{(ii)} if $\alpha = \infty$ and $\alpha' \in F_q$, then $\kappa\kappa' = 4$,
 the form of $V^{-1}V'$ is $(\alpha'-\kappa/2,- \kappa')$, and the eigenvalue $-\lambda'/\lambda$. 
 {\par}\item{(b)} 
 $q$ is odd. 
 
 \Rmk Part (b) means that if $q$ is even, then we obtain a contradiction to $\alpha \neq \alpha'$.

\Rmk The minus signs in the eigenvalues  are not misprints. 

\Rmk The prime symbol ($'$) does not mean transpose (and I wish people would stop using it for that purpose---superscript $T$ works perfectly well). 

\Pf (a) Suppose that $\alpha \neq \alpha'$.

\noindent (i) 
$$\eqalign{
V^{-1} V' & = \lambda^{-1}\(\I - \kappa \(\matrix \alpha \\ 1\\ \endmatrix\)\(\matrix 1&  -\alpha \\ \endmatrix\)\)
\lambda' \(\I + \kappa' \(\matrix \alpha' \\ 1\\ \endmatrix\)\(\matrix 1&  -\alpha' \\ \endmatrix\)\)\cr
& = \frac {\lambda'}{\lambda} \(
\I - \(\matrix  \kappa \alpha- \kappa' \alpha'\ & -\kappa \alpha^2 + \kappa'{\alpha'}^2 \\ \kappa- \kappa' & -\kappa \alpha + \kappa' \alpha' \\ \endmatrix\) - \kappa\kappa' (\alpha' - \alpha)\(\matrix \alpha & -\alpha\alpha'\\ 1 & -\alpha'\endmatrix \)
\).
}$$
The trace is $(\lambda'/\lambda)(2 + \kappa \kappa'(\alpha-\alpha')^2$, and this is nonzero. Assuming $V^{-1}V'$ is monopotent, let $\rho$ denote is eigenvalue. Then $\rho^2 = \det V^{-1}V' = (\lambda'/\lambda)^2$, so $\rho = \pm \lambda'/\lambda$. Then the trace is $2\rho = \pm 2\lambda'/\lambda$.

If $\rho = \lambda'/\lambda$, then $\kappa\kappa'(\alpha-\alpha')^2 = 0$ (taking traces), a contradiction. Hence $\rho = -\lambda/\lambda'$, and equating traces, we deduce $\kappa\kappa' (\alpha-\alpha')^2 = 4$.   The type of $V^{-1}{V'}$ is just the ratio of the $(1,1)$ entry of $\lambda'/\lambda + V^{-1}V'$ to to the $(1,2)$ entry, a big mess (but could be either zero or $\infty$), and and the second coordinate of the form(the new kappa)  is the negative of the  $(1,2)$ entry  $(\kappa' - \kappa)\lambda'/\lambda - 4/(\alpha-\alpha')$. 

\noindent (ii) This time
$$\eqalign{
V^{-1}V' & = \lambda^{-1} \(\matrix 1&- \kappa \\ 0 & 1\\ \endmatrix\) 
\lambda' \(\matrix 1 + \kappa'\alpha' 
&- \kappa'{\alpha'}^2 \\ \kappa' & 1- \kappa' \alpha'\\ \endmatrix\) \cr 
& = \frac{\lambda'}{\lambda }\(\matrix 1 + \kappa'\alpha' -\kappa\kappa' 
&- \kappa'{\alpha'}^2 -\kappa  + \kappa\kappa' \alpha' \\ \kappa' & 1- \kappa' \alpha'\\ \endmatrix\). \cr 
}$$

As before, let $\rho $ be the eigenvalue of $V^{-1}V'$ (assuming monopotence); then $\rho^2 = (\lambda'/\lambda)^2$ from the determinants. The trace is $(\lambda'/\lambda)(2 - \kappa \kappa')$. If $\rho =\lambda'/\lambda$, then, as before $\kappa\kappa' = 0$, a contradiction. Thus $\rho =-\lambda'/\lambda$, and equating traces yields $\kappa\kappa' = -4$, so $\kappa\kappa' = 4$. 

The type is $2/\kappa' + \alpha'- \kappa = \alpha' - \kappa/2$, and the new $\kappa$ is $-\kappa'$. 

\noindent (b) If $q$ is  even, then $2 = 0$. If both $\alpha,\alpha' \in F_q$, then $\kappa\kappa'(\alpha-\alpha')^2 = 4 = 0$, so 
$\kappa\kappa' = 0$ forcing one of $\kappa, \kappa'$ to be zero, contradicting their definition. 

If $\alpha = \infty \neq \alpha'$, we still obtain $\kappa\kappa' = 0$, again a contradiction.\qed

It follows that for even $q$, all the types of the monopotent matrices in $S$ must be equal, and so they commute. For odd $q$, essentially the same thing applies (if $q \geq 5$; $q = 3$ requires a little extra), but it takes longer to get there. 

Return now to $S = \brcs{\I, v(2), \dots, v(q)}$ satisfying (*), and also with $\rk v(i) = 2$ for $i < q$ and $v(q)$ has rank one. We have seen that $v(q)$ must be nilpotent, and all $v(i)$ must be monopotent, and moreover, $v(2)$ through $v(q-1)$ are not scalar matrices. In addition, we have that the eigenvalues of $v(1)$ through $v(q-1)$ are distinct (hence exhaust $F_q^*)$. We can conjugate $S $ so that the new $v(1) = \I$ and $v(2)$ is $\lambda_2\(\smallmatrix   1 & 1 \\ 0 & 1 \\\endsmallmatrix  \)$ for some $\kappa$, and all the just-stated properties still hold. 

 We have that $v(2)$ has form $(\infty, 1)$. We claim that if any  other $v(i)$ ($q > i > 2$) has type $\infty$, then they all do; this is automatic if $q$ is even, so we assume for now that $q$ is odd. Suppose that  $v(j)$ has form $(\infty, \kappa_j)$, but $v(k)$ has type $(\alpha, \kappa_k)$ for some $\alpha \in F_q$. From $v(j)^{-1} v(k)$, being (by assumption) properly monopotent, we have $\kappa_j \kappa_k = 4$. Similarly, $1\cdot \kappa_k = 4$ and thus $\kappa_j = 1$ (as $\kappa_k \neq 0$ from the definitions). But this means the forms of $v(2)$ and $v(j)$ are equal, and so $v(2)^{-1}v(j)$ is a scalar matrix, a contradiction. 

To apply, this we need $q \geq 3$. 

Recalling that  $v(q)$ is rank one, we already know that it must be nilpotent. We claim that in our current state ($v(1) = \I$, $v(2)$ has type $\infty$),  that $v(q)$ is upper triangular. Otherwise it is in the form $k\(\smallmatrix \alpha & -\alpha^2 \\ 1 & -\alpha \\ \endsmallmatrix \)$ for some $k \in F_q^*$ and $\alpha \in F_q$. Then 
$$\eqalign{
v(2)^{-1} v(q) &= \lambda^{-1}\(
 \matrix 1 & -\kappa \\ 0 & 1\\  \endmatrix\)k\(\matrix \alpha & -\alpha^2 \\ 1 & -\alpha \\ \endmatrix\)  \cr
& = \lambda^{-1} k \(\matrix \alpha- \kappa & * \\ 1 & \alpha \endmatrix \).\cr
}$$

In the new $S$ (the set $\brcs{v(i)^{-1}v(j)}_j$), $v(2)^{-1}v(q)$ is rank one, and since all the preceding applies to the new sequence, we must also have $v(2)^{-1}v(q)$ is nilpotent. But its trace is $-\lambda^{-1}k\kappa$, which is not zero, a contradiction. 

So by a further conjugation if necessary, we can assume $v (2) $ has type $(\infty , 1)$ and $v(q) = \(\smallmatrix 0 & k \\ 0 & 0 \\ \endsmallmatrix \)$.

At this stage, suppose there exists $j$ \st $v(l)$ has type $\alpha \neq \infty$. The new sequence $\brcs{v(l)^{-1}v(i)}$ contains the identity, and exactly one rank one matrix; by our previous reductions, this must be nilpotent. But 

$$
\lambda_l^{-1} \(\I - 4 \(\matrix \alpha & -\alpha^2 \\ 1 & -\alpha \\ \endmatrix\)\)  \(\matrix 0 & k \\ 0& 0 \\ \endmatrix \) =  \(\matrix 0 & * \\ 0 & -4k/\lambda \\ \endmatrix \);
$$
as the trace is not zero, the matrix is not nilpotent, a contradiction. 

So all of $v(2), \dots , v(q-1)$ are monopotent with the same type; and this is true for all $q$, and we may assume the type  $\infty$ (and thus the matrices are upper triangular), and $v(q)$ is upper triangular nilpotent. Moreover, all the eigenvalues must be distinct. Otherwise for some $i \neq j$, $v(\II_i^c \cap \II_j^c) \geq 1/q^2-1$ for any choice of maximal field $D$ and its normalized measure $\nu$, and since $\nu (\II_q^c) \leq (q/q^2-1) $, we deduce that 
$$\sum (1- \nu(\II_i)) \leq  \frac 1{q^2-1} + \frac {q-2}{q-1} + \frac q{q^2-1} = 1,$$
but now   Lemma \Asix(a) applies.

Thus, regardless of the parity of $q$, for $i = 2, \dots, q-1$, we must have $v(i)$ has form $(\infty, \kappa_i)$ where the $q-2$ $\kappa$s must be distinct (else, there is a nontrivial $v(i)^{-1}v(j)$ that is a scalar matrix). 

We claim there is $1 < j <q $ \st $\II_j^c \cap\II_q^c \neq \emptyset$ \wrt some choice of maximal subfield $D$ (the sets $\II_j$ depend on the choice of $D$). This is sufficient to reach a contradiction by the same remark as just previous. This requires two different arguments, one for odd $q$ and another for even $q$. 

We do not determine the maximal subfield $D$ until the last second. Suppose now that $q$ is odd. Consider the polynomial $f (x) = x^2 - k x = (x-k/2)^2 -k^2/4$. As is well-known, there are at least $(q-1)/2$ values of $T$ in $F_q$ \st $f(T)$ is not a square. As $q \geq 3$, $q-2 \geq (q-1)/2$. Thus there exists $v(j)$ with $1 < j < q$ with form $(\infty, \kappa)$ and eigenvalue $\lambda$, where $f(\kappa)$ is not a square. This also entails $k \neq \kappa$. Set $t = (\kappa^2 - k\kappa)/\lambda^2$. This is not a square in $F_q$, and thus $e:= \( \smallmatrix 0 & t \\ 1 & 0 \\\endsmallmatrix \)$ has irreducible characteristic polynomial.

Form $D = F_q [e] $ inside $M_2 F_q$; this is of course a maximal field. Every element of $D^*$ is uniquely representable in the form $r\I + se$ where $r,s$ run over $F_q$ and not both are zero. Set $r =-\lambda$ and  $s = \lambda^2/(k-\kappa)$, and define $d = r\I + se$. It is easily verifiable that \wrt the current choice of $D$, $\II_j^c \cap \II_q^c$ contains $d$, and so $\nu (\II_j^c \cap \II_q^c) \geq 1/(q^2-1)$. 

Now assume that $q$ is even; we wish to show that for a suitable maximal subfield $D$, that $\II_i ^c\cap \II_q^c$ is not empty.  We know that all the matrices $v(2)$ to $v(q-1)$
 have the same type, and thus commute, and by conjugation, we can assume $v(i) = \lambda_i \(\smallmatrix 1 & \kappa_i \\ 0 & 1 \endsmallmatrix \)$; we can also assume that the $\lambda_i$ are distinct  (else a pair of them has a common eigenvalue, so that $\II_i^c \cap \II_j^c \neq \emptyset$ for any choice of $D$ but some choice of $2 \leq i < j \leq q-1$. We may also assume that $\kappa_i \lambda_i$ are distinct (otherwise there is a pair with $v(i)^{-1}v(j) $ being scalar). 
 
 Consider the pair 
 $$
\( \matrix \lambda_i & \lambda_i\kappa_i \\ 0 & \lambda_i\\ \endmatrix\) \quad\text{and}\quad \( \matrix 0& k \\ 0 & 0 \\ \endmatrix\).
$$
If $k = \kappa_i \lambda_i$, let  $d = s\(\smallmatrix a & b\\ c& 0 \endsmallmatrix \) $, where $a, b,c ,s$ are to be determined so that $d$ has irreducible characteristic polynomial, and $d$ added to each of the displayed matrices results in noninvertible ones. The conditions $sa = \lambda$ and $sb = k$ (so that none of $s,a,b$ are nonzero)  take care of the first one; to the obtain the other one, we see that $s c (s b+ k) = 0$ is sufficient, and of course we have this. The characteristic polynomial of $e:= d/s$ is $x^2 + ax + bc$, so we only require that $bc/a^2 \notin B$.
$$
\frac{bc}{a^2} = sc\frac {k}{\lambda^2};
$$
since we can let $sc$ run over $F_q^*$, there exist plenty of choices for $s$ and $c$ so that $d/s$, and thus $d$, has irreducible characteristic polynomial.

 If $k \neq \kappa_i \lambda_i$, let  $e = s\(\smallmatrix 1 & b\\ c& h \endsmallmatrix \) $, where $ b,c,h,s $ are to be determined so that $d$ has irreducible characteristic polynomial, and $d$ added to each of the displayed matrices results in noninvertible ones. The conditions $s = \lambda_i$ and $b = \kappa_i$ take care of the first one. For the second, necessary and sufficient is the equation $s^2 h = sc(\lambda_i \kappa_i + k)$. As $s$ is nonzero, this reduces to $c = \lambda h/(\lambda_i \kappa_i + k)$. Then the requirement for irreducibility of $x^2 + (h+1) x + h + bc$ (it is enough to show $d := e/s$ has irreducible characteristic polynomial) is  $(h+bc)/(h+1)^2 \notin B$. Since
$$
\frac{h + bc}{(h+1)^2} = \frac h{h^2 +1}\frac {k}{k+ \kappa_i \lambda_i},
$$
by Lemma \Bten(b), with $\rho = k/(k+ \kappa_i \lambda_i)$, there are $q/2$ choices for $h$. 

This concludes the argument with exactly one $v(i)$ of rank one.

Now we assume that all the $v(i)$ are invertible, in addition (*), etc.
Condition (*)
applies, and so if either there is a triplet with a common eigenvalue, or there are two pairs (which need not be disjoint), we let $D$ be any maximal subfield of $\M_2 F_q$, and then Lemma \Asix(b) or (c) applies, and we may subtract $2/(q^2-1)$ from $\sum (1- \nu (\II_j)) \leq 1/(q^2-1) + (q-1)/(q-1) = 1 + 1/(q^2-1)$, and the result is less than one.

There are $q$ nonzero eigenvalues resulting from the $v(i)$, so there must be one pair $v(i), v(j)$ with the same eigenvalue---but by the previous paragraph there cannot be another pair. This precludes all the $v(i)$ from having two different eigenvalues, so that all the $v(i)$ are monopotent. 

By multiplying by a suitable scalar, then by the inverse of one of the elements, we may assume $v(1) = \I $ and $v(2) = \(\smallmatrix 1 & 1 \\ 0 & 1 \\ \endsmallmatrix\)$. In particular, for any choice of $D$, $\nu (\II_1^c \cap \II_2^c) = 1/(q^2-1)$, which a start.

Now consider the types of the remaining terms. If any of $v(2)$ throught $v(q)$ has type $\infty$, then they all do, that is, $S$ consists of elements in the commutative algebra of striped upper triangular matrices.

Otherwise, the remaining terms must have types other than $\infty$. If two have the same type, then each is polynomial in the other, and this property is preserved by conjugation, and upon conjugation (to render both with type $\infty$), and because the identity matrix is unchanged, and thus all of the conjugated terms must have the same type. In particular, after conjugation they all commute with each other, and this is preserved by unconjugation, so the previous system, with $v(2)  = \(\smallmatrix 1 & 1 \\ 0 & 1 \\ \endsmallmatrix\)$, all have type $\infty$. Thus in all cases, all the non-identity elements are of the form $\lambda_i  \(\smallmatrix 1 & k_i \\ 0 & 1 \\ \endsmallmatrix\)$ ($k_i$ is easier to type than $\kappa_i$), with $k_2 = 1$ and $\lambda_2 = 1$ . 

For $i \geq 2$, the form of $v(i)$ is $(\infty,k_i)$; if $k_i = k_j$, then $v(i)^{-1}v(j)$ would be a scalar matrix, a contradiction. Hence the $k_i$ are distinct, thus run over all of $F_q^*$. Moreover, the eigenvalues of $\brcs{v(2), \dots, v(q)}$ must be distinct by the earlier comment, so exhaust $F_q^*$.

Suppose that for some $i \neq j$, $\lambda_i k_i = \lambda_j k_j$. Then we show that  $\II_i^c \cap \II_j^c$  is nonempty \wrt a suitable maximal subfield $D$. Once we prove this, we can apply Lemma \Asix(b), and reduce to the case that each of the three sets $\brcs{\lambda_i}_{i\geq 2}, \brcs{k_i}_{i\geq 2}, \brcs{\lambda_i k_i}_{i\geq 2}$ is all of $F_q^*$).

Suppose $\lambda_i k_i = \lambda_j k_j$. Relabel (for convenience) $z = \lambda_i k_i$, $\lambda = \lambda_i$, and $\lambda' = \lambda_j$. The two matrices are then 
$$
\( \matrix \lambda & z \\ 0 & \lambda \endmatrix\) \quad\text{and}\quad \( \matrix \lambda'& z \\ 0 & \lambda' \endmatrix\).\tag1
$$
where $\lambda \neq \lambda'$ and $z\neq 0$. 
We show there exists a maximal subfield $D$ of $\M_2 F_q$ and a nonzero element $d$ of $D$ \st when $d$ is added to each of these matrices, neither  is invertible. In particular, \wrt this choice of $D$, $\II_i^c \cap \II_j^c $ is nonempty, and this is enough to apply Lemma \Asix(b) (observe that $v(1)$ and $v(2)$ have an eigenvalue in common, so that $\II_1^c \cap \II_2^c$ is not empty for any choice of $D$).  

First assume that $q $ is odd. We may find  $s $ in $F_q^*$ \st $(\lambda-\lambda')^2 - 4zs$ is not a square in $F_q$. Set $t = (4^{-1}s^{-2})((\lambda-\lambda')^2 - 4z)$, so that $t$ is also not a square. Set $r = -(\lambda + \lambda')/2$. Define  $e = \( \smallmatrix  0 & t \\ 1 & 0\\ \endsmallmatrix\)$. Since $t$ is not a square, the ring generated by $e$ is a quadratic extension field of $F_q$, and of course it is maximal; call it $D$. Set $d = r\I + se$; it is routine to verify that adding $d$ to each matrix results in matrices of determinant zero. 

Now we deal with the case $q$ is even. For $c \in F_q^*$, let $u_c:= \(\smallmatrix \lambda' & z \\  c & \lambda \\ \endsmallmatrix \)$. When  $u_c$ is added to each of the displayed matrices, the resulting matrices  have determinant zero. The characteristic polynomial (depending on $c$) of $u_c$ is $f_c:= x^2 + (\lambda + \lambda') x + cz$. Since $\lambda \neq \lambda'$ and $z$ is nonzero, we can choose $c$  so that $f_c$ is irreducible. For this choice of $c$, set $d = u_c$, and set $D = F_q[d]$. Then $d \in \II_i^c \cap \II_j^c$, so $\nu (\II_i^c \cap \II_j^c) \geq 1/(q^2-1)$, and we are done with this portion. 

\comment
As before, we consider the two matrices 
$$
\( \matrix \lambda_i & \kappa_i \\ 0 & \lambda_i\\ \endmatrix\) \quad\text{and}\quad \( \matrix 0& k \\ 0 & 0 \\ \endmatrix\).
$$
For $c \in F_q^*$, let $u_c:= \(\smallmatrix \lambda' & z \\  c & \lambda \\ \endsmallmatrix \)$. When  $u_c$ is added to each of the displayed matrices, the resulting matrices  are both not invertible. The characteristic polynomial (depending on $c$) of $u_c$ is $f_c:= x^2 + (\lambda + \lambda') x + cz$. This will be irreducible if $cz/(\lambda + \lambda')^2 \notin B$  ($\lambda \neq \lambda'$ and $z$ is nonzero). Thus we can choose $c$  so that $f_c$ is irreducible. For this choice of $c$, set $d = u_c$, and set $D = F_q[d]$. Then $d \in\II_i^c \cap \II_j^c$, so $\nu \II_i^c \cap \II_j^c) \geq 1/(q^2-1)$, and we are done with this portion. 
\endcomment

So we can assume that all the $\lambda_i k_i$ (as $i $ runs through $2$ to $q$) are distinct. Set $z_i = \lambda_i k_i$. Three sets $\brcs{\lambda_i}, \brcs{k_i}, \brcs{\lambda_i k_i} $ are all just $F_q^*$. If $q$ is odd, then $F_q^*$ is cyclic of even order, and this implies that if $\pi$ is a permutation of $F_q^*$, then $g\mapsto \pi(g)g$ cannot be a  permutation  \plainfootnote{$^1$}{This is an old result; see the comment {{https://mathoverflow.net/questions/508595/permutations-of-finite-groups}} due to  Darij Grinberg.}. Thus if $q$ is odd, we arrive at a contradiction. 

In the case that $q$ is even (so that $F_q^*$ is cyclic of odd order), this short circuit does not work, so we have to do yet another computation.

\comment
Since $z_i$ exhaust $F_q^*$, we may choose $i,j$ \st $1 - z_j/z_i$ is not a square. The corresponding matrices are  
 $$
\( \matrix \lambda_i & z_i \\ 0 & \lambda_i \endmatrix\) \quad\text{and}\quad \( \matrix \lambda_j& z_j \\ 0 & \lambda_j\endmatrix\).
$$
Set $s = (\lambda_i - \lambda_j)^2/(z_j-z_I)$ and $t = -z_i/s$; then $t = z_i^2 (1-z_j/z_i)$, and thus is not a square. Finally set $r = -\lambda_i$, and define  $e = \( \smallmatrix  0 & t \\ 1 & 0\\ \endsmallmatrix\)$ and the field  $D$ generated by $e$ as we did before. It is a two-line computation to that $d : = r\I + se$, when added to each of the matrices, yields noninvertible ones, and we are done with the case that $q$ is odd.

Now we show that for some $i < j$, for a suitable choice of $e$ with irreducible characteristic polynomial, and  $e $ and $D = F_q [e]$, that $\II_i^c \cap \II_j^c \neq \emptyset $. 
\endcomment

So assume $q$ is even. We wish to proceed as before, to show $\II_i ^c \cap \II_j^c$ (computed \wrt a to be determined maximal subfield $D$) is nonempty. We reduce to the following case (so $i = 2$), the matrices
$$
\( \matrix 1 & 1 \\ 0 & 1 \\\endmatrix\) \quad {\text{and}} \quad \( \matrix  \lambda_j & k_j \\  0 & \lambda_i\\\endmatrix\). 
$$
Here $\lambda_2 = k_2 = \lambda_2 k_2= 1$, and so none of $\lambda_j, k_j, \lambda_j k_j$ is $2$. Since ${\lambda_j}_j \geq 3$ consists of $q-2$ elements distinct from $\brcs{0,1}$, we may find $j$ \st $\lambda_j \notin B(F_q)$ (there are $q/2 -1$ choices for such, and since $q \geq 4$, at least one exists. Abbreviate $\lambda_j$ to $\lambda$, and $k_j$ to $u$, so I don't have to keep typing subscripts. We try $d = \(\smallmatrix \lambda & u \\ c & h\endsmallmatrix \)$ where $c,h$ are parameters so that $d$ plus both of the two matrices is noninvertible, and so that $d$ has irreducible characteristic  polynomial; then $D = F_q[d]$ is a maximal subfield of $\M_2 F_q$ and $d \in \II_2 ^c \cap \II_j^c$ (the latter computed \wrt $D$. 

We clearly have $d $ plus the second one is noninvertible for every choice of $c,h$. We may choose $h \neq \lambda,0$. We see that $d$ plus the first one is not invertible if and only if $(\lambda + 1)(h+1) = c(u+1)$. Since $u, \lambda \neq 1$, we have $c =( (h+1)\lambda  + \lambda + 1)/(u+1)$. The characteristic polynomial of $d$ is $x^2 + (h+\lambda) x + \lambda h + uc$, and thus the its characteristic polynomial is irreducible if $(\lambda h + uc)/(h+\lambda)^2 \notin B$. We have 
$$\eqalign{
\frac{\lambda h + uc}{(h+1)^2} & = \frac{h \( \lambda + \frac{\lambda+1}{u+1}\) + \frac{u(\lambda +1)}{u+1}  } {h^2 + \lambda^2}; \quad\text{replace $h \mapsto h+\lambda$; we obtain} \cr
\frac{h\(\frac {u+\lambda}{u+1}\) + \frac{\lambda^2 + u}{u+1}} {h^2}& := \frac{h\rho + \sigma}{h^2}. \cr
}$$
This defines $\rho$ and $\sigma$. If $\sigma = 0$ (that is, $\lambda^2 = u$), the result is $\rho/h$; as $u \neq \lambda$ (otherwise, one matrix would be a scalar times the other), $\rho \neq 0$. As $h$ varies over $F_q  \setminus \brcs{0,\lambda}$, there are at least $q/2 - 1$ choices for $h$ so that $\rho/h$ is not in $B$.

If $\sigma \neq 0$, let $\alpha$ denote the (unique) element of $F_q^*$ \st $\alpha^2 = \sigma$. Then we have
$$\eqalign{
\frac{h\rho + \sigma}{h^2} & = \(\(\frac{\alpha}{h}\)^2 + \frac{\alpha}{h}\) + \frac{\rho + \alpha}{h} \cr
& \equiv    \frac{\rho + \alpha}{h}   \mod B.
}$$
If $\rho + \alpha = 0$, then $\rho^2 = \sigma$, which converts to $(u+\lambda)^2 = (u+1){\lambda^2 + u}$. This simplifies to $\lambda^2 = u^2 + u$, implying that $\lambda^2 \in B$. But $B$ is closed \wrt squaring and taking square roots, so that $\lambda \in B$. This contradicts our initial selection of $\lambda$, so that $\rho + \alpha \neq 0$. Then there are (again) $q/2 - 1$ (that is, at least one) choices for $h$. 
\qed

It is highly likely that $\M_3 F_q$ (and thus all $\M_n F_q$ for $n \geq 2$) satisfies \gui {q+1}, but the inclusion-exclusion methods used here seem to have reached their limit. For $q = 3$, and $n =2$, the result is sharp (by Lemma \Antn).

\SecT Appendix C 

Here we show the following. 

\Lem Proposition \Bone. For all $n \geq 3$, $\Mn _n F_2$ satisfies \gui 3. 

This is a finite problem, since it is sufficient to show this only for $n = 4,5,6$, by Corollary \Anin. However, I am totally inept with computational software, so did it by hand. There are a disturbing number of cases.

Here is a straightforward variation on the proof of Proposition \Aeig. 

\Lem Lemma \Btwo. Let $R$ be an ring, and let $e$ be an idempotent of $R$. Suppose that $b, c$ are elements of $R$ \st $c = ec$, and there exists a relatively invertible $u_0$ in $eRe$ \st $u_0 + ece$ and $u_0 + ebe$ are both relatively invertible, and in addition, there exists relatively invertible $v_0$ in $(1-e)R(1-e)$ \st $v_0 + (1-e)b (1-e)$ is relatively invertible. Then there exists invertible $u$ \st both $u + b$ and $u +c$ are both invertible in $R$.

\Rmk As a consequence, if $c = ec$, $eRe$
 satisfies \gui {3}, and $(1-e) R(1-e)$ satisfies \gui 2, then for all $b$ in $R$, there exists invertible $u$ in $R$ \st $u+b$ and $u+c$ are invertible. 

\Pf Set $u = u_0 + v_0 - eb(1-e)$. If $u' = \overline{u}_0 + \underline{v}_0$, then $u'$ is clearly invertible and $u'u = 1 - (u_0 + v_0)eb(1-e) = 1 - eu_0 eb(1-e)$, and this is of the form $1 + n$ where $n^2 =0$, so is invertible. Similarly, $uu' = 1 - n'$ where $n' = eb(1-e)v_0(1-e)$. So $u$ is invertible. 

Next $u + c = (u_0 + ece) + v_0 - e(b-c(1-e))$, and we do the same process to show $u + c$ is invertible. 
Finally, 
$$\eqalign{
u + b &= u + ebe + eb(1-e) + (1-e)be + (1-e)b(1-e) \cr
& = (u_0 + ebe) + (v_0 + (1-e)b(1-e)) + (1-e)be.\cr
}$$
Multiplying on the left or right by $\overline{(u_0 + ebe)} + \underline{(v_0 + (1-e)b(1-e))}$ (which is clearly invertible) yields an element of form $1 + m$ where $m^2 = 0$, so once again, $u + b$ is invertible. \qed

\long\def\mat #1.{\(\matrix #1 \endmatrix \)}

Given $B, C \in \Mn_n F_2$ (where $n =3,4,5$), find $U \in \gl (n,F_2)$ \st $U + B$ and $U + C$ are invertible.  $\I$ will denote the identity matrix of whatever size we are currently working with (that is, $n$). And $\I_k$, $O_k$ will the denote the identity, zero matrices respectively, of size $k$. 

Let $\jnf (k)$ denote the Jordan block of size $k$ for the element $1$ (thus $\jnf (2) = \(\smallmatrix 1 & 1 \\ 0 & 1 \\\endsmallmatrix\)$), and let $\jnf_0 (k)$ be the size $k$ Jordan block for $0$.

We denote by $N_2 := \(\matrix 0 & 1 \\ 1 & 1 \\ \endmatrix\)$, the companion matrix of the one quadratic irreducible polynomial over $F_2$, $x^2 + x + 1$. It is conjugate (within $\gl(2,F_2)$) to its square.  The ring generated by $N_2$ is isomorphic to $F_4$.

\Lem Observation \Bthr. Let $C_0 \in \Mn_2 $ be one of 
$$
\brcs{\I_2,N_2,N_2^2, \mat 0 & 0 \\ 0 & 1. , \mat 1 & 0 \\ 0 & 0.  }
$$
Then $N_2 + C_0$ is invertible for the first, third, and fourth choices, while $N_2^2 + C_0$ is invertible for the first, second, and fifth choices. 
\vskip 10pt
We denote by $N_3$ (parameterized  by $a,b$), companion matrices
$$
N_3 = \mat 0 & 0 & 1 \\ 1 & 0 & b \\   0 & 1 & a \\  . .
$$
The characteristic polynomial of $N_3$ is irreducible precisely when $a \neq b$. Similarly, $N_4$ (parameterized by $a,b,c$) and $N_5$ (parameterized by $a,b,c,d$) are defined. We note that $N_k + \I$ is invertible if and only if the sum of the parameters ($a+b$, $a+ b + c$, \dots) is $1$. 

\comment 
We also observe that $\(\smallmatrix 0 & 1 \\  1 &0 \\ \endsmallmatrix\)$ is conjugate to $\(\smallmatrix 1 & 1 \\  0 &1 \\ \endsmallmatrix\)$ over $F_2$. 
\endcomment

Given $B$ and $C$ in $\Mn_n F_2$, we wish to find invertible $U$ \st both $U + B$ and $U+C$ are also invertible. 

\Lem Observation \Bfou. If $B = \I$ and the characteristic polynomial of $C$, $f$,  has no linear factors (that is, neither zero nor one are eigenvalues of $C$, then the unital subalgebra generated by $C$, $F_2 [C]$ is commutative and does not factor onto $F_2$ (since $F_2 [C] \iso F_2  [X]/(f)$), and being a finite-dimensional algebra (with simple factors $F_{2^n}$, for some value(s) of $n > 1$) thus satisfies \gui 3. In particular this applies to the pair $\I, C$. 

\noindent {\bf First case:} $\rk B  = n$. Here we can assume that $B = \I$ by the usual argument, with  $C$ updated.

\item{($\alpha$)} $C$ is invertible (and $B = \I$)

\item{($n=3$)} (i)   This means $0$ is not an eigenvalue, so we may assume that $1$ is (otherwise, \Bfou{} applies). Since we can conjugate $C$ by any element of $\gl (3,F_2)$ without affecting $B = \I$, we reduce to the following situations (we can of course discard the possibility that $C = \I$):
$$\matrix
C {\text{ is one of }} & N_2 \oplus (1) , &\mat 1 & 1 & 0 \\ 0 & 1 & 1 \\ 0 & 0 & 1 \\ . ,& \mat 1 & 1 & 0 \\ 0 & 1 & 0 \\ 0 & 0 & 1 \\ .  \\
{\text{Set }} U = & \mat  1 & 0 & 1\\ 1 & 1 & 1 \\ 0 & 1 & 1\\   . ,  & \mat  1 & 1 & 1\\ 0& 1 & 1 \\ 1 & 1 & 0\\   ., & \mat  0 & 0 & 1\\ 1 & 0 & 0 \\ 0 & 1 & 1\\   . \\
\endmatrix
$$

\item{($n=4$)} Again, we can assume that $1$ is an eigenvalue of $C$. This leaves us with the following possibilities: $N_3 \oplus (1)$, $N_2 \oplus \I_2$, $N_2 \oplus \( \smallmatrix 1 & 1 \\ 0 & 1 \\\endsmallmatrix\)$, and then the various Jordan normal forms associated to the eigenvalue $1$ having algebraic multiplicity four. 

We have the  the $4 \times 4$ companion matrices (parameterized by $a,b,c$ in $F_2$), 
$$
N_4 = \( \matrix 0 & 0 & 0 & 1 \\ 
1 & 0 & 0 & c \\ 
0 & 1 & 0 & b \\
0 & 0 & 1 & a \endmatrix\).
$$
 Its characteristic polynomial is $f(x) = x^4 + a x^3 + b x^2 + c x + 1$ (no need to deal with minus signs, since we are working in $F_2$!). If $a + b + c = 1$, then  $f(1) \neq 0$, that is, $1$ would not be an eigenvalue, and thus $W + \I$ would be invertible. The characteristic polynomial is irreducible if additionally, $(a,b,c) \neq (0,1,0)$.

\bu  $C = N_3 \oplus (1)$. There are two forms for $N_3$, one with trace $0$ and one with trace $1$; however, each is conjugate to the inverse of the other; since the pair $(\I,C)$ can be transformed into $(\I,C^{-1})$ (multiply the former by $C^{-1}$ and interchange) without affecting whether \gui 4{}
 holds), so we can assume $\tr N_3 = 1$.  
With $a = b = c = 1$, we check that $U = W^3$ satisfies $U + (N_3 \oplus (1)) $ is invertible. 

\bu  $C = N_2 \oplus I_2$. Using the observation, we see that we can set $U = N_2^2 \oplus N_2$.

\bu  $C = N_2 \oplus \jnf (2)$. Set $ U = W^2$ with $(a,b,c) = (1,1,1)$.

\bu  $C = \jnf (2) \oplus \jnf 2$. Set $B = W^2$ with $a= 1$ and $b = c =0$. 

\bu  $C= \jnf (4)$. Set $U = W$ with $a = 1$ and $b = c =0$. 

\bu  $C = \jnf (3) \oplus (1)$. Set $U = W$ with $a= 1$ and $b = c = 0$. 

\bu  $C = \jnf (2) \oplus \jnf (2)$.  Set $U = W$ with $a= 1$ and $b = c = 0$. 

\bu  $C = \jnf (2) \oplus\I_2$.  Set $B = W$ with $a= 0$ and $b = 1$, and $c = 0$.  This choice does not have irreducible characteristic polynomial (it is $(x^2 + x + 1)^2$), but since it doesn't have $1$ as an eigenvalue, $W + \I$ is also invertible. 

\item{($n= 5$)} $C$ still invertible. As before, we can assume $1$ is an eigenvalue.

\bu  $C = N_4 \oplus (1)$ or $\(\smallmatrix N_2 & X \\ 0_{2} & N_2 \endsmallmatrix\) \oplus (1)$.  The characteristic polynomials of the choices for $N_4$ are $f(x) = x^4 + s x^3 + tx^2 + u x + 1$; from  $1$ not being an eigenvalue, we have $A + B + C = 1$; being irreducible, $(s,t,u) \neq (0,1,0)$. 

On the other hand, the characteristic polynomial of the second matrix is $f(x) = (x^2 + x +1)^2$, which is precisely the case $(s,t,u) = (0,1,0)$, and the companion matrix of this polynomial is the second matrix (with appropriate choice of $X$). 

So we take $C = C_0 \oplus (1)$, where $C_0$ is the companion matrix of $f$ (and let the parameters $s,t,u$ vary subject to the one constraint, so there are four possible matrices. We can reduce this to three, by observing that if $s\neq u$, the the two corresponding matrices are conjugate to the inverse of each other, and the problem $\I,C$ is solvable for \gui {3} if and only if $\I, C^{-1}$ is (multiply both by $C^{-1}$ and interchange them). So we can assume $(s,t,u) \in \brcs{(0,1,0), (0,0,1), (1,1,1)}$. 

We use a matrix conjugate to the  generic invertible companion matrix (as it seemed more likely to give results relatively quickly), with parameters $a,b,c,d$ subject to $a+ b + c  + d = 1$ (to guarantee that $1$ is not an eigenvalue, so that $W + \I$ is invertible). Define 
$$
W = \( \matrix a & 1 & 0 & 0 &0 \\
b & 0 & 1 & 0 & 0 \\
c & 0 & 0  & 1 & 0\\
d & 0 & 0 & 0 & 1 \\
1 & 0 & 0 & 0 & 0 \\
  \endmatrix \)
$$

For $C$ with $(s,t,u) = (1,1,1)$, set $c = 1$ (and and any choice of $a,b,d$ in $F_2$ subject to $a+b + d = 0$); then $U = W$. If $(s,t,u) = (0,0,1)$, set $d = 1$ (and choice of $a,b,c$ with $a+ b + c = 0$, and set $U$ to be the corresponding $W$. 
Finally, if $(s,t,u) = (0,1,0)$, set $d = 0$ (and $a,b,c$ with $a + b + c =1$, with $U$ the corresponding $W$. 

\bu  If $C_0 = \I_2$ or $N_2$ is a matrix direct summand of $C$, then we can write $C = C_0 \oplus C_1$, where $C_1$ is $3\times 3$. From the results in the case that $n = 3$, there exists $U_1 \in \gl(3,F_2)$ \st $U_1 + \I_3$ and $U_1 + C_1$ are both in $\gl(3,F)2$; set $U = N_2^2 \oplus U_1$. 

\bu  $C = \jnf (5)$. In this case, we replace $C$ by its transpose (which is conjugate to it), and set $U = W$, with any choice of $a,b,c,d$ \st $a = c $ and $b+ d = 1$. 

\bu  $C= \jnf (4) \oplus (1)$. We can write this as a block upper triangular matrix with block sizes $3$ and $2$, $\(\smallmatrix A & X \\ 0_{2\times 3} & \I_2\\ \endsmallmatrix \)$, where $A = \jnf(3)$ (and $X$ doesn't matter). From the size three case, there exists $U_0$ \st $U_0 + \I_3$ and $U_I + A$ are both in $\gl  (3,F_2)$, and we set $U = U_0 \oplus N_2$. We see that $U + C$ is a block upper triangular matrix, with invertibles along the diagonal, hence is invertible. 

\bu  $C = \jnf(2) \oplus \jnf (2) \oplus (1)$. Similar to the previous case, there is a block upper triangular form $\(\smallmatrix A & X \\ 0_{2\times 3} & \I_2\\ \endsmallmatrix \)$, where $A = \jnf (2) \oplus (1)$, so the same argument applies. 

\bu  $C = \jnf (3) \oplus \jnf (2)$. Set $a = 0$ (and $b + c + d = 1$), and let $U$ be the corresponding $W$ (this one is particularly easy). 

\item{($\beta$)} $B $ is invertible, $C$ is not.

\noindent $\rk C = n-1$. 

\item{($n = 3$)} Since we can conjugate by any $g \in \gl(3,F_2)$, there are only a few possibilities with $C$ of rank two: $N_2 \oplus (0)$, $\I_2 \oplus (0)$, $\jnf_0 (3)$, $\inf_0 (2) \oplus (1)$, and $\jnf (2) \oplus (0)$.

\bu  $C =N_2 \oplus (0) $. Set $U = N_3^3$ with $a = 0$ and $b = 1$. 

\bu  $C=  \I_2 \oplus (0)$. Set $U=N_3^2$ with  $a = 0$ and $b = 1$. 

\bu  $C= \jnf_0 (3)$. Set $U=N_3$ with  $a = 0$ and $b = 1$. 

\bu  $C= \jnf_0 (2) \oplus (1)$. Set $U=N_3$ with  $a = 1$ and $b = 0$. 

\bu  $C=\jnf (2) \oplus (0)$. Set $U=N_3$ with  $a = 1$ and $b = 0$.

\item{($n = 4$)} The possibilities for $C$ are $N_3 \oplus (0)$, $N_2 \oplus \jnf_0 (2)$, $N_2 \oplus (1) \oplus (0)$, $\jnf_0(4)$, $\jnf_0 (3) \oplus (1)$, $ \jnf_0(2) \oplus \I_2$, $\inf_0 (2) \oplus \jnf (2)$, $\I_3 \oplus (0)$, $\jnf (3) \oplus (0)$, $\jnf (2) \oplus (1) \oplus (0)$. 

\bu  $C = N_3 \oplus 0$; there are two choices for $N_3$, one with trace zero, the other with trace $1$). Set $U = N_4$ with $(a,b,c) = (1, \tr N_3, \tr N_3)$. 

\bu  $C = N_2 \oplus \jnf_0 (2)$. Set $U = N_4^2$ with $a= b$ and $c = 1$.

\bu  $C \in \brcs{N_2 \oplus (1) \oplus (0), \I_3 \oplus (0), \jnf (2) \oplus (1) \oplus (0) }$. In each of these cases, we can take $U = N_2 \oplus N_2^2$, because of their $2 \times 2$ block structure.

\bu  $C = \jnf (3) \oplus (0)$. This can be written as block upper triangular matrices (with blocks of size $2$), with diagonal blocks, $\I_2$ and $\diag (1,0)$. Thus $U = N_2 \oplus N_2^2$ will do.

\item{($n = 5$)}.We can conjugate by any element of $\gl(5,F_2)$; this does not affect $B = \I$. 

\bu  $C= N_4 \oplus 0$ and also with $N_4$ replaced by the companion matrix of $(x^2 + x + 1)^2$. There are three parameters $s,t,u$ in the definition of $N_4$, and four, $a,b,c,d$ in the definition of $W$, with constraint $a + b + c + d = 1$.

If $s = 1$, set $(a,b,c,d) =( 1, 1+u,1 + u+ s,s)$, and let $U= W^2$. If $u = s \neq 0$, set $(a,b,c,d) = (0,0,1,0)$ and $U = W^2$. If $u = 0$ and $s= 1$, set $(a,b,c,d) = (0,0,0,1)$ and $U = W^2$. 

\bu  $C = N_3 \oplus \jnf_0 (2)$. There are two choices for $N_3$, but each is conjugate to the inverse of the other, so we can assume that the trace is $1$. Set $U = W^2$ with $a \neq b$ (and with $c + d = 0$).

\bu  If $N_2$ or $\I_2$ is a direct summand of $C$, say $C = C_0 \oplus C_1$ with $C_1$ being $N_2$ or $\I_2$, then we just apply the $3 \times 3$ result and matrix direct sum it with $N_2^2$. The same thing can be done if $(0) \oplus (1)$ is a matrix direct summand of $C$.   This reduces the number of cases substantially.

\bu   $C = \jnf_0 (5)$. Here we work instead with $C = \jnf_0 (5)^T$ (which of course is conjugate to $\jnf_0 (5)$). Take $U = W$ with $a = c$ (so $b + d = 1$).

\bu  $C = \jnf_0 (4) \oplus (1)$. This is a block upper triangular matrix with the lower right block being $(0) \oplus (1)$. If $C_0$ is the upper diagonal block, there exists $U_0$ in $\gl(3,F_2)$ \st $U_0  + \I_3$ and $U_0 \oplus C_0$ are both invertible, so that we simply set $U = U_0 \oplus N_2$;

\bu  $C = \jnf_0 (3) \oplus \jnf (2)$. Since we have solved the $n= 3$ case for $\I_3, \jnf_0 (3)$, say with $U_0$, and set $U= U_0 \oplus N_2$. 

\bu  $C =\jnf_0 (2)\oplus \jnf (3)  $. We solve the $n=3$ case for $C_0 = \jnf_0 (2) \oplus (1)$, giving $U_0$, and set $U = U_0 \oplus N_2$.

\bu  $C = (0) \oplus \jnf (4)$ or $(0) \oplus Z$ where $Z$ is any upper triangular matrix with  $1$ in its upper left position. This is block upper triangular, with the $2 \times 2$ block being $\diag (0,1)$; so we can set $U = N_2 \oplus U_0$, as in previous examples.

\bu  $C = \jnf (2) \oplus \jnf_0 (2)  \oplus  (1)$ (note the order in which they are written). This is block upper triangular with the lower $2 \time 2$ block being $\diag (0,1)$, and the by now standard technique applies. 

\noindent $\rk C \leq n-2$ (and $B = \I$). If $n \in \brcs{4,5}$, then Lemma \Btwo{} applies regardless of the choice of $B$. 
Hence we are reduced to $n = 3$ and $\rk C = 1$. Up to conjugation, there are only two possibilities for $C$.

\bu  $C = (1) \oplus 0_2$. Set $U = N_3$ with $\tr N_3 = 1$ (that is, $a=1$ and $b = 0$).

\bu  $C = \jnf_0 (2) \oplus (0)$. Set $U = N_3$ with $a=0$ and $b = 1$.

\noindent {{\bf Second case}} $\rk B = n-1\geq \rk C$. 

\item{($n=3$)}
 There exist $V,W $ in $\gl(3,F_2)$ 
\st $VBW = \I_2 \oplus (0)$; relabel the pair $(VBW,VCW)$ as $B = \I_2 \oplus (0)$ and $C$. We observe that $B$ is left invariant under conjugation by elements of the form $g \oplus (1)$ for $g \in \gl(2,F_2)$ and also the elementary operations of adding the third row to any of the others, and similarly adding the third column to the others.

Now assume $\rk C = 2$. 

\bu  $C_{3,3} = 1$. In this case, we can add if necessary, the third row to the others, and the third column to the others, reducing to (new) $C = C_0 \oplus (1)$, where $\rk C_0 = \rk C -1$. Thus, as $\rk C= 2$, $\rk C_0 = 1$, and by conjugating with suitable $g \oplus (1)$ for $g$ in $\gl (2,F_2)$, we reduce to the case that $C_0$ is either $(1) \oplus (0)$ or $\jnf_0 (2)$. 

\bu  For $C = (1) \oplus (0) \oplus (1)$ [and $B = \I_2 \oplus (0)$], set $U =N_3^4 $ with $a = 0$ and $b =1$. 

\bu  $C = \jnf_0 (2) \oplus (1)$, set $U = N_3^5 $ with $a = 1$ and $b = 0$. [This was the most worrying case of all.]

\bu  $C = C_0 \oplus (0)$. By conjugation, we may assume either $C_0 = \I_2$,  $C = N_2 \oplus (0)$, or $C = \jnf (2)$. 

\bu  If $C = \I_2 \oplus (0)$ or $N_2 \oplus (0)$, set $U = N_2^2 \oplus (1)$. 

\bu  If $C = \jnf (2)$, set $U = N_3^4$ with $a=1$ and $b =0$. 

\bu  If $C_{3,3} = 0$, but there exists $i$ with $C_{3,i} = 1$. Adding the $i$th column, if necessary to the $3-i$th, we can assume that $C_{3,3-i} =0$; by interchanging the first and second rows and columns (that is, conjugating with $\(\smallmatrix 0 & 1 \\ 1 & 0 \\ \endsmallmatrix\) \oplus (1)$, we can assume $i = 2$. Then adding the third row to the first and second rows (in necessary), we can assume that the second column
is $\(\smallmatrix 0 & 0 &1  \\ \endsmallmatrix\)^T$. Now we have two subcases. 

\noindent $\Delta$  The current $C$ has only zeros in the third column. In this case, 
$$
C = \(\matrix c_1 & 0 & 0\\ c_2 & 0 & 0 \\ 0 & 1 & 0 \\    \endmatrix \)
$$
where at least one of the $c_i$ is $1$. 

\bu  If $c_2 = 0$, set $U = N_3^2$ with $a = 0$, $b= 1$.

\bu  If $c_1 + c_2 = 0$, set $U= N_3^3$ with $a= 1$, $b = 0$. 

\bu  If $c_1 = 0$, set $U = N_3^5$ with $a = 1$, $b =0$.

\bu If $C_{i,3} = 1$ for some $i$, but $C_{3,3} = 0$. We can subtract the $i$th row from the $(3-i)$th, and so assume $C_{3-i,3} = 0$. Thus we are reduced to the case that 
$$\eqalign{
C &= \(\matrix d & 0 & 0\\ 0 & 0 & 1 \\ 0 & 1 & 0 \\    \endmatrix \) \quad {\text {OR}}\cr
C & = \(\matrix 0 & 0 & 1\\ e & 0 & 0 \\ 0 & 1 & 0 \\    \endmatrix \) \cr
}$$
But since $\rk C = 2$, we must have $d = e = 0$. For the second one, set $U$ to be the permutation matrix 
$\(\smallmatrix 0 & 1 & 0 \\ 0 & 0 & 1 \\ 1 & 0 & 0\\ \endsmallmatrix \)$ and for the first, take $U$ to be the transpose of the permutation matrix.

\bu  Bottom row consists of zeros, but $C_{i,3} \neq 0$ for some $i \in \brcs{1,2}$. The transpose can be applied here (this preserves invertibility and addition, and leaves $B$ unchanged), resulting in a case done above. 

\noindent   $\rk C  = 1$. 

\bu If $C_{3,3} = 1$, as before, we reduce to the case that all the other entries in the third row and column are zero, and now rank equalling $1$ entails $C = 0_2 \oplus (1)$. Then we can take $U$ to be the permutation matrix with ones  in  the $(1,2)$, $(2,3)$, and $(3,1)$ positions. 

\bu  If $C_{3,i} = 1$ for some $i \in \brcs{1,2}$, then we proceed as above (adding the third row to the others if necessary), and we obtain the possibilities, $C$ consists of zeros except possibly in the $(3,1)$ and $(3,2)$ positions---so the reduction is to three matrices, one of which has two nonzero entries. In any event, set $U$  to be $N_2 \oplus (1)$, which works for all three. 

\bu  If $C = C_0 \oplus (0)$, we can conjugate $C_0$ (by $\gl(2,F_2)$ to either $\jnf_0 (2)$ or $(1) \oplus (0)$. 

\bu  If $C = \inf_0 (2) \oplus (0)$, set $U$ to be the permutation matrix with nonzero entries in positions $(1,3), (2,1), (1,2)$.

\noindent  $\rk B = 1 = \rk C$. Left and right multiplying by invertibles, we can assume $B = (1) \oplus O_2$, and relabel the modified $C$, $C$. We can add either of the second or third rows (columns) to each other, and also to the first, and we can also conjugate by elements of $(1) \times \gl (2,F_2)$.  So we reduce to $C = (0) \oplus (1) \oplus (0)$ and $C = \jnf_0 (2) \oplus 0$. In the first case, set $U = \(\smallmatrix 0 & 1 \\ 1 & 0 \\\endsmallmatrix \)$, and in the second, just take the $90^{0}$ rotation of the identity matrix. 

This completes the proof of the proposition in the case that $n=3$. 

For both $n = 4$ and $n = 5$, we can apply Lemma \Btwo, and thus there is no problem if $\min \brcs{\rk B, \rk C } \leq n-1$. Hence in both these situations, it suffices to deal with the cases that $\rk B = \rk C = n-1$. 

\item{($n = 4$)} $\rk B= \rk C = 3$. As usual, we can assume $B = \I_3 \oplus 0$. We go through the same process as in the case $n=1$. 

\bu If $C = C_0 \oplus (0)$, there exists (by the result for $n = 3$) an element $U_0$ in $\gl(3,F_2)$ \st $U_0 + \I_3$ and $U_0 + C_0$ are both invertible (in $\Mn _3 (F_2)$). Set $U = U_0 \oplus (1)$. 

\bu If $C_{4,4} = 1$, we reduce (as in the case $n = 3$) to $C = C_0 \oplus (1)$, and now $C_0$ has rank two. 
We can conjugate with any element of the form $g \oplus (1)$  for $g \in \gl(3,F_2)$ (since $B = \I_3 \oplus 0$ is left unchanged), and thus have the following possibilities, $N_2 \oplus (0)$, $\I_2 \oplus 0$, $(1) \oplus (0) \oplus (1)$, $\jnf_0 (3)$, and $\jnf_0 (2) \oplus (1)$. 

\bu If $C = \I_2 \oplus 0 \oplus (1)$ or $(1) \oplus (0) \oplus (1)\oplus (1)$ or $N_2 \oplus (0) \oplus (1)$, we note the $2$-block diagonal structure and set $U = N_2^2 \oplus N_2$. 

\bu If $C = \jnf_0 (2) \oplus (1) \oplus (1)$, we conjugate (using $\gl(3,F_2) \times 1$) this to $(1) \oplus \jnf_0 (2) \oplus (1)$. This has the form of an upper triangular block matrix, and we can simply take $U = N_2^2 \oplus \(\smallmatrix 0 & 1 \\ 1 & 0 \\\endsmallmatrix\)$. 

\bu If $C = \jnf_0 (3) \oplus (1)$. Set $U = N_4$ with $(a,b,c) = (0,0,1)$. 

\bu  $C_{4,4} =0$. Then we can apply \Btwo, with 
$e = \I_3 \oplus (0)$, 
as we now know that $eRe = \Mn _3 F_2$ satisfies \gui 3; from $C_{4,4} = 0$, we have  $(1-e) C (1-e) = 0$, and obviously $(1-e)B(1-e) = 0$ (so we can set $V_0 = (1)$).

This concludes the proof for $n = 4$.

\item {($n= 5$)} We may assume $B = \I_4 \oplus (0)$. If $C_{5,5} =1$, adding the fifth column/row to  other columns/rows, we reduce to the case that $C = C_0 \oplus (1)$, and $C_0$ is rank three. We can conjugate $C_0$ by elements of $\gl(4,F_2)$, and so we further reduce to $C_0$ is one of 
$N_3\oplus 0$,  $N_2 \oplus \jnf_0 (2)$, $N_2 \oplus (1) \oplus (0)$$I_3 \oplus (0)$, $\jnf (3) \oplus (0)$, $\jnf (2) \oplus (1) \oplus (0)$, $I_2 \oplus \jnf_0 (2)$, $\jnf (2) \oplus \jnf_0 (2)$, and $\jnf_0 (4)$.
 
Now we observe that for all choices of $C_0$, the lower $2 \times 2$ block of $C$ (that is, $\(\smallmatrix C_{4,4} & C_{4,5}\\ C_{5,4} & C_{5,5}\\ \endsmallmatrix\)$) is  $\diag (0,1)$, and that of $\I_4 \oplus 0$ is $\diag (1,0)$, so that $\( \smallmatrix 0 & 1 \\ 1 & 0 \\ \endmatrix\)$ added to each will result in invertible matrices. Since $\Mn_3 F_2$ satisfies \gui 3, we can can apply Lemma \Btwo{} with $e = \I_3 \oplus 0_2$. 

Finally, suppose $C_{5,5} = 0$. Apply \Btwo{} with $e = \I_4 \oplus (0)$ (possible, since  $(1-e)B(1-e) = (1-e)C(1-e) = 0$). 
\qed  

Thus $\Mn _n F_2$ satisfies \gui {3} for all $n \geq 3$. How sharp is this? We know  from the rather flabby estimate (Lemma \Antn) that $\Mn _n F_2$ fails to satisfy \gui {n+1}. So $\Mn _3 F_2$ fails \gui {4}; but all we obtain for $n = 4$ is that $\Mn _4 F_2$ fails \gui {5}. So it is possible that it satisfies \gui 4. The amount of effort put into showing it satisfies \gui {3} suggests that it likely doesn't satisfy \gui 4.

\comment
If $R$ is Hermite, so is $M_n R$ 3.6

3.7 R commutative Hermite ring, let $(a_1, \dots, a_n)$; there exists n x n matrix with first row (a_1, \dots a_n) \st det of the matrix is $d$ and dR = \sum a_i R. 

3.5 R right Hermite, given $A$ there exists invertible U \st AU is lower triangular

3.4 Commutative $R$ w/o zero divisors is Hermite
iff sum and intersection of principal ideals is again principal. 

p 11 (474) PID (not necessarily commutative) is elementary divisor ring

Boolean algebra is direct limit of finite rings, etc. So M_n S satisfies \gui {3} iff $n \geq 3$
\endcomment

\SecT  Appendix D.  Failure of \gui 1.

Recall that a ring $R$ satisfies \gui 1{} if for all nonzero $a$ in $R$ \st $aR \neq R$, the intersection $aR \cap V(R)$ contains a nonzero element ($V(R)$ is the set of differences of invertible elements of $R$). We showed (Lemma \Atwh) that \gui 1 plays a (small) role in limiting properties of matrix rings. It is straightforward to show that   if $K$ is an algebraic extension of a finite field, then the polynomial ring with any nonconstant polynomial $P$ inverted, $K[x,P^{-1}]$, satisfies \gui 1. However, it turns out that for all other fields $K$, $K[x,P^{-1}]$ {\it fails\/} to satisfy \gui 1\ for all choices of nonconstant $P$. This amounts to inverting finitely many primes (irreducibles), and in some cases (e.g., if $K$ is uncountable), the same is true if we invert any countable set of primes in $K[x]$, and even some uncountable ones. 

If $W $ is a collection of irreducible elements of $K[x]$, we denote by $K[x,W^{-1}]$ the ring obtained from $K[x]$ by inverting all the elements of $W$. We   assume that $W$ consists of distinct monic polynomials. 

First we deal with the easy case. 

\Lem Lemma \Atwf. Suppose that $S$ is a commutative ring \st that all proper factor rings are finite, and in addition possesses a unit of infinite order. Then $S$ satisfies \gui 1. 

\Pf Let $u$ be a unit of infinite order, and form $J = aS$. As $S/J$ is finite, there exists an integer $n$ \st $u^n - 1$ belongs to $J$. As $u$ is of infinite order, $u^n -1$ is not zero. \qed

\Lem Corollary \Ctwo. All of the following rings satisfy \gui 1.
\item{(a)} $K[x,W^{-1}]$ for any nonempty set  of primes $W$ where $K$ is an algebraic extension of a   finite field.
\item{(b)}  Any order in a number field other than $\Q$ or a quadratic imaginary field.
\item{(c)} $\Z[1/n]$ for any positive integer $n$ exceeding $1$. 

\Pf If $K$ is finite, (a) follows from Lemma \Atwf. Every algebraic extension of a finite field is a direct limit of finite fields. As \gui 1{} is preserved by direct limits, (a) follows.

\noindent (b) All of these have an infinite order unit and all proper factor rings are finite, so \Atwf\ applies.

\noindent (c)  Lemma \Atwf\ applies. 
\qed

Now we will show the complementary result.

\Lem Proposition \Cthr. Let $K$ be a field that is not algebraic over a finite field. Then for all nonconstant  $P$ in $K[x]$, the ring $K[x,P^{-1}]$ does {\it not\/} satisfy \gui 1. 

We deal with the various cases separately. The strategies  are the same: given a nonconstant polynomial $P$, we factor it as $\prod P_i^{m(i)}$ in terms of irreducible polynomials, $P_i$ which we can assume monic, and we can also assume all the exponents are $1$. Then we pick a suitable element $a$  in $K[x]$ with appropriate zeros (depending on the choice of $P$), and obtain a contradiction by evaluating differences of units at the roots of $a$.  

The units of $K[x,P^{-1}]$ are precisely those elements of the form $\gamma \prod P_i^{C(i)}$ where $\gamma \in K\setminus \brcs{0}$ and $C (i)$ are integers (which can be negative). Hence a difference of units can be written in the form 
$$\eqalign{
\gamma_1 \prod P_i^{C(i)} - \gamma_2 \prod P_i^{D(i)} & = 
u \( 1 - \gamma\prod P_i^{A(i)} \)
}$$
where $u$ is a unit of $R_K = K[x,P^{-1}]$, $\gamma$ is a nonzero element of $K$, and the $A(i)$ are integers (not necessarily positive). We observe that if $z$ is an element of the algebraic closure of $K$ and is not a root of any of the $P_i$ (or equivalently, of $P$), then the assignment $x \to  z$ extends uniquely to a $K$-homomorphism $R \to \overline K$. Hence if $a \in K[x,P^{-1}]$ has zeros $\brcs{s_l}_{l=0}^k$  that are not zeros  of $P$, and $a = u\cdot (1- \gamma \prod P_i^{A(i)})$ as in the last display, then 
$$\eqalign{
\prod_i P_i^{A(i)} (s_l) &= \frac 1\gamma \quad \text{for all $l$, and thus} \cr
 \prod_i \(\frac{P_i(s_l)}{P_i(s_0)}\)^{A(i)}  &= 1 \quad \text{for all $l \neq 0$}.  \cr
}$$
The last form is valid in any field containing $\brcs{s_l}$,  and the coefficients of all the $P_i$ and of $a$. So if we can prove that no such set of equations hold in this field, then it fails in any overfield (with the same $a$, $P_i$, etc). But we can go a bit further. Let $K'$ be  a field containing $K$ and the coefficients of $P$ and the roots of $P$, and also suppose that $a$ was in $K$. If the equations are impossible in $K'$, then  $aR_{K'} \cap V(R_{K'}) = \brcs{0}$. But $aR_K \subset a R_{K'}$ and $V(R_{K} ) \subset V(R_{K'})$, so that $aR_K \cap a R_{K'} = \brcs{0}$, so $R_K$ fails \gui 1.

In all cases, $P$ is fixed, and we choose an $a$ depending on its prime factors.

\noindent {\it $K$ is algebraic over the rationals}

\noindent First, assume that $K$ is finite-dimensional over $\Q$.  Let $K'$ be the Galois closure of $K$ over $\Q$, and let $K_s$ be the splitting field of $P$ over $K'$; this is Galois over $\Q$. List the roots of $P$ in $K_s$, $\brcs{z_i}_{i=1}^m$. Let $\Z_s$ denote the ring of integers of $K_s$. Since the class group is finite, there exists a positive integer $r$ in $\Z$ (not merely in $\Z_s$) \st $T:= \Z_s[r^{-1}]$ is a principal ideal domain. 

We may thus write $z_i = \alpha_i/\beta_i$ where $\alpha_i, \beta_i$ belong to $T$ and are relatively prime (meaning, the ideal $(\alpha_i, \beta_i) = T$ for each $i$).

Pick an integer $n$ in $\Z$ that is relatively prime to $r$ \st $P(n) \neq 0$ (so in particular, $n \not\in \brcs{z_i}$).

Let $\Arrow N; K_s . \Q$ denote the norm. There exists a prime $p$ in $\Z$ that is relatively prime to $r$ as well as $N(n-z_i)$ for all $i$ and $N(z_i - z_j)$ for all $i \neq j$. Since the $n-z_i$ are fractions in $\Z_s$ this means that $p$ is relatively prime to each of the numerator and the denominator, that is to $N(n\alpha_i - z_i$ and $N(\alpha_i)$, and similarly to the numerator and denominator of all $z_i - z_j$. 

For each $i$, set $t_i = (p+1)z_i$ and $\delta = (x-n)\prod (x- t_)$. Let $a$ be the norm of this (the norm extended to the polynomial ring, $K_s [x]$, so that $a \in \Q[x] \subseteq K[x]$ and has roots (in $K_s$) consisting of $n,t_i$. 

If 
$$
 \delta K_s [x, W^{-1}] \cap V(K_s [x,W^{-1}]  = \brcs {0}
 $$ (where $W = \brcs{x-z_i}$), then it easily follows that $a K[x,P^{-1}] \cap V(K[x,P^{-1}] = \brcs{0}$. 

So suppose that the displayed intersection contains more than zero. Then there exists nonzero $q \in R_s = K_s [x, W^{-1}] $ \st 
$$\delta q = 
u\( 1 - \gamma\prod (x-z_i)^ {A(i)} \) .$$
By construction, the evaluation maps $x \mapsto n$ and $x \mapsto t_j$ extend to homomorphisms $R_s \to K_s$. The left side vanishes at each one of these, so we deduce 
$$\eqalign{
\prod_i (n-z_i)^{A(i)} &= \frac 1{\gamma} \cr
\prod_i (t_j - z_i)^{A(i)} & = \frac 1{\gamma} \quad \text{for each $j$, and thus} \cr
\prod_i \(\frac{t_j - z_i} {n-z_i}\)^{A(i)} & = 1 \quad \text{for each $j = 1,2, \dots, m$}.\cr
}$$

For $j$ fixed $t_j - z_j = p z_j$, whence $p$ divides $N(t_j - z_j)$, but for $i \neq j$, $t_j - z_i= pz_j + z_j - z_i$, the norm of which is rfelatively prime to $p$ (both denominator and numerator). This can only occur if $A(j) = 0$. Hence all the $A(i)$ are zero, forcing $a$ to be a unit, a contradiction.

Now assume that  {\it $K$ is algebraic over $\Q$}.
 (There is no reason to suspect that the property of not satisfying \gui 1{} is preserved by direct limits.) 
 
 Let $K'$ be the field generated over the rationals by the coefficients of $P$, and let $K''$ be the splitting field of $P$ over $K'$. The $a$ chosen above is in $\Q [x] $, and the equations cannot hold in $K''$, and thus in the compositum $KK''$,  $a KK''[x,P^{-1}] \cap V(K''[x,P^{-1}]) = \brcs{0}$, but this entails $a K[x,P^{-1}] \cap V(K[x,P^{-1}) = \brcs{0}$, and so $K$ fails \gui 1.
 \qed
 
 \def\RR{{\Cal R}}
 Next, we consider the case of transcendental extensions of the rationals or $F_p$. In some cases, we allow $W$ (a set of monic irreducible polynomials in $K[x]$ which we invert) to be infinite. 
 \comment
 For a monic irreducible polynomial $P$ in $W$ we define $\varepsilon (P) := (-1)^d P(0)$ where $d$ is the degree of $P$.  Then $\varepsilon (W)$ is defined to be $\Set {\varepsilon (P)}{P \in W}$, and 
 \endcomment
 Let $\RR (W)$ denote the set of all roots of elements of $W$ in a fixed algebraic closure of $K$; $F(\RR(W))$ denotes the field generated by $\RR(W)$. 
 
 If $K$ is a field, its {\it base field}, denoted $B(K) = B$, is $\Q$ if the characteristic is zero, and $F_p$ if the characteristic is $p$. 
 
 \Lem Lemma \Cfou. Suppose that $K$ is a field, $W$ is a set of monic irreducible polynomials in $K[x]$, and {\it either\/} of the following hold:
 \item{(a)}$ F(\RR(W)) \cap K  \not \subset \RR(W)$  and  there exists $t$ in $K$ that is transcendental over $B (\RR(W))$.
 \item{(b)} There exists a subset $\brcs{c,t} \subset K$ that is algebraically independent over $F(\RR(W))$. {\par} 
 \noindent Then $K[x,W^{-1}] $ fails to satisfy \gui 1. 
 
\Pf If (a) holds, select $c $ in $F(\RR(W)) \cap K \setminus \brcs{\RR(W)}$. In either case,  form $a = (x-c)(x-t) \in K[x]$. Assume that $a K[x,W^{-1}] \cap V(K[x,W^{-1}])$ contains a nonzero element. Then there exists a finite subset $\brcs{P_i} $ of $W$, together with $\gamma \in K\setminus \brcs{0}$, a unit $u$ of $K[x, \prod P_i^{-1}]$, and $q \in K[x]$
\st 
$$a q = u\cdot (1- \gamma \prod P_i^{A_i}). $$
We can evaluate at both $c$ and $t$ since neither is a root of any element of $W$. This yields (after dividing one relation by the other, to eliminate the $\gamma$)
$$
\prod \(\frac{P_i (t)}{P_i (c)}\)^{A(i)} = 1.
$$
The $A(i)$ are only integers (not necessarily nonnegative integers), so we set $I_+ = \Set{i}{A(i) > 0}$ and $I_-= \Set{i}{A(i)< 0}$, and obtain
$$
\prod_{i \in I_+} \(\frac{P_i (t)}{P_i (c)}\)^{A(i)} = \ \ \prod_{i \in I_-}\(\frac{P_i (t)}{P_i (c)}\)^{|A(i)|} .
$$
If (a) applies, we observe that $c \in F(\RR(W))$, so all of the denominators lie in $F(\RR (W))$; in either case, we obtain a contradiction.
\qed

If $F $ is a subfield of $K$, then $\tr\! \deg_F K$ ({\it transcendence degree\/}) denotes the cardinality of any maximal algebraically independent set over $F$ in $K$. 

\Lem Corollary \Cfiv. If $|W| +1 < \tr\! \deg_B K  $, then $K[x,W^{-1}]$ fails to satisfy \gui 1. In particular, if $W$ is finite and $\tr\!\deg_B K$ is infinite, then again $K[x,W^{-1}]$ fails \gui 1.

This includes the case that $W$ be infinite, e.g., if $K$ is uncountable, then the result holds for  countably infinite $W$. This contrasts starkly with Examples \Afou(ii, iia). 

There remains the final case. 

\noindent {\it $K$ has nonzero but finite transcendence degree over the base field $B$.}

\noindent There exists a finite subset $\brcs{y_1,\dots, y_n}$ of $K$ that is algebraically independent over $B$ and $K$ is algebraic over $Q:= B(y_1, \dots, y_n)$. First, assume that $K$ is finite-dimensional over $Q$ (we do not assume that the extension $Q \subset K$ is separable). 

Form the principal ideal domain $\Z_Q: = B(y_1, \dots, y_{n-1}) [y_n]$ (its main feature, that it has infinitely many prime ideals), and let $\Z_K$ be the integral closure of $\Z_Q$ in $K$. Write  $W= \brcs{P_1, P_2, \dots, P_m}$ (consisting of irreducible monic polynomials---in $x$---of $K[x]$). As before, $\RR(W)$ is the set of roots in a fixed algebraic closure of $K$, and let $K' \subset \overline {K}$ such that $K'$ a finite extension field of $K$ that contains $\RR(W)$; correspondingly, $\Z_{K'}$ will be the integral closure of $\Z_Q$ in $K'$; by transitivity of integrality, this is also the integral closure of $\Z_K$ in $K'$. 

For each $i$, let $S_i$ denote the {\it list\/} of roots of $P_i$, together with their multiplicities (since we permit nonseparable extensions, roots may have multiplicities, $m(s)$, divisible by the characteristic), so that $P_i$ factors as $\prod_{s\in S_i} (x-s)^{m(s)}$. 

If $r,s$ are elements of a commutative ring $A$, we say that {\it ${r,s}$ is relatively prime\/} or {$r$ is relatively prime to $s$\/} if $Ar + As = A$, that is, there exist $t,u$ in $A$ \st $tr + us = 1$. The use of this term does not imply any sort of unique factorization in $A$.

Since $\RR(W)$ is finite and $\Z_Q$ is not, there exists $b$ in $\Z_Q \setminus \RR(W)$ (that is, $b$ is not a zero of $P = \prod P_i$). There exists nonzero $\alpha$ in $\Z_Q$ \st for all $s \in \RR(W)$, $\alpha s \in \Z_{K'}$.

 Because the set of prime ideals of $\Z_Q$ is infinite, there exists an irreducible element $v$  of 
$\Z_Q$ (monic and irreducible as a polynomial  in $y_n$ with coefficients in $B(y_1, \dots, y_{n-1})$) with the following properties:
\item{(i)} $v$ is relatively prime (in $\Z_Q$) to $b$ in $\Z_Q$;
\item{(ii)} $v$ is relatively prime (in $\Z_{K'}$) to $s\alpha $ for each $s $ in $\RR(W)$;
\item{(iii)} $v$ is relatively prime to each $(s- s') \alpha$ for $s,s'$ varying over $\RR(W)$ with $s \neq s'$;
\item{(iv)} $v$ is relatively prime to $(b - s)\alpha$ for each $s$ in $\RR(W)$;
\item{(v)} $(v +1)s $ is not in $\RR(W)$ for all $s $ therein.
\vskip 6pt 

Let $d(i)$ denote the degree of $P_i$ (in $x$ of course). Define the polynomials
$$
Q_i  (x) = \frac 1{(v+1)^{d(i)}} P_i\( \frac x{v+1}\),
$$
so that $Q_i$ is a monic polynomial in $K[x]$, and its roots are  of the form $(v+1) s$ as $s$ varies over $S_i$ (and the multiplicities are preserved as well). 

Now set $a = (x-b) \prod_i Q_i$ in $K[x]$. Its zeros are precisely $b$ and $\brcs{(v+1)s}_{s \in \RR(W)}$, and by construction, none of these are in $\RR(W)$.

Now we show that $a K'[x, W^{-1} ] \cap V(K'[x, W^{-1}])= \brcs{0}$, which implies the corresponding result with $K$ replacing $K'$. If the intersection contained something other than zero, we can evaluate it at any of the zeros of $a$. At $x \mapsto b$ and $x \mapsto (v+1)s$, $s \in \RR(W)$ (taking each $s$ only once, regardless of its multiplicity), we obtain, as in the previous cases, 
$$\eqalign{
\prod P_i(b)^{A(i)} &= \frac 1{\gamma} \cr
\prod P_i ((v+1)s)^{A(i)} & = \frac 1 {\gamma}, \quad \text{and so} \cr
\prod_i \(\frac{P_i ((v+1)s)}{P_i(b)}\)^{A(i)}  &= 1 \quad \text{for all $s$ in $\RR(W)$.} \cr
}$$

Select $s$ in $\RR(W)$, and suppose it belongs to $S_j$ ($s$ is a root of $P(j)$). Then 
$$
\frac{P_j ((v+1)s)}{P_j(b)} = \(\frac{vs}{b-s}\)^{m(s)} \cdot \prod_{t \in S_j \setminus \brcs{s}} \( \frac{vs + (s-t}{b-t}\)^{m(t)}, 
$$
while for 
$i \neq j$, 
$$
\frac{P_i ((v+1)s)}{P_i(b)} = \prod_{z \in S_i}\( \frac{vs + (s-z)}{b-z}\)^{m(z)}.
$$
Multiplying numerators and denominators by the same power of $\alpha$ (enough to make them integral). Then we see that all the terms in the numerators and denominators except the power of $v\alpha s$ are relatively prime to $v$ (one justification for this is to localize 
$\Z_Q$ at $v$, creating a local ring; its integral closure in $K'$ is semilocal, and of course all of its maximal ideals contain $v$, so that all the things relatively prime to $v$ are units). This forces $A(j) =0$ (as in the preceding argument, where we rewrite the equation as an equality of polynomials). Since we can pick $s$ in any $S_i$, this forces $A(i) = 0$ for all $i$, a contradiction. 

If $K$ algebraic but infinite dimensional over $Q$, define the subfield $K'$ generated by the coefficients of all the irreducible factors on $P$, and then observe that the equations leading to a contradiction lead to a contradiction, and the conclusion is still valid in any field larger than $K
'$ (providing $a \in K'[x]$). 

This concludes the proof of Proposition \Cthr. 
\qed 

It is plausible that even if the field $K$ is finite (or algebraic over a finite field), then $K[x,P^{-1}]$ fails to satisfy \gui 2{} (where $P$ is any nonconstant polynomial). This suggests that the same negative result is true for any $\Z[1/n]$ ($n$ is an arbitrary positive integer). 

\SecT Acknowledgment

I would like to thank Ottmar Loos for informative discussions and
providing  additional references. Pace Nielsen also suggested additional references.

\SecT References

\item{[B]} H Bass, Algebraic K-theory, WA Benjamin Inc (1968) New York.

\item{[C]} PM Cohn, {\it On the structure of the $\gl (2,R)$ of a ring,} Pub math de l'IHES 30 (1966), 50--53.

\item{[GH]} KR Goodearl \& DE Handelman, {\it Metric completions of
partially ordered abelian groups,} Indiana U Math J 29 (1980) 861--895.

\item{[H]} DE Handelman, {\it Fixed points of two-sided fractional matrix
transformations,} Fixed point theory and applications (2007) doi:10.1155/2007/41930

\item{[HL]} JD Hews \& L Livshits, {Groups of matrices that act monopotently,} Electronic J of Linear Algebra 32 (2017)  423--37 doi: https://doi.org/10.13001/1081-3810.3479

\item{[HS]} P de la Harpe \& G Skandalis, {\it D\'eterminant associ\'e a une
trace sur une algebre de Banach,} Ann Inst Fourier (Grenoble) 34 (1984)
241--260

\item{[Ka]} 
Irving Kaplansky, {\it Elementary divisors and modules,} Trans Amer Math Soc 66 (1949) 464--491.

\item{[KKN]} 
Anjana Khurana{, Dinesh Khurana, \& Pace P Nielson}, {\it Sums of units in self-injective rings,} J Alg App 13(6):1450020 (2014) 7 pages.

\item{[Ko]} M  Koecher, {\it Uber eine Gruppe von rationalen Abbildungen,}
Invent\. Math\. 3 (1967), 136--171.

\item{[K]} M Kolster, {\it General symbols and presentations of general
linear groups,} Crelle's Journal 353 (1984) 132--164.

\item{[L]} O  Loos, {\it Steinberg groups and simplicity of elementary
groups defined by Jordan pairs,} J\,Algebra 186 (1996) 207--234.

\item{[L2]} ---\!\!---\!\!---\!\!--, {\it Elementary groups and
stability for Jordan pairs,} K-Theory 9 (1995) 77--116.

\item{[MM]} P Menal \& J Moncasi, {\it K$_1$ of von Neumann regular rings,} J Pure \& App Algebra 33 (1984) 295--312.

\item{[MV]} P Menal \& LN Vaserstein, {\it On the structure of GL over stable range one rings,}
J Pure App Alg 64 (1990) 149--62. 

\item{[SS]} Ferroz Siddique \& Ashish K Srivastava, {\it Decomposing elements of a right self-injective ring},
J Alg Appl 12 (2013). 

\item{[S]} Wolfgang Schmidt, Diophantine approximations and Diophantine equations,
Springer-Verlag Lecture Notes in Mathematics 1467 (1991). 

\item{[V]} P V{\aac}mos, {\it 2-good rings}, Quart J Math 56 (2005) 417--30. Behind a paywall---shame on you, Quarterly Journal of Mathematics!

\item{[WR]} Yao Wang \& Yanli Ren, {\it 2-good rings and their extensions},
Bul Korean Math Soc 50 (2013) 1711--1722. 

\item{[W]} JHM Wedderburn, {\it On continued fractions in noncommutative quantities},
Annals of Math (1913--14) 15 101--5. 

\vskip 6pt \noindent  rochelle2\@sympatico.ca

\end 
Introduction
Can be regarded as an enlargement of a comment of Kolster
or as development of noncommutative continued fractions,
or simply as an interesting nonstable invariant.

introduce AG(R), AG_i (R) (i=1,2) for all rings, based on  
concrete invariants defined for some rings, G(R), G_i(R).
main results If 6 invertible, then AG_2 (R) is perfect and
equals the commutator subgroup of AG(R).

Results on lengths of elements (in terms of the natural set of
generators) in AG(R), expressed in terms of the noncommutative
continued fraction polynomials...if 1 is in the stable range,
then every element has order 21/2 or less, etc, 

codntion xx

and also, if
1 in stable range and some modest additional conditions (usually 
satisfied) hold, then AG_2 (R) is simple and is the minimal normal
subgroup of AG(R).

The quotient AG(R)/AG_2 (R) is natually isomorphic to GL(R)/Z(R)^2\cdot
GL(R)' if xx holds (generally, the former is only a homomorphic image of
the latter). When G(R) is defined, an alternative proof is given, and
this also applies to some dense subrings of $R$.

assignment is not functorial (centre gets in the way), in appendix, we
discuss some related invariants for limits of finite dimensional
semisimple algebras and
their C* completions.

\comment
ThE PROOF IS WRONG:

\Lem Proposition. Suppose that $S$ is an elementary divisor ring satisfying \gui 2. Then for all $n$, $R = \Mn_n S$ satisfies \gui 3.
\Pf By the usual reduction result, we can assume $B = \diag (b_1, b_2, \dots, b_n)$ is diagonal; $C$ is an arbitrary element of $R$. 

First we examine the set of invertibles $U$ in $R$  \st $U+B$ is invertible. Included are the following
$$
\( \matrix -b_1 & u_1 & 0 & \dots & 0 \\
0 & * & u_2 & 0  & \dots & 0 \\
0 & * & * & u_3 & * \dots & 0 \\
 & \ddots & \ddots  & \ddots &  \\
 u_n & * & * & \dots & * \\
\endmatrix \), \tag 1
$$
where the $u_i$ appearing in the $(i,i + 1)$  (identifying $n+1$ with $1$) entries are units of $S$, and $*$ is completely arbitrary. To see that $U$ is invertible is straightforward: elementary row operations can be performed to show the columns span $F^{n\times 1}$ and column operations show that the rows span $F^{1 \times n}$. 

When we form  $F + B$, the $(1,1)$ becomes zero and again it is easy to check invertibility. 

Now with $C$, we first  replace the $*$ terms with the negative of the corresponding $c_{ij}$. The resulting matrix has bottom row $(u_n + c_{1,n}, 0, \dots 0)$. For positions in the cyclic permutation, that is, $(i,i+1)$, we may find, for each $i$ a unit $u_i$ \st $u_i + c_{i,i+1}$ is also invertible (this of course uses \gui {2} $n$ times). 

The resulting $U+C$ has bottom row a $(u',0, \dots , 0)$ more or less arbitrary stuff in the rest of the first column, and all the other entries on or below the main diagonal are zero. The entries corresponding to the permutation are units, and once again, it is straightforward to check that this is invertible (after multiplication on either side by the relevant permuation matrics, the result is triangular, with the diagonal entries all units). 
\qed

The proposition does not apply to $\Z$, since it does not satisfy \gui 2. However, it does apply if $S = \M_m \Z$ for any $m > 1$, since the elementary divisor ring property is preserved by taking matrix rings, and $\M_m S$ satisfies \gui 2. 
The upshot is the following weird result. 

THis is wrong too
\Lem Corollary. Let $S$ be an elementary divisor ring. If $n$ is composite, then $\Mn_n S$ satisfies \gui 3. 

We can reduce the set of possibly exceptional $n$ to four, $n = 3, 5, 7, 11$ (note that $\M_2 \Z$ does not satisfy \gui {3}---since $\M_2 F_2$ does not). 

This is replaced by a better result. 
\Lem Lemma. Let $S$ be an elementary divisor ring \st $\Mn_n S$ and $\M_m S$ both satisfy \gui 3. Then $\M_{n+m} S$ satisfies \gui 3. 

\Pf We can assume that $n < m$ (the case that $n= m$ is covered by the proposition).  Take $B,C \in \M_{n+m}S$; as before, we can replace $B$ by a diagonal matrix. Now parfition both matrices in four pieces corresponding to $m,n$. The diagonal matrix $B$ is a direct sum of two diagonal matrices, and $C$ consists of the usual  four pieces, 
$$
C = \(\matrix M & Y \\ X & N \endmatrix \),
$$
where $M, N$ are square of sizes $n,m$ respectively, $X$ is $n \times m$ and $Y$ is $m \times n$. 

If $V \in \gl(n,S)$ and $U \in \gl(m,S)$, then we can conjugate  $B$ and $C$ (this will not change the problem) by the block diagonal matrix $W:= V \oplus U$; the resulting $B$ will no longer be diagonal, but will still be a direct sum of two matrices. We have 
$$
W C W^{-1} = \(\matrix VMV^{-1} & VYU^{-1} \\ UXV^{-1} & U N U^{-1} \\\endmatrix \),
$$
The matrix $X$ is $m \times n$ where $m > n$. An elementary property of elementary divisor rings is that if there exist invertible matrices $U \in \gl(m,S)$ and $V' \in \gl (n,S)$ \st all the entries below the first $n$ of $UXV'$ are zero (and the remaining $n \times n$ matrix is diagonal); set $V = (V')^{-1}$. That is, the first $n$ columns of the current matrix
$$
\(\matrix  VMV^{-1} \\ D \\ {\pmb 0} \\
\endmatrix\)
$$
where $D$ is a diagonal square matrix of size $n$, and the zero matrix is $(m-n) \times n$.

Write the current version of $B$ (that is, after conjugation by the block diagonal invertible) as $E \oplus F$, and let $G$ denote the lower $m \times m$ block of the conjugated $C$. 

Applying the \gui {3}, there exists $R\in \gl(n,S)$ \st both $E + R$ and $VMV^{-1} + R$ are invertible; there exists $T \in \gl(m,S)$ \st both $F+ T$ and $G + T$ are invertible. 

Now set $Z$ to be  of $R \oplus S$ together with $-D$ in the same block position as $D$. Then $Z$ is lower block triangular, as is $(E \oplus F) + Z$, both with invertible block entries along the (block) diagonal. It easily follows that $Z$ and $(E \oplus F) + Z$ are invertible. Finally, $WCW^{-1} + Z$ is block {\it upper\/} with invertible diagonal blocks, and thus is also invertible. 
\qed

This is wrong:
\Lem Corollary. Let $S$ be an elementary divisor ring. Then for all $n$ other than $n=2,3,5,7, 11$, $\Mn_n S$ satisfies \geq 3. 

\Pf By xxx, for $\Mn_n S$ satisfies \gui {3} for composite $n$. Every prime $13$  or more is expressible as a sum of two composites: $p - 9$  is even and exceeds two, so is composite. 
\qed

For $S=\Z$, $\M_2 \Z$ does not satisfy \gui {3}. But the status of each of $\Mn_n \Z$ for $n \in \brcs{3,5,7, 11}$ is (currently) unknown.

\endcomment

Kaplansky [Ka, p 465] defines a ring $S$ to be  (right) {\it Hermite\/} if for every $a,b$ in $S$, there exists $Q$ in $\gl (2,S)$ and $d$ in $S$ \st $(\matrix a \\ b \\\endmatrix)Q =(\matrix d \\ 0 \\\endmatrix) $; left Hermite is similarly defined with the matrix on the left, and the row replaced by a column. Right Hermite implies (since $aR + bR = dR$) all finitely generated right ideals are principal. 

We say a matrix (usually, but not necessarily square) over a ring $S$  is {\it \paren{upper, lower respectively} block triangular\/} if it has the form
$$
\(\matrix M  & {\pmb 0} \\ Y & N\) \qquad \(\matrix M  & Y\\ {\pmb 0}  & N\)
$$
where $M$ and $N$ are square, and the remaining matrices are rectangular to fit (that is, if the whole matrix is $K \times K$, then $M$ is $m  \times m$, $N$ is $n \times n$ with $K = m+n$, the lower left matrix is $n \times m$ and and the upper right matrix is $m \times n$. We implicitly do not permit $m $ or $n$ to be zero. 

\comment
wrong proof
\Pf Since $S$ is an elementary divisor ring, we may assume $B = \diag (b_1, b_2, \dots , b_n$ and $C = (c_{ij})$ is arbitrary. For each position $(i,i-1)$ (including $(1,n)$), there exists (by \gui 2) an invertible $u_{i,i-1}$ \st $u_{i,i-1} + c_{i,i-1}:= u'_{i,i-1}$ is also invertible. 

Define $U $ via 
$$
U_{ij} = \cases -b_{i} & \text{if $i = j$ and $i < n$}\\
 -c_{n,n}& \text{if $i = j = n$}\\
 -c_{i,n} & \text{if $j = n$ and $i < n$}\\
-c_{i,j}& \text{if $i \geq 2$ and $j leq i-1$}\\ 
0 & \text{if $i <  j < n$}\\
\endcases
$$

Then we have 
$$
U = \(\matrix -b_1 & 0 & 0 & \dots &&0 & u_{1,n} \\
u_{2,1} & -b_2 & 0 & \dots &&0 & -c_{2,n} \\
-c_{3,1} & u_{3,2} & -b_3 & 0&\dots &0 & -c_{3,n} \\
\ddots & \ddots & \ddots & 0&0 & & \vdots \\
-c_{n-1,1} & -c_{n-1,2} & \dots & -c_{n-1,n-3}& u_{n-1,n-2}& -b_{n-1} & - c_{n-1,n} \\
-c_{n,1} & -c_{n,2} & &\dots & -c_{n-1,n-2}& u_{n,n-1}& - c_{n,n} \\
\endmatrix\)
\qquad U+B = blah \qquad U + C = blah,blah
$$

In the case that $n=3$, we have
$$
U = \(\matrix
-b_1 & 0 & u_{1,3} \\
u_{2,1} & -b_2 & -c_{2,3}\\
-c_{3,1} & u_{3,2} & -c_{3,3} \\
\endmatrix
$$
AArgh

\endcomment

\Lem Proposition. Suppose that $S$ is a very weakly Hermite ring satisfying \gui {3}. Then for all $n$, $\Mn_n S$ satisfies \gui {3}. 

\Pf Pick $B,C$ in $\Mn_n S$; we may replace $B$ and $C$ by $UBV$ and $UCV$ respectively, so by the very weak Hermite property, we can assume $B$ is upper triangular (the argument for lower triangular is almost identical). 

We use the standard subscript notation for matrices, for example, $B = (b_{ij}) $ (or $b_{i,j}$ if there is potential ambiguity). 

For each $i$, there exists invertible $u_{ii}$ in $\gl(1,S)$ \st both $u_{ii} + b_{ii}$ and $u_{ii} + c_{ii}$ are invertible. 

Define the upper  triangular matrix $W$ via 
$$
W_{ij} = \cases u_{ii} & \text{if $i = j$}\\
-c_{ij} & \text{if $ i < j$} \\
0 & \text{if $i > j$}
\endcases
$$
A triangular matrix with invertibles along the main diagonal is (two-sided) invertible, as is easy to see. Thus $W$ is invertible. As $W + B$ is upper triangular with diagonal entries $u_{ii} + b_{ii}$, it is invertible. Finally, $W + C$ is lower triangular with diagonal entries $u_{ii} + c_{ii}$, so is also invertible. 
\qed

It is straightforward to show that if the diagonal blocks, in this case, $M$ and $N$, are invertible (as always, {\it invertible\/} means two-sided invertible), then the matrix is invertible. The converse does {\it not\/} hold; that is, over some rings, there are invertible block upper triangular rings whose diagonal blocks are not invertible (and the inverses are not block triangular either). 

\Lem Lemma. Suppose that $S$ is a ring \st $M_n S$ and $\M_m S$ satisfy \gui 3, and set $k = m+n$.  Let $B \in \M_k S$ be block triangular with blocks of size $m$ and $n$, and let $C$ an arbitrary element of $\M_k S$. Then there exists $U \in \gl (k,S)$ \st $U + B$ and $U + C$ are both invertible. 

\Pf Assume $B$ is lower triangular (the proof for upper triangular is essentially the same, but note that we cannot just take the transpose, since in general, the transpose of an invertible matrix need not be invertible!), say 
$$
B = \(\matrix M  & Y\\ {\pmb 0}  & N \\ \endmatrix \),
$$
where $M$ is $m \times m$, $N$ is $n\times n$, etc
and partition $C$ accordingly, 
$$
C =  \(\matrix P  & Z\\ W  & Q\\ \endmatrix\).
$$
That is, $P$ is $m \times m$, $Q$ is $n \times n$. 
By hypothesis, there exist $U \in \gl(m,R)$ and $V \in \gl(n,R)$ such that $U + M$ and $U +P$ are both invertible ($m\times m$ matrices), and $V + N$ and $V+Q$ are both  (invertible $n\times n$ matrices). 

Now form the block lower triangular matrix,
$$
E = \(\matrix U  & -Z\\ 0  & V \\ \endmatrix \)
$$
Then $E$ is invertible, $E + B$ is lower triangular with invertible diagonal blocks, so it invertible, and $E + C$ is upper triangular with invertible diagonal blocks, so is also invertible. 
\qed

\Lem Corollary. Let $S$ be a very weakly Hermite ring  and suppose that for some integers $m$ and $n$, $\M_m S$ and $\Mn_n S$ both satisfy \gui 3. Then so do $\M_r S$  for all $r$ in the additive subgroup of $\N$ generated by $\brcs{m,n}$. 

\Pf Giiven $B\in \M_{m+n} S$, there exist invertible matrices $G, H$ \st $GBH$ is in upper block diagonal form. From the preceding, $\M_{m+n} S$ satisfies \gui 3, and induction yields the result. 
\qed

As a particular result, if $S$ is an Hermite \plainfootnote1{Charles Hermite, after whom this property is named, was French, so the {\it H\/} is not pronounced; thus {\it an,} not {\it a.} Kaplansky  did this.}%
satisfying \gui 3, then so do all of its matrix rings. This provides yet another proof that if $q \geq 4$, then for all $n$, $\Mn_n F_q$ satisfies \gui 3.

Division rings satisfy everything in sight, specifically, they are Hermite rings (and much more). Moreover, if $D$ is finite, then it is commutative, and thus of the form $F_q$. If $D$ merely has an infinite centre, then $\Mn_n D$ satisfies  \gui {k} for all $k$, by a trivial argument. This leaves the case that $D$ is a noncommutative (hence infinite) division ring with finite centre.

\Lem Corollary. Let $S$ be a right or left artinian ring. Then $S$ satisfies \gui {3} if and only if it has no images isomorphic to any one of $F_2$, $\M_2 F_2$, or $F_3$.

\Pf Necessity of the conditions is clear. We reduce to the semisimple case by factoring out the Jacobson radical. Since products behave well, we need only deal with $\Mn_n D$ where $D$ is a division ring, and as mentioned above, we can assume $D$ is infinite. 

Obviously, $D$ is an elementary divisor ring, and since it is infinite, it satisfies \gui {k} for all $k$. By Proposition xxx (above) all matrix rings satisfy \gui 3. 
\qed

Appendix 
\item{[Ko]} M  Koecher, {\it Uber eine Gruppe von rationalen Abbildungen,}
Invent\. Math\. 3 (1967), 136--171.
\comment
direct limit examples

The assignments $R \mapsto \Ag(R)$ and its relatives are not functorial,
because the centre can fail to be mapped to the centre under ring
homomorphisms. This means that we cannot simply write down the values of
these invariants on direct limits, say of finite-dimensional semisimple
algebras (over a field). We present some  interesting direct limit
examples. Let
$R_m =
\Mn {3^m} \C $, and define unital homomorphisms $\Arrow \phi_{m,m+1}; R_m
. R_{m+1}$ via $a \mapsto a \oplus a\oplus \overline{a}$; let $R = \lim
R_m$. This is a very well known object; it is a simple real algebra (the
maps are not complex algebra homomorphisms), with centre the real
scalars. We observe  for $a$ in $R_m$,  that $\det a/|\det a| = \det
(\phi_{m,m+1} (a))/|\det (\phi_{m,m+1} (a))|$, where the determinants are
the unnormalized ones in the relevant matrix rings. This induces a group
homomorphism $\Arrow D; \GL R . \Tt$ (the unit circle), simply by sending
$a$ in $R_m$ to $\det a/|\det a|$ (the $m$ has only to be chosen large
enough that the representative of the element of $\GL R$ is itself
invertible). (The function  $D$ is {\it not\/} continuous \wrt the norm
topology on $R$.) If
$D (r) = 1$, that is,
$r$ is represented by
$a$ in
$R_m$
\st $\lambda:= \det a$ is positive, then we can write $a =
\lambda^{1/3^{m}}a_0$, where $\det a_0 = 1$, hence $a_0$ is a commutator,
and thus  $r$ belongs to $\Z(R)^2 \cdot [\GL R, \GL R]$ (the square comes
from the fact that $\lambda$ and its roots are positive).

Recall the definition of $\Arrow \Psi_R; \GL R/Z(R)^2 \GL {R}'.
\Ag(R)/\Ag_2 (R)$ from section 8. The following property for a (unital)
ring 
$R$ is sufficient to show that $\Psi_R$ is an isomorphism.

\item{($\Ss$)} for all $k =1,2, \dots$, for all $r_1$, $r_2$, \dots,
$r_{2k}$ in $R$
\st $Q_{2k} \equiv Q_{2k}(r_1, \dots , r_{2k})$ is invertible, the element
$Q_{2k}^{-1} Q_{2k}\Op$ is in the commutator subgroup of $\GL R$.

We notice that if $S = \Mn N F$ for any field at all, then $S$ satisfies
($\Ss$); the property is preserved by finite direct products (and infinite
direct products if there is a bound on the number of commutators required
in each of the components). It is obviously preserved by direct limits
(in contrast, the condition that $\Psi_R$ be an isomorphism is not a
priori preserved by direct limits, which motivates the definition here).
It is preserved by quotients if
$S$ has stable rank 1, or more generally, if invertibles lift modulo
two-sided ideals.

Since $R$ is a limit of matrix algebras over a field and each of these
has property ($\Ss$), so does $R$. Hence $\Psi_R$ is an isomorphism, and
it is compatible with $D$.

It is worthwhile investigating other invariants for this and a class of
complex algebras. In this case, $\GL R /[\GL R , \GL R] \iso \C^*$,
$\text{K}_0 (R) = \Z [\frac 13]$,
$Z(R)^2 = \R^+$, and
$\G(R)/\G_2 (R) \iso \C^*/\R^+ \iso \Tt$. (However, $\Tt$ should
be regarded only as an abstract abelian group, not as a
topological group.) Moreover,
$-1$ is not in
$Z(R)^2 \cdot [\GL R,\GL R]$, so the image of $-1$ under $\phi_R$ is not
trivial, and thus $\G_2 (R) \neq \G_1 (R)$, although $\G (R)/\G_1 (R) \iso
\Tt/\brcs{\pm1} \iso \Tt$ again.

With complex algebra homomorphisms, the invariants are somewhat
different. Suppose this time that $R = \lim R_m$ where each $R_m$ is a
product of $s(m) $ matrix algebras over the complexes, and the maps
$\Arrow \phi_{m,m+1};R_m. R_{m+1}$ are unital complex algebra
homomorphisms. Each $R_m$ satisfies ($\Ss$), and thus $\G(R)/\G_2 (R)$ is
naturally isomorphic to $\GL R/(Z(R)^2 \cdot [\GL R,\GL R])$. If we
assume that $R$ is simple, then $Z(R) = Z (R)^2 = \C^*$. Also, $\GL
R^{\text{ab}}$ commutes with direct limits, so that $\GL R^{\text{ab}}$
is the limit $\lim (C^*)^{s(m)} \to (C^*)^{s(m+1)}$, the maps induced by
the maps between the $\text{K}_0$ groups. The copy of $\C^*$ coming from
$Z(R)^2$ is just $\Delta (\C^*)$, the diagonally embedded copy, and we
obtain   $\G(R)/\G_2 (R) \iso \lim (C^*)^{s(m)}/\Delta_{s(m)}
(\C^*)
\to (C^*)^{s(m+1)}/\Delta_{s(m+1)} (\C^*)$ (of course, the groups are
regarded as multiplicative groups).

Now we consider $\GL R/\overline {[\GL R,\GL R]}$, where the
overline denotes the {\it relative\/} closure (that is, within $R$, not
of course within the completion of $R$). It turns out under modest
assumptions, that as an abelian group, this is the real span of the image
of
$\text{K}_0 (R)$ in its affine representation (formulated this way, it
does have an ordering, which corresponds to orientations of loops; traces
on
$R$ correspond to determinants on $\GL R$). 

For an algebra written as a unital direct limit of unital semisimple
finite dimensional algebras, we adopt the following notation. The ring
$R  =
\lim R_m$, where $R_m = \Mn {a(m,1)} \C \times \Mn {a(m,2)} \C \times
\dots \times \Mn {a(m,s(m))} \C$, i.e., $R_m$ is a product of $s(m)$
matrix algebras of sizes $\brcs{a(m,i)}_{i=1}^{s(m)}$. Denote the maps
$R_m \to R_{p}$ (where $p > m$) $\phi_{m,p}$ and the maps $R_m \to R$ by
$\phi_m$. Finally, let $a(m) = \inf_{i \leq s(m)} \brcs{a(m,i)}$, the
minimum of the matrix sizes. It is known that the following are
equivalent:

\item{(a)} $a(m) \to \infty$
\item{(b)}  $R$ has no nonzero finite dimensional images
\item{(c)} the natural map $\text{K}_0 (R) \to \text{Aff} ( T(R))$ has
dense range.

\noindent These properties (for the corresponding AF algebras; there is no
harm in extending the definition) are known as {\it approximate
divisibility.} If
$R$ is simple but not finite dimensional, then $R$ is approximately
divisible. The only implication that is not obvious is (b) $ \implies$
(c); this was shown in [GH, Theorem 4.8].

\Lem Proposition \Aone. Let $R = \lim R_m$ be approximately divisible. If
$r = (r(1),r(2), \dots, r(s(m)))$ is an element of $R_m$ \st $|\det\,
r(i)| = 1$ for all $i$, then $\phi_m (r)$ belongs to $\overline{[\GL R,
\GL R]}$.

\Pf First, for each $p > m$, the condition on the determinants of the
components (that they all have absolute value one) is 
satisfied by $\phi_{m,p} (r)$, the image of $r$ in $R_p$. For each $p$,
let $z(p)$ be the $s(p)$-tuple scalar matrices in $R_p$, $z(p) =
(z(p,1), z(p,2), \dots, z(p,s(p)))$  where $z(p,j) $ is the determinant
of the $j$th component of $\phi_{m,p} (r)$. Define complex numbers
$w_{p,j}$with the following properties:

\item{(i)} $w_{p,j}^{a(p,j)} = z(p,j)$ 
\item{(ii)} $0 \leq \Arg w_{p,j} \leq 2\pi/a(p,j)$

Of course, we can always arrange for the roots of $z(p,j)$ to satisfy
(ii) by taking the root closest  to $1$ with positive imaginary part.
Property (i) entails that if $w_p = (w_{p,j})$, then the determinant of
each component of $w_p^{-1}\phi_{m,p} (r)$ is $1$. A particular
consequence of (ii) is that $|w_{p,j} - 1| \leq 2\pi/a(p,j)$ (here
viewing the $w_{p,j}$ as complex numbers rather than scalar matrices, but
it makes no difference), and thus $\| w_p - 1 \| \leq \max_j
\brcs{2\pi/a(p,j)}$ (here viewing $w_p$ and $1$ as elements of $R_p$), and
by assumption,
$\max_j \brcs{1/a(p,j)}= 1/a (p)  \to 0$ as $p\to \infty$. Hence the
sequence of elements of $R$ given as $\brcs{\phi_p (w_p)}$ converges to
$1$ (now viewed as an element of $R$). The $w_p$ are of course  unitary,
so
$\brcs{\phi_p (w_p)^{-1}} \to 1$ as well.

Since the determinant of each component of $r_p:= w_p^{-1}\phi_{m,p} (r)$
is
$1$, it follows that $r_p$ is a commutator in $R_p$, hence $\phi_p (r_p)$
is a commutator in $\GL R$. However, $\phi_p (r_p) = \phi_p (w_p)^{-1}
\phi_m (r)$, so as $p \to \infty$, these converge to $\phi_m (r)$. Hence
$\phi_m (r)$ is a limit of a sequence of commutators of $\GL R$.
\qed

For the rest of this section, we assume the reader is familiar with
traces, determinants, ordered $\text{K}_0$ of direct limits of finite
dimensional semisimple algebras, AF C*-algebras, etc. The arguments are in
outline form. De la Harpe and Skandalis [HS] have calculated a number of
invariants associated to AF algebras, in particular, if $S$ is an
approximately divisible AF algebra, then $\GL S/\overline{[\GL S, \GL
S]}$  is naturally isomorphic to $\Aff T(R)$. The direct limits
considered here are dense subalgebras, and it is rather surprising (to
me) that  something similar holds, namely, if $R$ is a limit of complex
 finite dimensional semisimple algebras (not requiring complex algebra
homomorphisms between them), then $\GL R/\overline{[\GL R, \GL R]}$ is
naturally isomorphic to the real span of the image of $\text{K}_0 (R)$ in
$\Aff T(R)$ (keeping in mind that the closure, $\overline{[\GL R, \GL
R]}$, is the {\it relative\/} closure, calculated within $\GL R$).

There is a family of little-studied maps $\text{K}_0 (R) \to
\text{K}_1(R)$ (for any ring $R$ with $2$ invertible), of which the
simplest is induced by the assignment $p \mapsto 1+p$ where $p$ runs over
all idempotents in all matrix rings of $R$. We now fix $R$ to be our
standard direct limit. The assignment $p
\mapsto 2^p$ for all projections in $R$---obviously $2^p$ is in $R$ 
(this map is slightly different from $p \mapsto 1+p$;  it is defined only
on projections in
$R$, not idempotents in all matrix rings; these differences are not
significant)---is a little easier to work with. This extends uniquely to a well defined group homomorphism
$\text{K}_0 (R) \to
\text{K}_1(R)$, and since we are interested in $\GL R /\overline{[\GL R,
\GL R] }$ (which is a quotient of $\text{K}_1 (R)$, as $R$ has stable rank
one), it is natural to look at the kernel of the composite map $\text{K}_0
(R) \to\text{K}_1(R) \to \GL R /\overline{[\GL R,
\GL R] }$. The K-theory of this paragraph
could have been dispensed with; however, it motivates the following.

\Lem Lemma \Atwo.
Suppose $R = \lim R_m$ is an approximately divisible limit of finite
dimensional semisimple $\C$-algebras. Let $p$ and $q$ be projections in
$R$ \st for all normalized traces $\tau$ of $R$, $\tau (p) = \tau (q)$.
Then
$2^{p-q}$ is a limit of single commutators in $R$, hence is in the
relative closure of
$\overline{[\GL R,\GL R]}$.

\Pf (outline) Without loss of generality, we may assume $p$ and $q$
commute with each other, and are represented by $p_0$ and $q_0$ 
respectively, in $R_m$. The tracial condition means that for all
$\epsilon > 0$, there exist projections $p_{\epsilon} < p$ and
$q_{\epsilon} < q$ in $R$ \st $\tau (p-p_{\epsilon}), \tau
(q-q_{\epsilon}) < \epsilon$ for all normalized traces $\tau$, and
$p_{\epsilon}$ is conjugate to $q_{\epsilon}$. 

This conjugacy is implementable in some $R_n$ (where $n \equiv n(\epsilon)
\geq m$). We deduce that for all $\epsilon > 0$, there exists an integer
$l \equiv l(\epsilon)$ \st with $v_l:= \phi_{m,l} (p) - \phi_{m,l}(q)$,
every component of $v_l = (v_{l,i} )$ in $R_l$ satisfies $|\Tr_i
v_{l_i}|/a(l,i) \leq \epsilon$, where $\Tr_i$ is the standard unnormalized
trace on the $i$th matrix algebra. (This is a fairly standard computation
involving limits.) It follows that the unnormalized (that is,
the usual) determinants of the components satisfy $\det_i 2^{v_{l,i}} <
2^{a(l,i)\epsilon}$ (the determinant of an exponential of a difference of
idempotents is always positive, so there is no need for absolute value
signs.

Define $d_i$ to be the positive $a(l,i)$ root of $\det_i 2^{v_{l,i}}$,
so that if $d_{(l)} := (d_i)$ is the corresponding central element of
$R_l$, then $2^{v_l}d_{(l)}^{-1}$ has determinant one at each of its
components, and is thus a commutator. On the other hand, $\log_2 d_i \leq
\epsilon a(l,i)/a(l,i) = \epsilon$, so that $\| d_{(l)} - 1\| < \epsilon$
in $R_l$. Hence $\brcs{\phi_l (d_{(l)})} \to 1$ (as $\epsilon \to 0$), so
$\brcs{\phi_l (d_{(l)})^{-1}} \to 1$. Thus $2^{p-q}$  is the limit in $R$
of $\brcs{\phi_l (2^{v_l}) \phi_l(d_{(l)})^{-1}}$, hence is a limit of
(single) commutators. 
\qed

 It follows readily that if the direct limit is approximately divisible,
then the kernel of the map $\text{K}_0 \to \text{K}_1$ is exactly the set
of infinitesimals (those elements sent to zero by every trace).
\endcomment
\comment

\Lem Proposition . Suppose $R$ is a simple unital ring generated by its
invertibles. Let $g$ be an element of $\Ag(R) $ with $\ord  (g) < 4$.
{\par}
\item{(i)} 
$\Ag_2 (R) \subseteq \Nn (g)$.{\par}
\item{(ii)} If $R$ is generated by the commutator subgroup of $\GL R$ and
the centre, and $g$ belongs to $\Ag_2 (R)$, then $g$ is contained in no
proper normal subgroups of $\Ag_2 (R)$.
\endcomment

\SecT 6 Normal subgroups

We show that $\Ag_2 (R)$ is the commutator  subgroup of $\Ag (R)$ and
is  perfect, under very modest conditions. Our first step is
obtain bounds (somewhat sloppy) on  the number of $e_{a}$ terms are
necessary to express suitable $m_{r,s}$ as a product thereof; the bound
itself is not that important, but the result does yield that a lot of
elements of this form ($m_{r,s}$) are in $\Ag_2 (R)$. The quotient group,
$\Ag (R)/\Ag_2 (R)$, is also quite interesting and we will investigate
that as well. 

\Lem Proposition \sixone. Suppose that $r$ and $s$ are invertible
elements of $R$, and there exist (multiplicative) commutators
$u_1$, $u_2$, \dots, $u_k$ together with a central invertible
element $\lambda$ \st $s^{-1}r =
\lambda^2 u_1
\cdots u_k$. Then there exist $\brcs{a(1), a(2), \dots, a(8k+4)}$ \st
$m_{r,s} = e_{a(1)} e_{a(2)}\cdots e_{a(8k +4)}$. In particular,
$m_{r,s}$ belongs to $\Ag_2 (R)$. 

\Pf First, let $u = x^{-1}y^{-1} xy$ be a commutator of elements of $R$.
Set
$a =  x^{-1}y^{-1}x$ and observe that $m_{a,a} m_{ y,{x^{-1}yx}} = m_{ 
{u},{x^{-1}yxa}} = m_{{u},1}$. Since $a$ is invertible, $m_{a,a}$ is  
of  the form $e_{\bullet} e_{\bullet}e_{a-1}e_{1}$. 

  Next we observe that everything of the form $m_{y,{x^{-1}yx}} $ can be 
expressed as a product of four $e_a$s---set $a_1 = x$ and $a_2 =
x^{-1}(y^{-1}-1)$.   Then  $a_2 a_1 =  x^{-1}y^{-1}x -1$ and $a_1 a_2 = 
y^{-1}-1$. Thus $y =  (1+a_1 a_2)^{-1}$ and $x^{-1}yx = (1+a_2
a_1)^{-1}$. Since $m_{(1+a_1  a_2)^{-1}, (1+a_2 a_1)^{-1}}$ is
a product of four $e_a$s, this step is done. 

  Thus $m_{{u},1}$ is a product of eight $e_a$s for any 
commutator $u$. It follows immediately that if $v$ is a product of $k$ 
commutators, then $m_{v,1}$ can be expressed as a product of $8k$ $e_a$s. 
 Since we can replace the pair
$(r,s)$ by $(\lambda r,\lambda s)$ without changing $m_{r,s}$, we can
assume $\lambda =1$. Now $m_{r,s} = m_{rs^{-1},1} m_{s,s}$,  expresses
$m_{r,s}$ as a product of two expressions, one of which is a product of
$8k$ $e_a$s, and the other is a product of four $e_a$s. 
\qed

Consider the map $\Arrow \psi_R ; \gl (1,R) ; \petwo$ given by  $r \m r,1$. This is a group homomorphism and induces one
 $$
 \Arrow \Psi_R;\gl(1,R)/Z(R)^2 \gl(1,R)' . \pe /\petwo.
 $$ 
 (Recall that $Z(R)^2$ consists of the squares of the invertible central elements.) That $\Psi_R$ is well-defined follows from Proposition 
. 

More is true: $\Psi_R$ is onto. Since $\pe$ is generated by $\petwo$ and $\brcs{\m x,y}$, it suffices to show that every $\m x,y$ can be written in the form $\psi(r) e_{a(1)} \dots e_{a(2k}}$ for some choice of invertible $r$ and elements $a(i)$. However, $\m x,y = \m xy^{-1},1 \m y,y$. The latter factor is a product of four $e$s, and the former is $\psi (xy^{-1}$. 

At this point, we have that $\Psi$ is an onto group homomorphism. A particular consequence is that $\pe/\petwo$ is abelian. This will be improved (Proposition 
)

Maintaining the notation of the proof, define $T_v = F^n \setminus (\brcs{0} \cup \setminus S_v)$, so that $T_v$ consists of the nonzero elements of $G$, $d$, \st $d+v$ is  invertible. It follows that 
$$
\eqalign{
|T_v| & \geq q^n -1 -  \frac{q^n-1 - (q^{n-r}-1)}{q-1}
\cr
&= \(1 - \frac 1q\)(q^n-1) + \frac{q^{n-r}-1}{q-1}. \cr \tag 1
}$$
(Recall that $r$ is the rank of $v$.) Define $\alpha(q,n) = |T_v|/|G| = |T_v|/q^n$ (although not explicit in the notation, the cardinality of $T_v$ depends on the choice of $G$; if there is ambiguity about which $v$ is referred to, we use the notation $\alpha_v (q,n)$). 

It is  easy to verify the following from (1). 

\item{(a)} If $q \geq 3$ and $n \geq 2$, then $\alpha(q,n) > 1/2$.
\item{(b)} For $q = 2$ and $n \geq 2$
$$
\alpha(2,n) \geq \frac 12 - \frac 3{2^{n+1}} + \frac 1{2^r}.
$$
In particular, if $r = n$, $\alpha(2,n) \geq 1/2 - 2^{-(n+1}$; if $r = n-1$, \alpha(2,n) = 1/2 + 2^{-(n+2})$, and for all smaller $r$, $\alpha(2,n)  \geq 1/2 + 5/2^{n+1}$. 

\Lem Theorem xxx. With $F = F_q$, the ring $R = M_n F$ satisfies \gui {3} if and only if either $q > 3$, or $q =3$ and $n\geq 2$, or $q = 2$ and $n \geq 3$.

\Pf Necessity  is clear (we have already seen that  $M_2 F_2$, $F_2$, and $F_3$ fail \gui 3), and if $q \geq 4$ and $n=1$, that $F_q$ satisfies (4) is immediate.

So assume $n \geq 2$, and let $b,c$ be nonzero elements of $R$. Fix a choice of $G$ as in the statement of the Lemma. It is sufficient to show that the subset of $G$ given by  
$T_b \cap T_c$ is nonempty.

This will occur if $\alpha_b (q,n) + \alpha_c (q,n) > 1$. By part (a), above, this will hold when $q \geq 3$. So we are reduced to the case that $q = 2$, and thus $n \geq 3$.

We see from the computations in (b), that the only way $\alpha_b (q,n) + \alpha_c (q,n) $ can be less than or equal to $1$ is if $\rk b + \rk c \geq 2n-1$---that is, either both $b$ and $c$ is invertible, or one of them is invertible  and the other has rank $n-1$. 

So we reduce to the special case, that $\rk b = n$, that is, $b$ is invertible.

\comment
Let $\alpha = (\lambda \mu -1)^{-1}$, and observe 
$$\eqalign{
\(\matrix \lambda & 0 \\ 0 & \mu^{-1} \\ \endmatrix\) 
\(\matrix 1 & \alpha a\\ 0 & 1\\ \endmatrix\) 
\(\matrix \lambda^{-1} & 0 \\ 0 &
\mu \\ \endmatrix\)
\(\matrix 1 & -\alpha a \\ 0 & 1 \\ \endmatrix\) 
 & = \(\matrix 1 & \lambda \mu\alpha \\ 0 & 1 \\\endmatrix\) 
\(\matrix 1 & -\alpha a \\ 0 & 1 \\ \endmatrix\)  \cr
& = \(\matrix 1 & a \\ 0 & 1 \\ \endmatrix\) \cr
 }$$
(b) Set $\lambda = 2$, $\mu = 1$, and apply (b).
(c) Standard method working in size three and higher matrix rings.\qed
\endcomment

For each $v(i)$ that is rank one (if one exists in $S$), define $\Cal B(i) = \Set{\lambda \in F^*}{\lambda v(i)\in S}$.   Pick $v(i)$ to maximize $\left| \Cal B(i) \right| $; say the maximum value is $m$. Obviously, $ m \leq q-1$. 

\noindent (b) Suppose $m = q-1$. Then we may assume, after applying $P,Q$ as usual, that $v(1) = \(\smallmatrix 1 & 0 \\ 0 & 0 \\ \endsmallmatrix\)$, and $v(2)$ through $v(q-1)$ vary over $\lambda v(1)$ as $\lambda$ runs over $F^* \setminus \brcs{1}$. 

Now we manipulate $v(q)$ as we did above. We may multiply the second row or column by a nonzero scalar or add the second row or column to the first without changing the set $\brcs{\lambda v(1)}$, and we can also multiply the first row/column by a nonzero scalar. If $v(q)$ is rank $1$, then we reduce to one of $W_1 $ or $W_2$ above; if $v(q)$ is rank two, we reduce to one of $\I$ and $\(\smallmatrix 0 & 1 \\ 1 & 0 \\\endsmallmatrix\)$. For each of these four matrices, we write down a matrix $M$ \st $M+v(i)$ is invertible for all $i$ (where we have relabelled the last one---after transformation---$v(q)$. 

\noindent (c) Suppose $m \leq q-2$; then there are at least two disjoint pairs $(i,j)$ and $(i',j')$ \st each of $\brcs{ v(i), v(j)}$ and $\brcs{ v(i'), v(j')}$ is linearly independent (this obviously requires $q \geq 4$. 

Following up the earlier remark (replacing $\ker v$ by $\ker (d_0 + v)$, it is possible to verify that 
if $n \geq 2$ and $q \geq 3$,  and $v_1, v_2, \dots, v_q$ is a list of elements of $\Mn_n F_q$ with at least one of the $v_i$  invertible, then there exists invertible $d$ \st all $d+ v_i$ are invertible. That is, $\Mn _n F_q$ {\it almost\/} satisfies \gui {q+1}, and probably does for sure; for $n=2$, $\Mn _2 F_3$ {\it does\/} satisfy \gui 4.